\documentclass[twoside,11pt]{article}
\usepackage{amssymb,amsmath,amsthm,amsopn} 
\usepackage{amscd} 
\usepackage{psfig,latexsym,graphicx}
\usepackage{fullpage}

\newenvironment{Quote}
{\begin{list}{}{%
   \vspace{-.8cm}
   \setlength{\rightmargin}{7mm}
   \setlength{\leftmargin}{7mm}}
    \item[]}
{\end{list}}

\refstepcounter{equation}

\newtheorem{lem}{LEMMA}
\newtheorem{theo}[lem]{THEOREM}

\newtheorem{coro}[lem]{COROLLARY}
\newtheorem{prop}[lem]{PROPOSITION}

\theoremstyle{definition}
\newtheorem{definition}[lem]{\indent Definition}
\newtheorem{rem}[lem]{\indent Remark}
\newtheorem{conj}[lem]{\indent Conjecture}
\newtheorem{ex}[lem]{\indent Example}
\refstepcounter{lem}

\renewcommand{\descriptionlabel}[1]%
     {\hspace{\labelsep}\textsf{#1}}

\newcommand{\mycaption}[1]{\begin{quote} 
 \caption{\small#1\normalsize}\end{quote}\vspace{-3em}} 

\newcommand\F{{\mathcal F}}
\newcommand\tr{\protect\psfig{figure=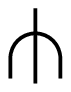,height=.1in}}
\newcommand\C{{\mathbb C}}
\newcommand\wC{\widehat{\mathbb C}}

\newcommand\Z{{\mathbb Z}}
\newcommand\R{{\mathbb R}}

\newcommand\bP{{\mathbb P}}

\newcommand\B{{\mathcal B}}

\newcommand\cC{{\mathcal C}}

\newcommand\ssm{{\smallsetminus}}
\font\bit=cmssi12 at 11truept

\begin{document}

\title{\bf Elliptic Curves as Attractors in ${\mathbb P}^2$\\
 Part 1: Dynamics}

\author{\bf {Araceli Bonifant$^{*}$,}\,\, {Marius Dabija\thanks{Died June 22, 2003.}}\,\,\, and \,\, {John Milnor$^{**}$}\\
\\
$^{*}$Department of Mathematics, University of Rhode Island,
Kingston, RI. 02881-0816.\\
bonifant@math.uri.edu\\
$^{**}$Institute for Mathematical Sciences, Stony Brook University,
 Stony Brook, NY. 11794-3660\\
jack@math.sunysb.edu}

\setcounter{footnote}{1}
\setcounter{equation}{0}

\date{}

\maketitle

\thispagestyle{empty}
\def\IMSmarkvadjust{0 pt}
\def\IMSmarkhadjust{0 pt}
\def\IMSmarkhpadding{0 pt}
\def\IMSpubltext{Published in modified form:}
\def\SBIMSMark#1#2#3{
 \font\SBF=cmss10 at 10 true pt
 \font\SBI=cmssi10 at 10 true pt
 \setbox0=\hbox{\SBF \hbox to \IMSmarkhpadding{\relax}
                Stony Brook IMS Preprint \##1}
 \setbox2=\hbox to \wd0{\hfil \SBI #2}
 \setbox4=\hbox to \wd0{\hfil \SBI #3}
 \setbox6=\hbox to \wd0{\hss
             \vbox{\hsize=\wd0 \parskip=0pt \baselineskip=10 true pt
                   \copy0 \break%
                   \copy2 \break%
                   \copy4 \break}}
 \dimen0=\ht6   \advance\dimen0 by \vsize \advance\dimen0 by 8 true pt
                \advance\dimen0 by -\pagetotal
	        \advance\dimen0 by \IMSmarkvadjust
 \dimen2=\hsize \advance\dimen2 by .25 true in
	        \advance\dimen2 by \IMSmarkhadjust

%
%
  \openin2=publishd.tex
  \ifeof2\setbox0=\hbox to 0pt{}
  \else 
     \setbox0=\hbox to 3.1 true in{
                \vbox to \ht6{\hsize=3 true in \parskip=0pt  \noindent  
                {\SBI \IMSpubltext}\hfil\break
                \input publishd.tex 
                \vfill}}
  \fi
  \closein2
  \ht0=0pt \dp0=0pt
 \ht6=0pt \dp6=0pt
 \setbox8=\vbox to \dimen0{\vfill \hbox to \dimen2{\copy0 \hss \copy6}}
 \ht8=0pt \dp8=0pt \wd8=0pt
 \copy8
 \message{*** Stony Brook IMS Preprint #1, #2. #3 ***}
}

\SBIMSMark{2006/01}{January 2006}{}

\begin{abstract}
A study of rational maps of the real or complex projective plane of degree
two or more, concentrating
on those which map an elliptic curve onto itself, necessarily by an
expanding map. We describe relatively
simple examples with a rich variety of exotic dynamical behaviors which are
perhaps familiar to the applied dynamics community but not to specialists
in several complex variables. For example, we describe
smooth attractors with riddled or
intermingled  attracting basins, and we observe ``blowout'' bifurcations
when the transverse Lyapunov exponent for the invariant
curve changes sign. In the
complex case, the elliptic curve (a topological torus)
can never have a
trapping neighborhood, yet it can have an attracting basin of large measure
(perhaps even of full measure). We also describe examples where there
appear to be Herman rings (that is topological cylinders mapped to themselves
with irrational rotation number) with open attracting basin.
In some cases we provide proofs, but in other cases the discussion
is empirical, based on numerical computation.

\vspace{.3cm}


\vspace{.3cm}

\end{abstract}


\setcounter{lem}{0}
\section{Introduction}

 In the paper
``Self-maps of ${\mathbb P}^ 2$ with invariant elliptic curves'', Bonifant
and Dabija (see \cite{BD}) constructed a number of examples of rational maps $f$ of
the real or complex projective plane of degree $d\ge 2$ with an elliptic curve
${\;\mathcal C}=f({\mathcal C})$ as invariant subset.
This case of a curve of genus one is of particular interest, since genus zero
examples are easy to construct while higher genus examples cannot exist.
(Compare Remarks \ref{rm:g0} and \ref{rm:highg}.)
The present paper studies the extent to which such an invariant
elliptic curve ${\;\mathcal C}\subset{\mathbb P}^2$ can be an ``attractor''.
Here we must distinguish between possible degrees of attraction. By definition,
a compact set $A=f(A)\subset\bP^2$ 
will be called:

\begin{quotation}\noindent
$\bullet$ a {\bit measure-theoretic attracting set\/} if its {\bit attracting
basin\/}, consisting of all points whose orbits converge to $\;A$, has
positive Lebesgue measure, and\smallskip

\noindent$\bullet$ a  {\bit trapped attracting set\/} if there it has a compact
{\bit trapping neighborhood\/} $\;N\;$ such that
$\;f(N)\subset N\;$ and $\;A=\bigcap_n f^{\circ n}(N)$.
\end{quotation}

\noindent  In both cases,
the word {\bit attractor\/} will be reserved for an
attracting set $A$ which contains a dense orbit.\footnote{For some purposes
it is important to have a more restrictive definition, as in \cite{M1}.
Compare the discussion of ``Milnor attractors'' in \cite{Ko} or \cite{ABS}.
For other concepts of attractor, see  Remark \ref{rm:attr}, as well as
\cite{AuBS}; and for exotic examples see \cite{AKYY}, \cite{Kn},
\cite{MMP}, \cite{ABS}, \cite{OS}, \cite{OSAKY}. For
dynamics in $\bP^2(\C)$ see \cite{FS2}, \cite{FW}, \cite{JW}, \cite{Si},
\cite{U2}.}
A measure-theoretic attractor
will be called a {\bit global attractor\/}
if its attracting basin has full measure in the ambient space $\bP^2$.


\smallskip

In both the real and complex cases, we provide examples in which a smooth
elliptic curve $\;\cC\;$ is
a measure-theoretic attractor. In fact there are examples in which there
are two distinct smooth
measure-theoretic attractors whose attracting basins are thoroughly
intermingled, so that they have the same topological closure.
We provide an example of a singular real elliptic curve which is a
trapped attractor; but we prove that a complex elliptic curve
can never be
a trapped attractor. In fact, it seems likely that the attracting basin
of a complex elliptic curve cannot have
interior points, so that the set of points {\it not\/} attracted to
$\;\cC\;$ must be everywhere dense. (Compare Lemma \ref{lem:noint}, as well
as Proposition \ref{prop:bd}.)

We describe examples in which it seems possible that the elliptic curve is
a global attractor so that its attracting basin has full measure.
We also provide a family of examples where there
appears to be a pair of Herman
rings as attractor, with an open neighborhood as attracting basin.
However, in many
cases the dynamical behavior is sufficiently confusing that we are not even
sure what to conjecture.
\medskip

{\bf An Outline of the Following Sections.} Section \ref{s:texp} describes
some basic ideas, including the {\it transverse
Lyapunov exponent\/} along an invariant elliptic curve,

\noindent
 which is a primary
indicator of whether or not the curve is attracting. Methods for actually
computing this transverse exponent will be described in Part 2, the sequel
to this paper; however, the conclusions of these computations are often
quoted below. Section \ref{s:mfi} describes the very restrictive class
of rational maps with a first integral. These are used in Sec.
\ref{s:De} to construct a three parameter family of
more interesting rational maps of degree four. Section \ref{s:ex} studies
eight explicit examples, with conjectured descriptions
based on numerical computation. In the first three examples,
all observed randomly chosen orbits seem to converge to the Fermat curve
$x^3+y^3+z^3=0$, both in the real case with ambient space
$\bP^2(\R)$ and in the complex case with ambient space
$\bP^2(\C)$. This suggests that the real or complex
Fermat curve may be a global attractor, with attracting basin
of full measure. (However, other attractors with basins of extremely small
measure could easily be missed by such random samples.)
Example \ref{e4} suggests that a cycle of two Herman rings can be a
measure theoretic
attractor in $\bP^2(\C)$ (perhaps even be a global attractor)
for such a map with invariant elliptic curve.
In Example \ref{e19}, there is an attracting fixed point at the ``north
pole'', while the ``equator'' forms a measure theoretic attractor, while in
Example \ref{e8} a typical orbit bounces between three near-attractors.
This section concludes with examples of lower degree maps which have a smooth
invariant elliptic curve. In particular,  Example \ref{dg3map} describes
a degree $3$ map of $\bP^2(\C)$
which appear to have the Fermat curve as a global attractor.
All of these conclusions have been empirical.
However, Secs. \ref{s:rb}
and \ref{s:sa} provide cases with explicit
proofs.  Example \ref{e7} describes maps with three different attractors with
thoroughly intermingled basins, all of positive measure. (Compare \cite{AKYY}.)
Two of these basins are dense in the Julia set, while the third basin,
which is everywhere dense, is equal to the Fatou set.
Theorem \ref{th:cas} provides examples of
singular real elliptic curves which are trapped attractors under suitable
rational maps; while Theorem
\ref{th:notrap} shows that a complex elliptic curve
can never  be a trapped attractor.
(We don't know whether non-singular real curves can be trapped attractors.)
Section \ref{s:hr} provides a more general discussion of Herman rings.
The transverse Lyapunov
exponent for a (complex 1-dimensional) Herman ring or Siegel
disk in $\bP^2(\C)$ provides a strict criterion for attraction or repulsion.
This exponent is no longer constant, as it was in the case of an
elliptic curve, but is rather a convex piecewise linear function
on the ring or disk, constant on each invariant circle.
We prove the persistence of invariant circles in $\bP^2(\R)$ under suitable
hypotheses, but our results are not strong enough to prove the conjecture
that the associated Herman rings in $\bP^2(\C)$ are also persistent.
Section \ref{s:?} concludes the discussion by
providing a brief outline of open problems.

We will usually concentrate on the complex case,
although many of the illustrations  will necessarily
illustrate the real case.

\begin{rem}{\bf Computation.}\label{rm:comp} 
Numerical computations are extremely delicate near the invariant curve
${\;\mathcal C}$. Thus it is essential to work with multiple precision
arithmetic; but even so, numerical simulation of the dynamics
must be understood as a hint of the true state of affairs, rather than a
definitive answer. One surprising aspect of these maps is that in some cases
orbits tend to spend quite a bit of time {\it extremely~} close to
${\;\mathcal C}$ even when the transverse exponent is positive.
(Compare Figs. \ref{e13} and \ref{e24}.)
In a similar situation, Maistrenko, Maistrenko and Popovich
( see \cite[p.~2713]{MMP}) report that, in the presence of a small positive value
of the transverse exponent: ``{\it
a trajectory may spend a very long time in the neighborhood of
the invariant subspace. From time to time, the repulsive character of the
chaotic set manifests itself, and the trajectory exhibits a burst in which
it moves far away from the invariant subspace, to be reinjected again into the
proximity of this subspace. $\cdots$ $[$The$]$ positive value of the Lyapunov
exponent applies over long periods of time. For shorter time intervals,
the net contribution $\cdots$ may be negative, and the trajectory is
attracted to the chaotic set.\/}'' (Similar behavior was described in
\cite{PST}.)
\end{rem}

\begin{rem}{\bf Genus Zero.}\label{rm:g0}
Attracting curves of genus zero are easy to construct.
(Compare Example \ref{ex:hr1}.) Thus, for the map
$(x:y:z)\mapsto(x:y:z/2)$, the line $z=0$ is an attracting curve with
the region $|z|^2<|x|^2+|y|^2$ as trapping neighborhood. Similarly,
if $(x:y)\mapsto(f_1(x,y):f_2(x,y))$ is any rational map of $\bP^1$ of
degree $d\ge 2$, then the line $z=0$ is a trapped attracting curve
for the map $(x:y:z)~\mapsto(f_1(x,y):f_2(x,y):z^d)\,. $
In particular, if we start with a map of $\bP^1$ which has a dense orbit
(for example a Latt\`es map---compare Remark \ref{rem:Lat}), then we obtain a trapped attractor.
\end{rem}

\setcounter{lem}{0}
\section{Rational Maps and the Transverse Lyapunov Exponent.}\label{s:texp}
Let $f$ be a rational map of $\bP^2=\bP^2(\C)$. We can write
$~f\,:\,\bP^2\ssm{\cal I}~\to~\bP^2\,$, where
$$f(x:y:z)~=~\big(f_1(x,y,z):f_2(x,y,z):f_3(x,y,z)\big)~, $$
using homogeneous coordinates $(x:y:z)$ where $(x,y,z)\in\C^3\ssm\{(0,0,0)\}$
are well defined up to multiplication by a common non-zero factor.
Here $f_1\,,\,f_2\,,\,f_3$ are to be
homogeneous polynomials of the same degree $d=d(f)\ge 2$
with no common factor, and
${\cal I}={\cal I}(f)$, the {\bit indeterminacy set\/},  is the finite set
consisting of all common zeros of $f_1\,,\,f_2\,,\,f_3$. By definition, $d$
is the {\bit algebraic degree} of $f$. Such a rational map
$f$ is called {\bit regular\/} if
${\cal I}_f$ is vacuous so that $f$ is an everywhere defined holomorphic map
from $\bP^2$ to itself. 
The topological degree of $f$ as a map
from ${\mathbb P}^2$
to itself is then equal to $d^2$, while the algebraic number of fixed points is $d^2+d+1$.
However, if ${\cal I}\ne\emptyset$, then there will be fewer fixed points
(or sometimes an entire curve of fixed points),
and a generic point will have fewer than $d^2$ preimages.

Now consider an algebraic curve $\;\cC\subset\bP^2$, defined by a homogeneous
equation $\Phi(x,y,z)=0$ of degree ${\rm deg}(\cC)
\ge 1$. We will say that a rational map
 $f$ is {\bit well defined\/} on $\;\cC$ if the intersection
$\;\cC\cap{\cal I}_f$ is empty so that $f$ is defined and holomorphic throughout
a neighborhood of $\;\cC$. It then follows that the image $\;\cC'=f(\cC)$
is itself an algebraic curve of degree ${\rm deg}(\cC')\ge 1$. Furthermore,
the degree of the restriction $f|_{\cC}:\cC\to\cC'$ is determined by the
relation\footnote{Proof Outline: A generic line $L$ intersects $\;\cC'$ in
${\rm deg}(\cC')$ distinct points, none of which is a critical value
of $f|_\cC$. Each of these has ${\rm deg}(f|_\cC)$ preimages in $\;\cC$.
On the other hand by Bezout's Theorem, the curve $f^{-1}(L)$
of degree $d(f)$ will intersect
$\;\cC$ in $\,d(f)\,{\rm deg}(\cC)$ points, counted with multiplicity.
In fact, for generic $L$, each of these intersections will be transverse,
and it follows that $d(f)\,{\rm deg}(\cC)=d(f|_{\cC})\,{\rm deg}(\cC')$,
as required.}
\begin{equation}\label{eq:degrees}
  d(f)\;{\rm deg}(\cC)~=~d(f|_{\cC})\;{\rm deg}(\cC')\,.
\end{equation}

{\bf Definition.}
An algebraic curve $\;\cC\subset\bP^2$ will be called {\bit invariant\/}
under $\,f\,$ if
$f$ is well defined on $\;\cC$ and if $f(\cC)=\cC$.
It then follows from  Eq. (\ref{eq:degrees})
that the degree of the restriction $f|_\cC:\cC\to\cC$
is precisely equal to $d(f)$, the degree of the equations which define $f$.

\begin{rem}{\bf Curves of Higher Genus.}\label{rm:highg} As one consequence
of this discussion, it follows that a curve of genus $g\ge 2$ can never be
invariant under a map of $\bP^2$ of
degree $d\ge 2$. For it follows from the Riemann-Hurwitz formula
(see for example \cite{M3})
that a curve of genus $\ge 2$ does not admit any self-maps of degree
$\ge 2$.
\end{rem}

On the other hand, if $\;\cC$ contains points of indeterminacy then these
remarks break down. For example, $\;\cC\ssm({\cal I}\cap\cC)$ may consist
entirely of fixed points, or may map to a single point under $f$.
In cases where
${\cal I}\cap\cC$ may be nonempty, but the image $f\big(\cC\ssm
({\cal I}\cap\cC)\big)$ is contained in $\;\cC$, the curve $\;\cC$ will be called
{\bit weakly $f$-invariant\/}. If $\;\cC$ is smooth and weakly invariant,
then $f$ extends
uniquely to a holomorphic map from $\;\cC$ to itself; but the degree on $\;\cC$ may
be smaller than $d(f)$. (Compare Remarks \ref{rem:fp} and
\ref{rem:degen} below.)

Let ${\mathcal D}(n)=n\,{\mathcal D}(1)$ be the divisor class on $\;\cC$
which is obtained by intersecting $\;\cC$ with a generic curve of degree $n$
in  $\bP^2$. {\it A given holomorphic map
$g:\cC\to\cC$ of degree
$d>0$ extends to a rational map of $\bP^2$ which
is well defined on $\;\cC$ if and only if}
$$g^*{\cal D}(1)~=~{\cal D}(d)\,.$$
(See \cite[\S2]{BD}.)
If $d<{\rm deg}(\cC)$ then this extension is unique;
but if $d\ge{\rm deg}(\cC)$ then there exists a
$3\big(\begin{smallmatrix}d-d(\Phi)+2\\
 2 \end{smallmatrix}\big)$ dimensional family of such extensions
since we can always replace the associated homogeneous polynomial map
$F:\C^3\to\C^3$ by $F+\Phi H$ where $\Phi=0$ is the
defining equation for $\;\cC$ with $d(\Phi)={\rm deg}(\cC)$,
and where $H:\C^3\to\C^3$ is any homogeneous map
of degree $d-d(\Phi)$. 
In this case, a {\it generic\/} extension will be
regular (that is, defined and holomorphic everywhere).

One particular virtue of an elliptic curve ${\mathcal C}\subset\bP^2$
is that any holomorphic self-map must be
linear in terms of suitable coordinates. In fact such a
curve is conformally isomorphic to some flat torus ${\mathbb C}/\Omega$, where
$\Omega$ is a lattice. More precisely, there is a holomorphic uniformizing
map
$$\upsilon:{\mathbb C}/\Omega~\to~{\mathcal C} $$
which is biholomorphic in the case of a smooth elliptic curve, and is
one-to-one except on finitely many singular points in the case of a singular
curve.
Any holomorphic self-map of ${\mathcal C}$ lifts to a holomorphic self-map of
${\mathbb C}/\Omega$ which is necessarily linear, $t\mapsto at+b$.
It follows easily that the normalized Lebesgue measure on ${\mathbb C}/\Omega$
pushes forward
to a canonical smooth probability measure $\lambda$ on ${\mathcal C}$ which
is invariant under every non-constant self-map.
The derivative $a$ of the linear map on $\C/\Omega$ will be called the
{\bit multiplier\/} of $f$ on $\;\cC$. Note that the product $a\,\Omega$ 
is a sublattice of finite index in $\Omega$, and that $|a|^2$ is equal to
the index of this sublattice. Equivalently, $|a|^2$ is the
topological degree of $f$ considered as a map from $\;\cC$ to itself.
In particular, $|a|^2$ is equal to the algebraic degree $d$  of $f$ whenever
$\;\cC\subset\bP^2$ is invariant under the rational map $f$.  Since
we always assume that $d\ge 2$, it follows that this
canonical measure $\lambda$ is ergodic.

In the case of an elliptic curve defined by equations with real coefficients,
the real curve
$${\mathcal C}_{\mathbb R}~=~{\mathcal C}\cap{\mathbb P}^2({\mathbb R})$$
has at most two connected components. In the case where there are two
components, if $\cC_\R$ is mapped into itself by a rational map $f$ of
${\mathbb P}^2({\mathbb R})$,
then at least one of the two components, say
$\cC_\R^0$,  must map onto itself under $f$ (or at least under $f\circ f$).
In this case we have a uniformizing map
${\mathbb R}/{\mathbb Z}\to{\mathcal C}_{\mathbb R}^0$, such that $f$
(or $f\circ f$)
corresponds to a map on ${\mathbb R}/{\mathbb Z}$
which is linear with constant integer multiplier. Such an invariant component
${\mathcal C}_{\mathbb R}^0$ has a canonical invariant probability measure.
\smallskip

{\bf The Transverse Lyapunov Exponent.} Let $\;\cC$ be a smooth curve, invariant
under the rational map $f$. We will describe an associated real number
which is conjectured to be negative if and only if $\;\cC$ is a measure
theoretic attractor.
To fix ideas we will concentrate on the complex case, but constructions in
the real case are completely analogous. 
The notation $T{\mathbb P}^2|_{{\mathcal C}}$ will be used for the
complex 2-plane bundle of vectors tangent to ${\mathbb P}^2({\mathbb C})$ at
points of the submanifold ${\mathcal C}$, and the abbreviated notation
$T_{\tr\;}\cC$ will be used for the ``transverse''
complex line bundle over ${\mathcal C}$ having the quotient vector space
$$T_{\tr}({\mathcal C}, p)=T({\mathbb P}^2,p)/T({\mathcal C},p)$$
as typical fiber. In other words, there is a short exact sequence
$~  0\;\to\;T\cC\;\to\;T\bP^2|_\cC\;\to\;T_{\tr\;}\cC\;\to\;0~ $
of complex vector bundles over $\;\cC\,$.
It is sometimes convenient to refer to $\;T_{\tr\;}\cC\;$ as the
``{\bit normal bundle\/}'' of $\;\cC\,$, although that designation isn't
strictly correct. If $\;f:\bP^2\to\bP^2\;$ with $\;f(\cC)\subset\cC\,$,
then $f$ induces linear maps
\begin{equation}\label{eq:n-der}
    f'_{\tr}(p):T_{\tr}({\mathcal C},p)~\to~ T_{\tr}({\mathcal C}, f(p))~,
\end{equation}
and these linear maps collectively form a  fiberwise linear self-map
$f'_{\tr}:T_{\tr\;}{\mathcal C}\to T_{\tr\;}{\mathcal C}~.$

Now choose a metric on this complex normal bundle. That is, choose a norm
$\|\vec{v}\,\|_{\tr}$ on each quotient
vector space $T_{\tr\;}{\mathcal C}$ which depends
continuously on $\vec{v}$, vanishes only on zero vectors, and satisfies
$\|t\vec{v}\|_{\tr}=|t|\,\|\vec{v}\|_{\tr}$.
Then the linear map $\;f'_{\tr}\;$ of Eq. (\ref{eq:n-der}) has an
operator norm
\begin{equation}\label{eq:norm}
    \|f'_{\tr}(p)\|
  ~ =~ \|f'_{\tr}\vec{v}\;\|_{\tr}\;/\|\vec{v}\|_{\tr}
\end{equation}
which is well defined and satisfies the chain rule. Here
$\vec{v}$ can be any non-zero vector in the fiber $T_{\tr\;}({\mathcal C},p)$
over $p$. By definition, the {\bit transverse Lyapunov exponent\/} along the
 invariant curve $\mathcal C$ is equal to the rate of exponential growth
$${\rm Lyap}_{\,\cC}~=~
 \lim_{k\to\infty}\;(1/k)\,\log\|(f^{\circ k})_{\tr}'(p)\| $$
for almost every choice of initial point $p\in{\mathcal C}$.
By the Birkhoff Ergodic Theorem, this coincides with the average
value
$${\rm Lyap}_{\,\cC}(f)~=~\int_{\mathcal C}
\log\|f'_{\tr}(p)\|\,d\lambda(p) ~.$$
Using the fact that the measure $\;\lambda\;$ is invariant
under $\;f\,$, it is not hard to check directly that this average value
is independent of the choice of metric.

Thus a negative value of ${\rm Lyap}_{\,\cC}$ means that under iteration of
$f$ almost any point which is
``infinitesimally close'' to ${\mathcal C}$ will
converge towards ${\mathcal C}$. A key role in this case is played by the
stable sets of the various points $\;p\in\cC\,$. By definition, the
{\bit stable set\/} of $\;p\;$ is the union of all connected sets containing
$\;p\;$ for which the diameter of the $\;n$-th forward image tends to
zero as $\;n\to\infty\,$. Many such stable sets are smooth curves.
With a little imagination, some of these
are clearly visible in Figs. \ref{fig:min} through \ref{e21}. It is natural
to expect that negative values of $~{\rm Lyap}_{\,\cC}~$ will imply that
${\mathcal C}$ is a measure-theoretic attractor. (A sketch of a proof is
given in \cite{AKYY}; and a very special case
will be proved in Theorem \ref{th:el} below.)
On the other hand, if ${\rm Lyap}_{\,\cC}>0$
then almost any ``infinitesimally close'' point will be pushed away from
${\mathcal C}$. It seems natural to conjecture that positive values
of ${\rm Lyap}_{\,\cC}$ should imply that the
attracting basin of ${\mathcal C}$ has measure zero. However, this seems
like a difficult question. (Compare Remark \ref{rm:Lpos}.)

 The term
{\bit blowout bifurcation\/} has been introduced in \cite{OS} for a transition
in which a transverse
Lyapunov exponent crosses  through zero. (Compare \cite{MMP} or \cite{AAN}.)

\setcounter{lem}{0}
\section{Maps with First Integral.}\label{s:mfi}
By definition, a {\bit first integral\/} for a dynamical system is a
non-constant function which is constant on each orbit. In particular, by a
{\bit first integral\/} for a rational map $f:\bP^2\ssm{\cal I}_f\to\bP^2$ we
will mean a rational function $\;\eta:\bP^2\ssm{\cal I}_\eta\to\bP^1\,,$
with values in the projective line, which satisfies
\begin{equation}\label{eq:int}
\eta\big(f(x:y:z)\big)~=~\eta(x:y:z)
\end{equation}
whenever both sides are defined. Identifying $\bP^1$ with the Riemann sphere
$\wC=\C\cup\{\infty\}$, we will write
$$\eta(x:y:z)~= ~\Phi(x,y,z)/\Psi(x,y,z)\,\in\,\wC~, $$
where
$\Phi$ and $\Psi$ are homogeneous polynomials of the same degree without common
factor. It follows that $\bP^2$ has a somewhat singular ``foliation'' by
algebraic curves of the form
\begin{equation}\label{eq:fol}
\alpha\,\Phi(x,y,z)+\beta\,\Psi(x,y,z)~=~0~,
\end{equation}
which intersect only in the finite set ${\cal I}_\eta$ consisting of common
zeros of $\Phi$ and $\Psi$, and such that each of these curves is weakly
invariant under $f$ as defined in \S\ref{s:texp}.
(A point can be contained in two such curves only if it is either periodic,
or a point of indeterminacy for $f$.) Such maps are exceedingly special.
For example we have the following.
\smallskip

\begin{Quote}
\begin{lem}{\bf Some Maps with First Integral.}\label{lem:fimaps}
Let $f$ be a rational map with first integral, such that
a generic point of
$\bP^2$ is contained in an elliptic curve which is mapped to itself
with degree $d\ge 2$. Then:

$\bullet$ There are no dense orbits, since every orbit is contained in a
weakly invariant curve.

$\bullet$ Periodic points, repelling along this curve,
are everywhere dense. Hence the {\bit Fatou
set\/}, that is the union of open sets on which the iterates of $f$ form a
normal family, is empty.

$\bullet$ For most values of $n$ there are {\it infinitely many\/} fixed
points of $f^{\circ n}$, with at least one in each weakly invariant curve.
Hence there must be an entire algebraic curve of such points.

$\bullet$ The indeterminacy set ${\cal I}_f$ is necessarily non-empty.
\end{lem}
\end{Quote}

The last statement follows
since there are infinitely many points of fixed period $n$, or since a generic
point has $d$ preimages in the invariant curve which passes through it,
but also just $d$ (rather than $d^2$) preimages in all of $\bP^2$.
The other statements are easily verified.\qed\smallskip

Here is a class of examples which generalize a
construction due to A. Desboves (see \cite{De}). For any smooth cubic curve
$\;\cC\subset\bP^2$, there is a canonical map $f:\cC\to \cC$
called the {\bit tangent process\/}, constructed as follows. For any
point $p\in \cC$, let $L_p\subset\bP^2$
be the unique line which is tangent to $\;\cC$ at $p$. Then the image
$f(p)$ is defined by the equation
$$L_p\cap\cC~=~\{p\}\,\cup\,\{f(p)\}~.$$
In fact if we choose the parametrization $\upsilon:\C/\Omega\to\cC$ of
\S\ref{s:texp} correctly, then three points $t_j$ of $\C/\Omega$
will correspond to
collinear points of $\;\cC$ if and only if $t_1+t_2+t_3=0$. In our case, since
there is a double intersection at $p$, we obtain the equation $2t_1+t_3=0$
or $t_1\mapsto t_3=-2t_1$. Thus $f$ has multiplier $-2$ and degree 4.

Now start with two
distinct cubic curves in $\bP^2$, described by homogeneous equations
$\Phi(x,y,z)=0$ and $\Psi(x,y,z)=0~.$
Then there is an entire one-parameter family of such curves (Eq. (\ref{eq:fol}))
filling out the projective plane. In fact, any point of $\bP^2$ which is
not a common zero of $\Phi$ and $\Psi$ belongs to a unique curve
$$\label{eq:ic}
\Phi/\Psi~=~{\rm constant}~=~-\beta/\alpha~\in~\widehat\C~. $$
If a generic curve in our one parameter family is smooth, then a generic
point $p\in\bP^2$ belongs to a unique smooth curve $\;\cC_p$ in the family.
Applying the tangent process at $p$,
we obtain a well defined image point $f(p)\in
\cC_p$. Since $p$ is generic, this extends to a uniquely defined rational map
of $\bP^2$ which carries each curve of our family into itself.

\begin{figure} 
\centerline{\psfig{figure=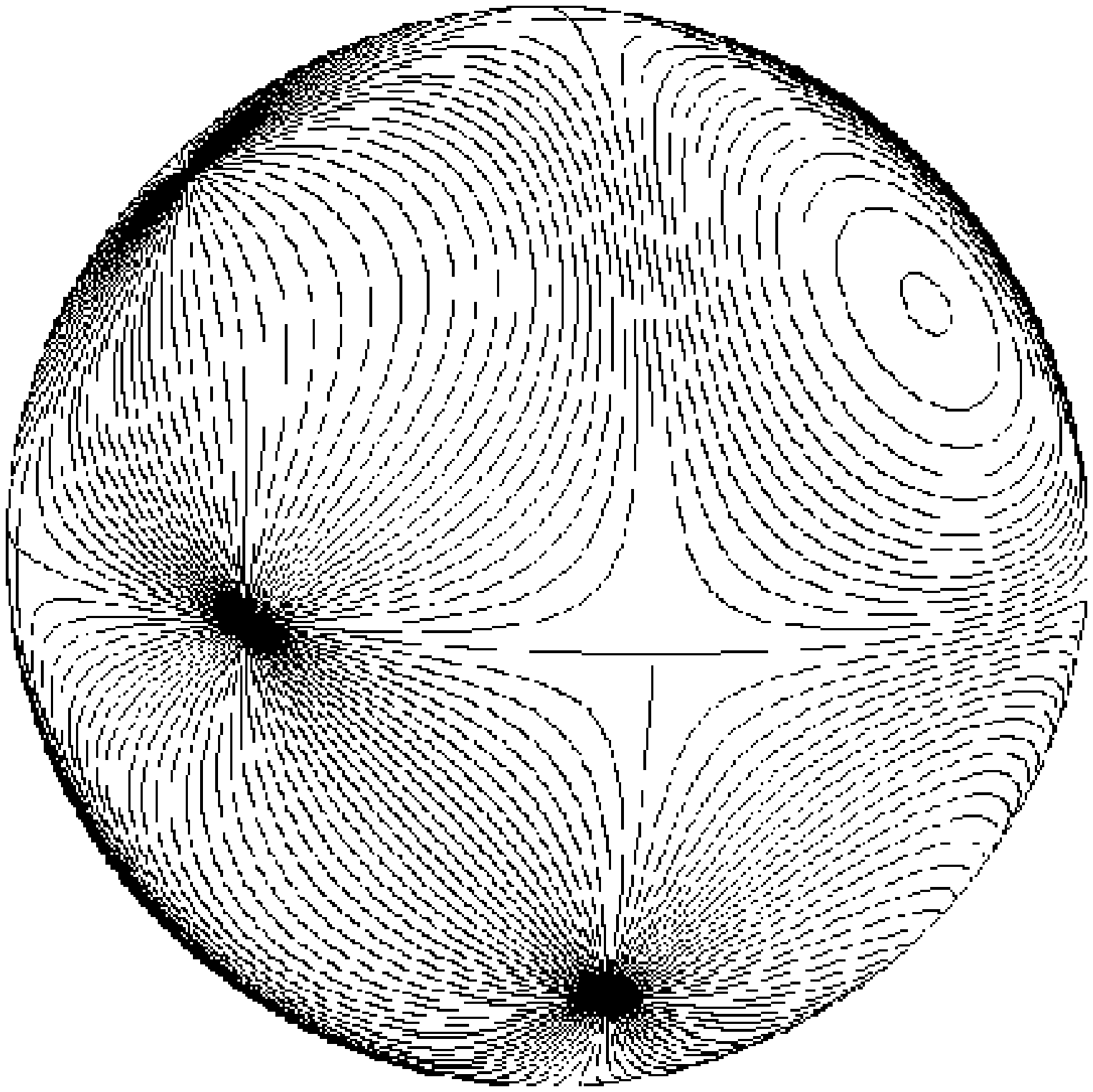,height=3in}}
\mycaption{ \bit \label{fg:rpp} ``Foliation'' of the real projective plane by the family of elliptic curves
$\;\Phi/\Psi=k\in\widehat\C\,$, where $\;\Phi=x^3+y^3+z^3\;$ and $\;\Psi=3xyz\,$.
Here $\;{\mathbb P}^2({\mathbb R})\;$ is represented as a unit 2-sphere
with antipodal points identified. These real curves
intersect only at their three common inflection points, which look dark in
the figure. In the limiting
case as $\;k\to\infty\,$, the curve $\Phi=k\Psi$
degenerates to the union $\;xyz=0\;$ of the
three coordinate lines, which intersect at the points
$\;(-1,0,0)\,,\,(0,1,0)\;$ and $\;(0,0,1)\,$ respectively near
the left, top, and center of the figure.}
\end{figure}

Let us specialize to the classical example, with
$\; \Phi(x,y,z)=x^3+y^3+z^3\;$ and $\;\Psi(x,y,z)=3xyz~,$ as studied by
Desboves (where the factor 3 is inserted for later convenience).
The corresponding foliation of the real projective plane  $\bP^2(\R)$ by
 the curves of Eq. (\ref{eq:fol}) is illustrated in Fig. \ref{fg:rpp}.
This foliation has three kinds of singularities, all represented
in the figure. There are:

{\bf (a)} Three singularities where two of the three
coordinates $x,y,z$ are zero. These all lie in the real plane $\bP^2(\R)\subset
\bP^2(\C)$.

{\bf (b)} Three singularities in the real plane (or nine in the complex plane)
where all of these curves
intersect. Each of these lies along just one of the three coordinate axes.

{\bf (c)} One real singularity (or nine complex singularities) where
$x^3=y^3=z^3$, represented by the center in the upper right of Fig.\ref{fg:rpp}.
\smallskip

According to  Desboves, the tangent processes for these various
curves $\,\Phi/\Psi=k\,$ fit together to yield a well defined rational map
$f_0:\bP^2\ssm{\cal I}\to\bP^2$ which is given by the formula
\begin{equation}\label{dbm}
f_0(x:y:z)~=~\big(x(y^3-z^3):y(z^3-x^3):z(x^3-y^3)\big)~.
\end{equation}
The indeterminacy set ${\cal I}={\cal I}(f_0)$ for this {\bit classical
Desboves map\/} consists of the twelve
points of type {\bf (a)} and {\bf (c)}, as listed above.
This particular example has the advantage (as compared with a generic choice
for $\Phi$ and $\Psi$) that most 
of the curves in our one parameter family are non-singular, and contain
no points of indeterminacy. The only exceptions are the
curves $\Phi=k\Psi$ with $k^3=1$, which are singular at points of
indeterminacy of type $\bf(c)$, and the degenerate case
$\Psi=0$ (corresponding to $k=\infty\in\widehat\C$)
with singular indeterminacy points of type $\bf(a)$.
The foliation singularities of type {\bf(b)}, where all of the curves
intersect, are all fixed points at which the value $f_0(p)=p$ is well defined.

For further examples of rational maps of ${\mathbb P}^2$
with first integral, see Example \ref{dg3map} and Remark \ref{rem:degen}.

\setcounter{lem}{0}
\section{ The Desboves Family.}\label{s:De}
Let $\Phi(x,y,z)$ be the homogeneous polynomial $x^3+y^3+z^3$.
The {\bit Fermat
curve~} $\F$ is defined as the locus of zeros $\Phi(x,y,z)=0$ in the
projective plane $\bP^2$.
(Here we can work either over the real numbers or over the complex numbers.)
Most of the examples in \S\ref{s:ex} will belong to a family of
$4^{\rm th}$-degree rational maps of $\bP^2$ which carry this Fermat curve
into itself, as introduced in \cite[\S6.3]{BD}. We will call these
{\bit Desboves maps\/}, since they arise from a simple perturbation of
the classical Desboves map $f_0$ of Eq.~(\ref{dbm}). Evidently $f_0$
lifts to a homogeneous polynomial map
$$F_0(x,y,z)~=~\Big(x(y^3-z^3)\,,~y(z^3-x^3)\,,~z(x^3-y^3)\Big) $$
from $\;{\mathbb C}^3\;$ to itself.
Geometrically, $f_0$ is defined by the property that the line from $\;p\;$ to
$\;f_0(p)\;$ is tangent to the elliptic curve $\,(x^3+y^3+z^3)/3xyz=k\,$
which passes through the point $\;p\,$. Its set
of fixed points on each smooth curve in our family coincides with the
intersection
$$x^3+y^3+z^3=3xyz=0~, $$
and can also be identified with the set of points of
inflection on any one of these curves, or as the set of points where
all of these curves intersect. This map $\;f_0\;$ is not everywhere
defined: it has a twelve point set of points of indeterminacy as
described at the end of \S\ref{s:mfi}. However, for any specified curve
$\Phi_k(x,y,z)=x^3+y^3+z^3-3kxyz=0\;$ in our family,
if we replace $\;F_0\;$ by the sum
$$ F_L(x,y,z)~=~F_0(x,y,z)\,+\,L(x,y,z)\,\Phi_k(x,y,z) $$
where $L$ is any linear map from
${\mathbb C}^3$ to itself, then we obtain a new map $f_L$ of ${\mathbb P}^2(
{\mathbb C})$ which coincides with $f_0$ on the particular
curve $\Phi_k(x,y,z)=0\,$.
For a generic choice of $L$, the resulting map $f_L$ of
${\mathbb P}^2({\mathbb C})$ is well defined everywhere.\smallskip

To simplify the discussion, we will restrict attention to the case $k=0$,
taking $$\Phi(x,y,z)~=~\Phi_0(x,y,z)~=~x^3+y^3+z^3~,$$
and will take a linear map $L$ which is described by a diagonal matrix,
$$~ L(x,y,z)~=~(a\,x\,,\,b\,y\,,\,c\,z)~.$$

\begin{definition}{\bf Desboves Maps.} The resulting
3-parameter family of  maps of the real or complex projective plane
will be called the family of {\bit Desboves maps\/}.
These maps $f=f_{a,b,c}$ are given by the formula
\begin{equation}\label{e11}
 f(x:y:z)=\Big(x(y^3-z^3+a\,\Phi)\,:\,y(z^3-x^3+b\,\Phi)\,:\,
z(x^3-y^3+c\,\Phi)\Big)~,
\end{equation}
where $a,b,c$ are the parameters. Each such $f$ maps the Fermat
curve $\F$, defined by the equation $\Phi(x,y,z)=0\,$, into itself.
Furthermore, each $f$
maps each of the coordinate lines $x=0$ or $y=0$ or $z=0~$ into itself.

For special values of the parameters, the map $f$ may have points of
indeterminacy (but never on the curve). However, for a  generic choice
of parameters $f$ is everywhere defined. More explicitly,
it is not hard to see that $f$ is an everywhere defined holomorphic map from
$\bP^2$ to itself if and only if we avoid a union of seven hyperplanes in
the space $\C^3$ of parameters, defined by the equation
\begin{equation}\label{eq:bad}
abc(a+b+c)(a+1-b)(b+1-c)(c+1-a) ~=~ 0~.
\end{equation}\end{definition}

\begin{rem}{\bf Fixed Points.}\label{rem:fp}
Generically, each Desboves map has 21 distinct fixed points
(nine with $xyz\ne 0$, nine on the Fermat curve with just one of the
coordinates equal to zero, and three with two coordinates equal to zero).
However, there are two kinds of exception:

\begin{quotation}
\noindent$\bullet$
If the product $~a\,b\,c\,(a+b+c)~$ is zero, then one or more of the fixed
points will be replaced by an indeterminacy point.\smallskip

\noindent$\bullet$ If one or more of the differences $b-a\,,~c-b$, or $a-c$
is equal to $+1$, then there is not only an indeterminacy point
but also an entire line of fixed points.
\end{quotation}

\noindent In any case, there are exactly nine fixed points
on the complex Fermat curve $\F$, forming the
intersection of $\;\F\;$ with the locus $\;xyz=0\,$. Consider
for example the three points $(0:{\root 3\of {-1}}:1)$ obtained
by intersecting the curve ${\F}$ with the invariant line  $x=0$. If we
introduce the coordinate $~Y=y/z\in\widehat\C~$ on this line, then the
restriction of the map $f$ to this line is a rational map given by the formula
$$  Y~\mapsto~ Y\,\frac{b\,Y^3+(b+1)}{(c-1)\,Y^3+c}~,$$
with fixed points at $Y=0$, at $Y=\infty$, and at the three points
$Y^3=-1$. 
A brief computation shows that the derivative of this one variable
map at these five fixed points is respectively
\begin{equation}\label{eq:tr-der}
(b+1)/c\,,~ (c-1)/b 
\,,~{\rm and}~~ 3(c-b)-2\,,\,\,
{\rm where~the~last~is~ counted~three~times~.}
\end{equation}
(Something very special occurs when $c=b+1$. In that case, all five derivatives
are $+1$, and in fact it is easy to check that {\it every\/} point
on the line $x=0$ is fixed under $f$.)
Similarly, permuting the coordinates cyclically, we obtain corresponding
formulas for the invariant lines $\;y=0\;$ and $\;z=0\,$.
\end{rem}

\begin{rem}{\bf  The Attracting Basin.}\label{rem:basin}
According to \cite[Theorem 5.4]{BD}, if $~\cC=f(\cC)~$ is any invariant elliptic
curve, then the set of iterated preimages of any point of $\;\cC$ is everywhere
dense in
the Julia set $J(f)$ (the complement of the Fatou set).
It seems likely that the following further statement is true:

\begin{Quote}
\begin{conj} The entire attracting
basin $\;{\mathcal B}(\cC),$\; consisting of all points
whose forward orbits converge to $\;\cC$, is contained in  the Julia set $J(f)$;
hence the closure $\overline{{\mathcal B}(\cC)}$ is precisely equal to
$J(f)$.
\end{conj}
\end{Quote}

\noindent An immediate consequence would be the following.

\begin{Quote}\begin{lem}{\bf No Interior Points?}\label{lem:noint}
If this conjecture is true, then for any map $f_{a,b,c}$ in the Desboves
family, the attracting basin  $\B(\F)$ has no interior points. In other words,
the closure of the complementary set
\noindent
 $~\bP^2(\C)\ssm {\mathcal B}(\F)\,$,
consisting of points which are not attracted
to $\;\F\,$, is the entire plane $\bP^2(\C)$.
\end{lem}\end{Quote}

{\bf Proof.}
The conjecture asserts
that points outside of the Julia set are not in the basin
$\B={\mathcal B} (\F)$, so we need only show that every point of $J$
can be approximated by points outside of $\B$. 
Since the average of the differences $c-b\,,~b-a\,,~a-c$ is zero,
it follows from (\ref{eq:tr-der}) that the average of the transverse
derivatives at the fixed points of $f$ in $\F$ is $-2$. Hence at least one
of these fixed points is strictly repelling. Suppose for example that the
point $(0:-1:1)$ is repelling within the line $x=0$. Then the intersection
of the basin ${\mathcal B}$ with this line consists only of countably
many iterated preimages of $(0:-1:1)$. Therefore there are points $(0:y:z)$
arbitrarily close to $(0:-1:1)$ which are {\it not} in this basin. Since every
point of $J$ can be approximated by iterated preimages of $(0:-1:1)$, it can
also be approximated by iterated preimages of such points $(0:y:z)$,
as required. \qed

\end{rem}

On the other hand, some of these fixed points on $\F$ may be attracting
in the transverse direction. For example, if $~|3(c-b)-2|<1\,$
then each of the three fixed points
where $\F$ intersects the line $x=0$ is a {\bit saddle\/}, repelling along the
Fermat curve, but
attracting along this line, which intersects it transversally. The
{\bit stable manifold\/} for such a saddle point can be identified with its
immediate attracting basin within the line $x=0$. It is not hard to check
that this stable manifold is contained in the Julia set. Hence its
iterated preimages must be dense in the Julia set.

The attraction within this stable manifold
will be particularly strong if $c-b=\frac{2}{3}$, so that the
transverse derivative $3(c-b)-2$ is zero, or in other words so that the
associated fixed point is transversally superattracting.
Similarly, the transverse derivative at the three points where $y=0$
(or where $z=0$) is zero if and only if $a-c=\frac{2}{3}$ (or respectively
$b-a=\frac{2}{3}$).

\begin{figure}
\centerline{\psfig{figure=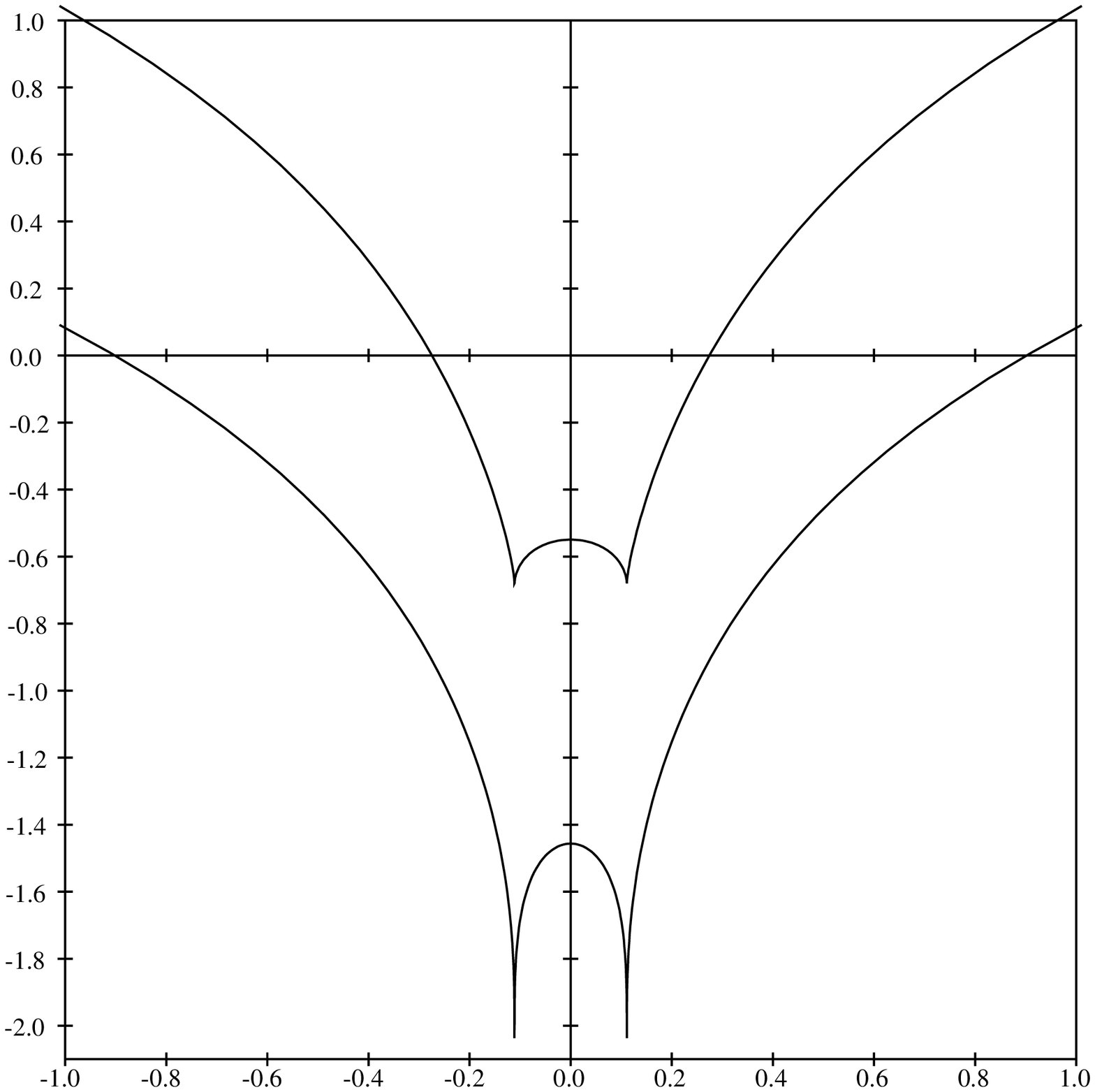,height=3in}}
\mycaption{\label{fig:2/3}\bit
Graph of the transverse exponents along the real and complex Fermat curves
as functions of the middle parameter $b$ for the ``two-thirds family''
of Definition \ref{def:2/3}, with Desboves parameters
$(b-\frac{2}{3}\,,~b\,,~b+\frac{2}{3})\,$. The lower graph represents
the real case, with a transverse exponent which is strictly smaller
(more attracting). In both cases the function is even, with a sharp
minimum at $b=\pm 1/9$. The box encloses the region $-1\le b\le 1$ with
$ -2.1\le{\rm Lyap}_{\F}\le 1$.
}
\end{figure}

\begin{definition}{\bf The Two-Thirds Family.}\label{def:2/3}
 We will say that $f$ belongs to the {\bit two-thirds
family~} if two of the three differences $\;b-a\,,~c-b\;$ and $\;a-c\;$
are equal to $\;\frac{2}{3}\,$, or equivalently if two out of three of the
fixed points on the Fermat curve are transversally superattracting. Since
the average of the three values of the transverse derivative is always $-2\,$,
if two out of three values are zero then it follows that the remaining
one is $-6$ (rather strongly repelling). To fix our ideas, let us suppose that
\begin{equation}\label{eq:2/3}
(a,b,c)~=~(b-\textstyle{\frac{2}{3}}\,,~b\,,~b+\textstyle{\frac{2}{3}})\,,
\end{equation}
so that the transverse derivative is zero when $x=0$ or $z=0$. The associated
transverse Lyapunov exponent, plotted as a
function of $b$, is shown in Fig.~\ref{fig:2/3}. In the real case,
this transverse exponent is negative (that is attracting)
if and only if $|b|<0.901\cdots$, while in the complex
case is is negative if and only if $b<0.274\cdots$. See Part 2 of this
paper for such computations.
\end{definition}

\begin{rem}{\bf Symmetries.}\label{rem:sym} We conclude this section
with some more technical remarks.
If we permute the three parameters $(a,b,c)$ cyclically, then clearly
we obtain a new map $f_{b,c,a}$ which is holomorphically conjugate to
$f_{a,b,c}$. We can generalize this construction very slightly by allowing
odd permutations also, but changing signs. If $S_3$ is the symmetric
group consisting of all permutations $i\mapsto\sigma_i$ of the three symbols
$\{1,2,3\}$, then $S_3$ acts as a group of rotations of $\R^3$ or $\C^3$ as
follows. For each $\sigma\in S_3$, consider the {\bit sign-corrected
permutation of coordinates\/}
\begin{equation}\label{eq:sgroup}
\widehat\sigma(z_1\,,\,z_2\,,\,z_3)~=~{\rm sgn}(\sigma)\,\big(z_{\sigma_1}\,,\,
z_{\sigma_2}\,,\,z_{\sigma_3}\big)~.
\end{equation}
Then a brief computation shows that the homogeneous map
 $F_{a,b,c}$ of $\R^3$ or $\C^3$ is linearly conjugate to the map
$$ F_{\widehat\sigma(a,b,c)}~=~ \widehat\sigma\circ F_{a,b,c}\circ\widehat
\sigma^{-1}~. $$
It follows that the associated map $f_{a,b,c}$ of the projective plane
is holomorphically conjugate to the map $f_{\widehat\sigma(a,b,c)}$. One can check
that these are {\it the only\/} holomorphic conjugacies between Desboves maps
(for example, by making use of the eigenvalues of the first derivative of $f$
at the 21 fixed points).

In the complex case, note also that each Desboves map $f$ commutes with a
finite group $G\cong{\mathbb Z}/3\times{\mathbb Z}/3~$ of symmetries
of the projective plane. In fact
$~f\circ g=g\circ f~$ for each $g$ in the group $G$ consisting of all
automorphisms
\begin{equation}\label{eq:sym}
 g(x:y:z)~=~(\alpha x:\beta y:z)\qquad{\rm with}\qquad \alpha^3~=~\beta^3~=~1~.
\end{equation}\end{rem}

\begin{rem}{\bf A Closely Related Map.}\label{rem:Lat}
It is sometimes convenient to eliminate these last symmetries by passing
to the quotient space $\bP^2/G$ which is
isomorphic to $\bP^2$ itself, but with coordinates $(x^3:y^3:z^3)$.
If we introduce variables $\;\xi=x^3\;,~\eta=y^3\,$, and $\;\zeta=z^3\,$,
and set $\;\phi=\xi+\eta+\zeta\,$,
then the map $\;(x:y:z)\mapsto(\xi:\eta:\zeta)\;$
transforms the Fermat curve $\;\F\;$ to a line $\;\phi=0\,$.
Under this transformation,
the Desboves map (\ref{e11}) is semiconjugate to a different rational map
$$  (\xi:\eta:\zeta)~\mapsto~\Big(\xi(\eta-\zeta+a\phi)^3:
\eta(\zeta-\xi+b\phi)^3:\zeta(\xi-\eta+c\phi)^3)\Big)~,$$
also of degree four. Evidently
this new map carries the line $\;\phi=0\;$ into itself by a {\bit Latt\`es
map\/}, that is the image of a rigid torus map under a holomorphic
semiconjugacy. (Compare \cite[2006a]{M2}.)
\end{rem}

\setcounter{lem}{0}
\section{ Empirical Examples.}\label{s:ex}

 This section will provide empirical discussions of six examples from the
Desboves family of degree four maps, as described in \S\ref{s:De}, plus two
examples of lower degree maps. 
Four of the six examples from the Desboves family,  belong to the
``two-thirds'' subfamily of Definition \ref{def:2/3}.\smallskip

\begin{figure} 
\centerline{\psfig{figure=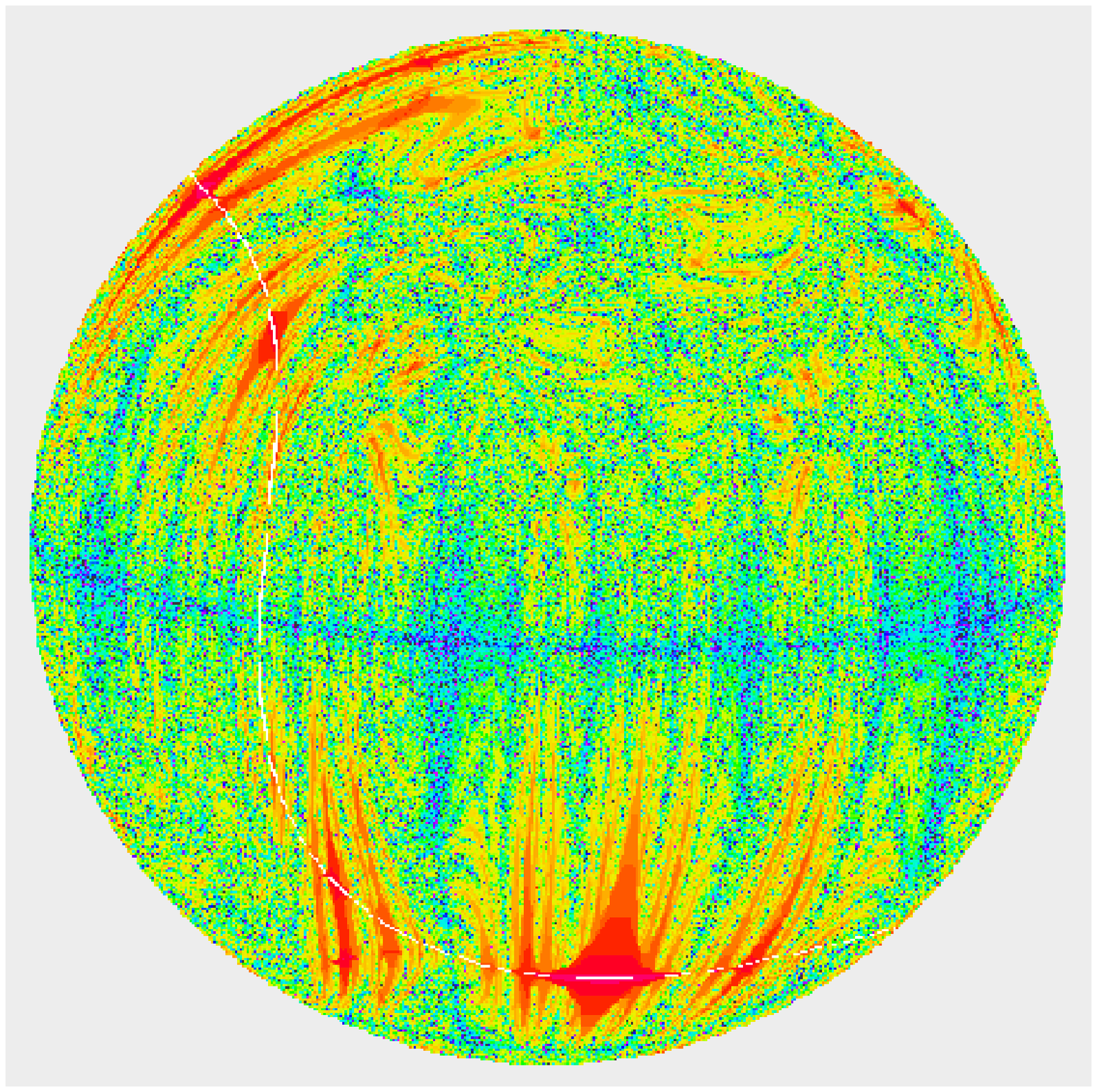,height=3in}\quad
\psfig{figure=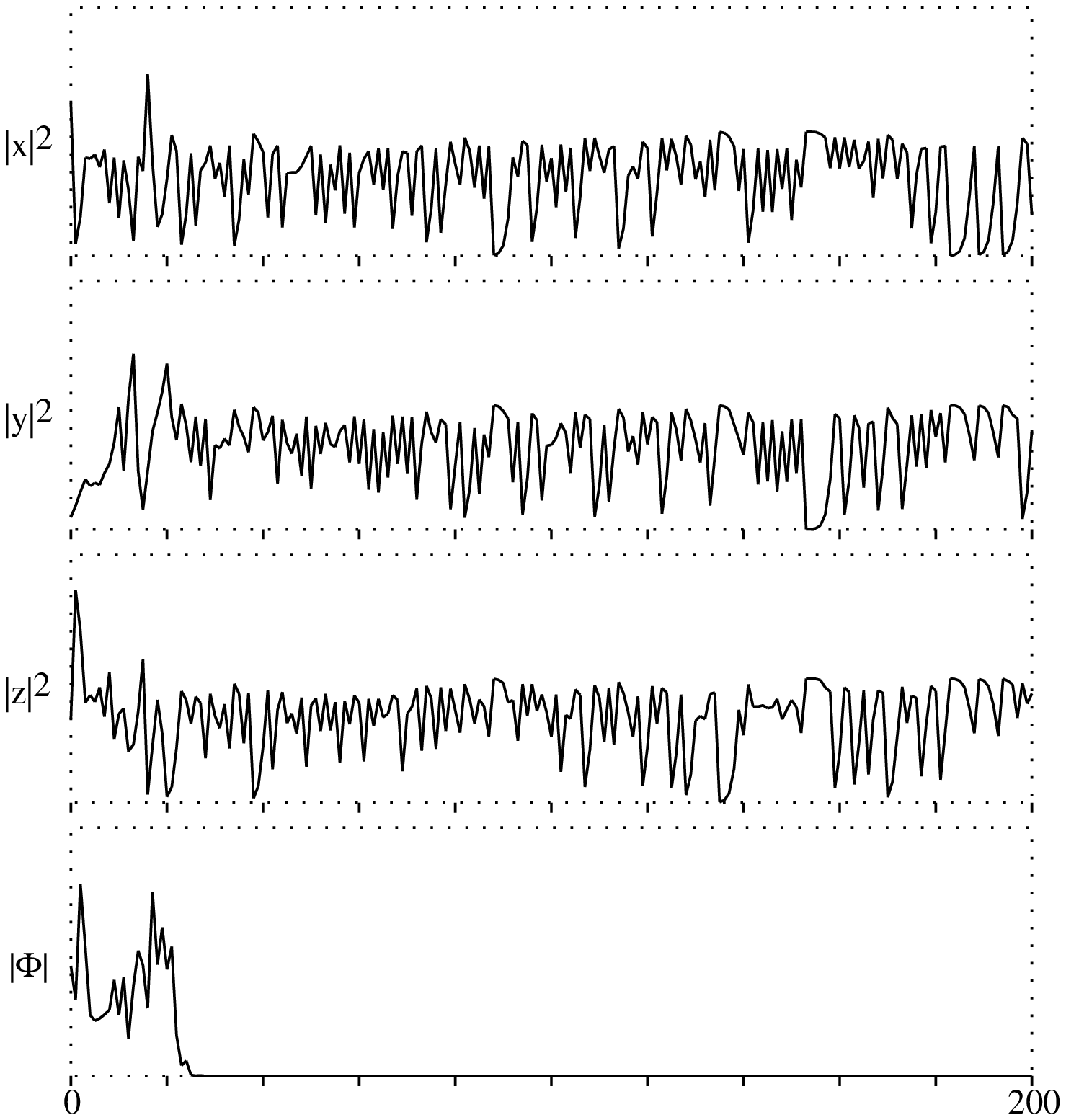,height=3in}}
\mycaption{\bit \label{fig:min} $($Example \ref{ex:2}.$)$ On the left:
dynamics on the real projective plane for the Desboves map in the
two-thirds family 
with parameters $\;(a,b,c)=(-\frac{5}{9}\,,\,\frac{1}{9}\,,\,\frac{7}{9})\,$.
The sphere is oriented as in Fig.~\ref{fg:rpp}.
On the right: plot of $\;|x|^2\,,~|y|^2\,,~|z|^2\;$ and $\;|\Phi|\;$
as functions of the number of iterations for a typical randomly chosen
complex orbit.
Here each orbit point $\;(x:y:z)\;$ has been normalized so that $\;|x|^2+|y|^2
+|z|^2=1\,$. In this run, it took 23 iterations to come close enough to
the Fermat curve so that $\;|\Phi|\;$ appears to be zero on the graph.}
\bigskip \bigskip \bigskip

\centerline{\psfig{figure=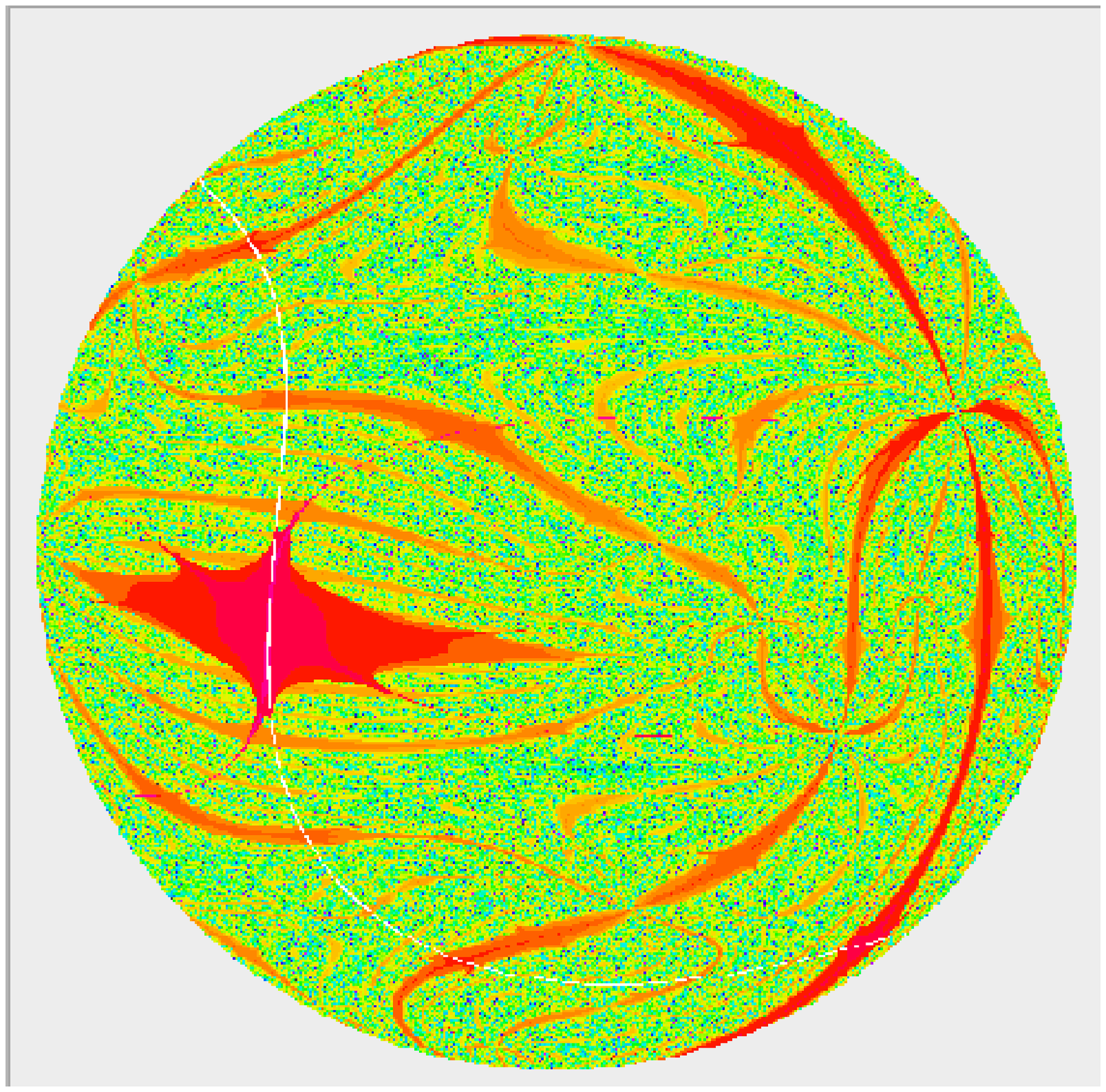,height=3in}\qquad
\psfig{figure=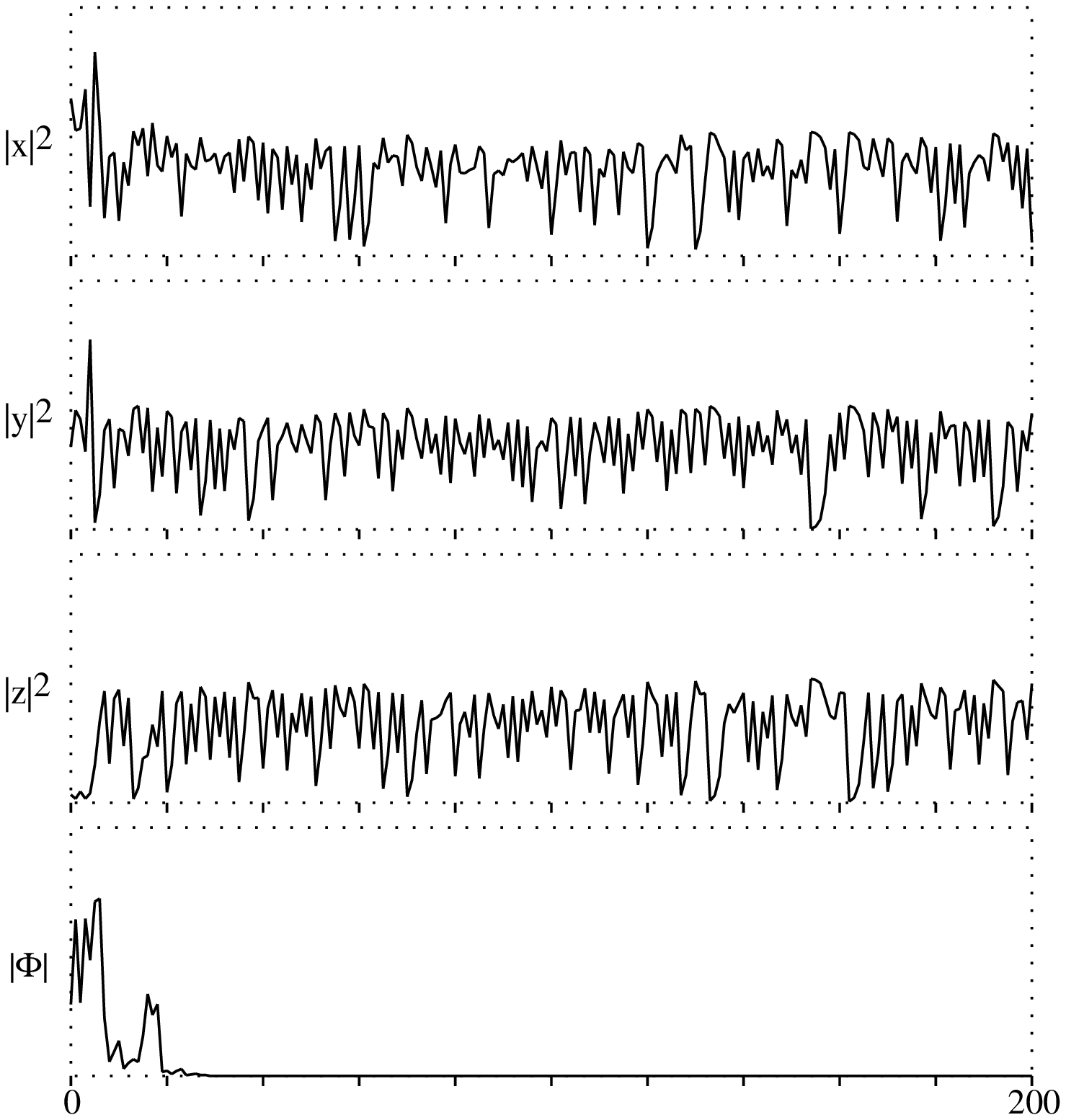,height=3in}}
\mycaption{\bit\label{fig:+-} $($See Example \ref{ex:+-}.$)$
For the map with Desboves parameters $(\frac{1}{3}\,,\,0\,,\,-\frac{1}{3})$,
the Fermat curve
again seems to attract all or nearly all orbits in both the real and
complex cases.}
\end{figure}

{\bf Note on Pictorial Conventions.}\; Each of the color pictures
which follow shows   the
real projective plane, represented as a unit 2-sphere with antipodal points
identified, oriented as in Fig.~\ref{fg:rpp}. Thus the $x$-axis, pointing
towards the right,
and the $y$-axis, pointing almost vertically, are close to the
plane of the paper, while the $z$-axis points up out of the paper.
(Because of this choice of orientation, we will sometimes refer to the
coordinate point $(0:1:0)$ near the top of the picture as the {\bit north
pole\/}.)
The Fermat curve $x^3+y^3+z^3=0$ is traced out in white. In Figs.
\ref{fig:min} and \ref{fig:+-},
other points are colored from red to blue according as their orbits converge
more rapidly or more slowly towards this Fermat curve, and subsequent figures
use various modifications of this scheme.
As an example, in Fig. \ref{fig:min} the {\bit equator\/} $y=0$ shows up as a
blue circle, since orbits in the invariant line $y=0$ cannot converge to $\F$,
hence orbits near this line cannot converge rapidly towards~$\,\,\F$.


\begin{figure}
\centerline{\psfig{figure=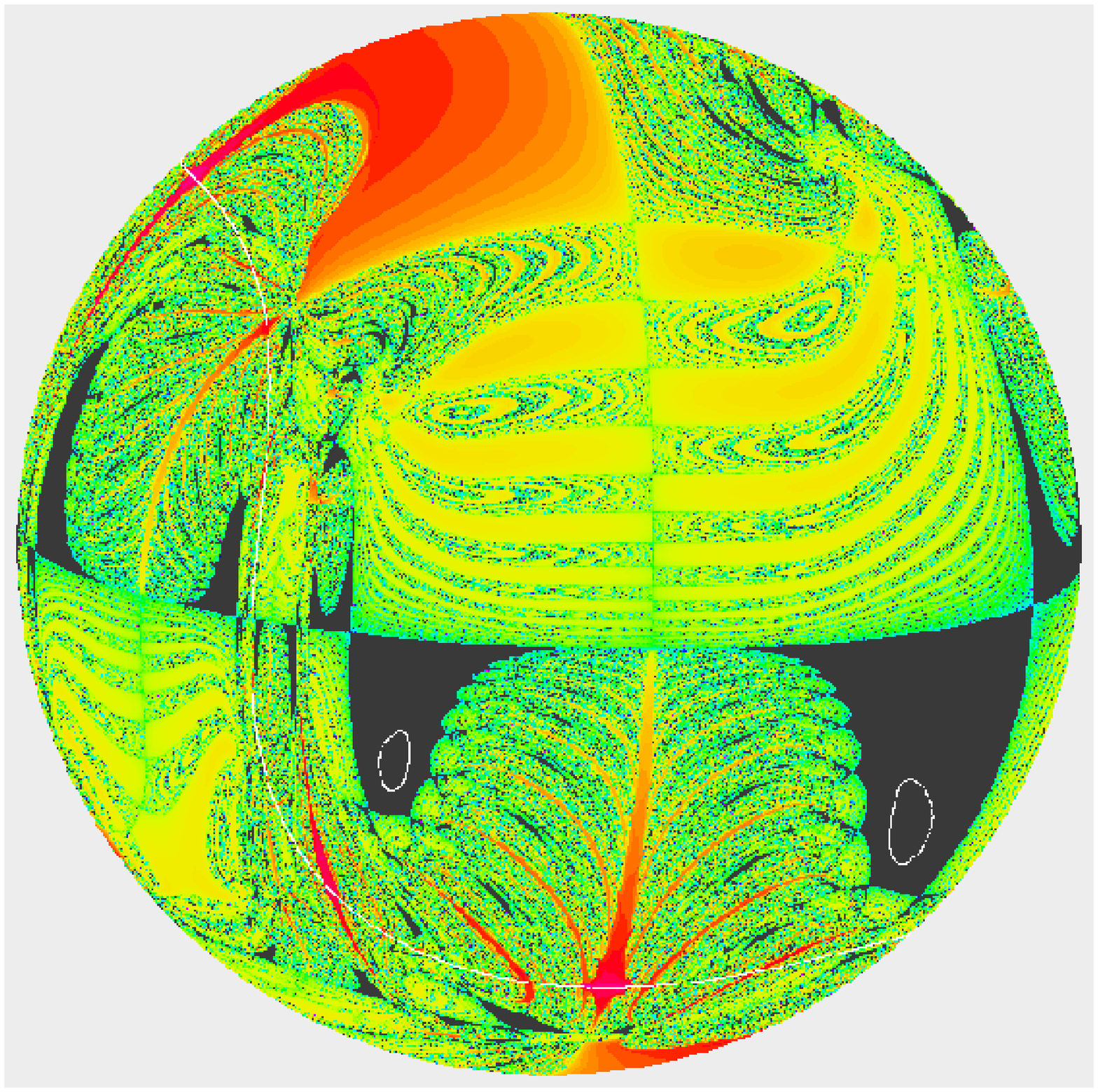,height=3in}\quad
\psfig{figure=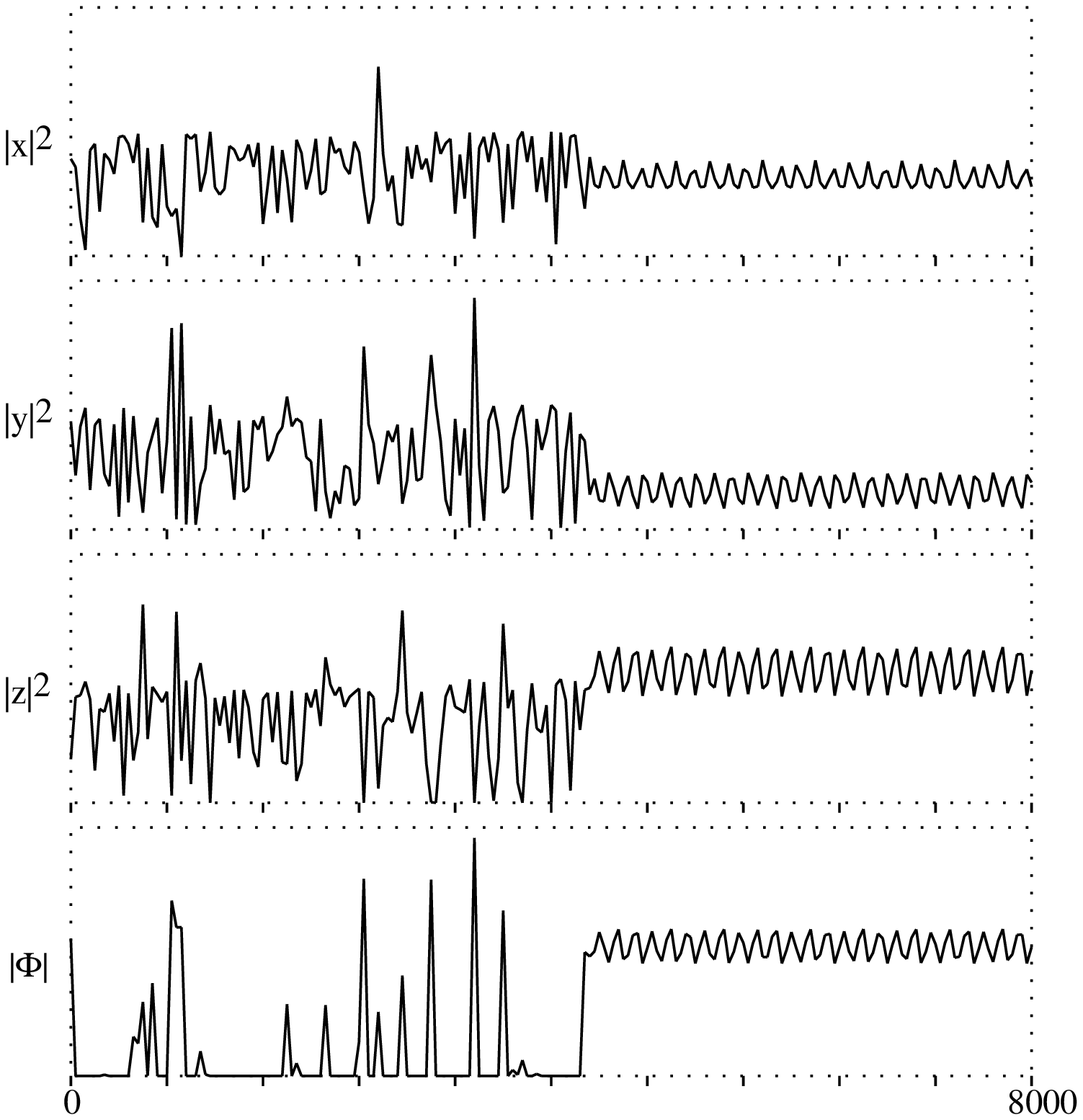,height=3in}}
\mycaption{ \bit\label{e13}  $($Example \ref{e4}.$)$
Dynamics for the parameters $(-\frac{1}{5}\,,\, \frac{7}{15}\,,\,
\frac{17}{15})$. Left: In the real case there are
two attractors. The basin of the Fermat curve is
colored as in Figs. \ref{fig:min}, \ref{fig:+-}. However, the two small
white circles also
form an attractor. The corresponding basin is shown in dark grey. Right:
A typical randomly chosen orbit for the complex map. This orbit often
comes very close to the
Fermat curve during the first 4000 iterations, but then seems to converge
to a cycle of two Herman rings.}
\bigskip \bigskip \bigskip

\vspace{.5cm}

\centerline{\psfig{figure=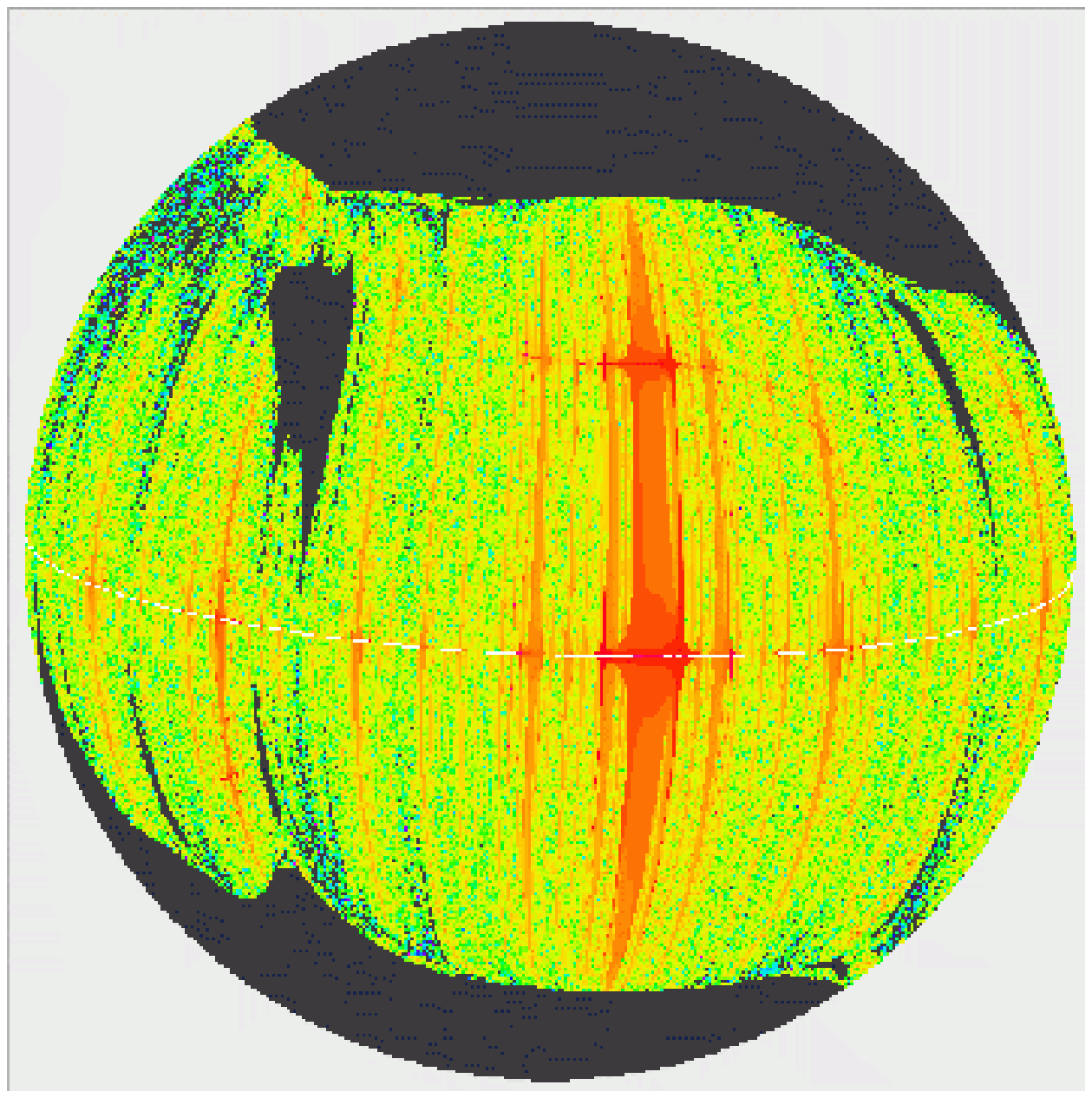,height=3in}\qquad
\psfig{figure=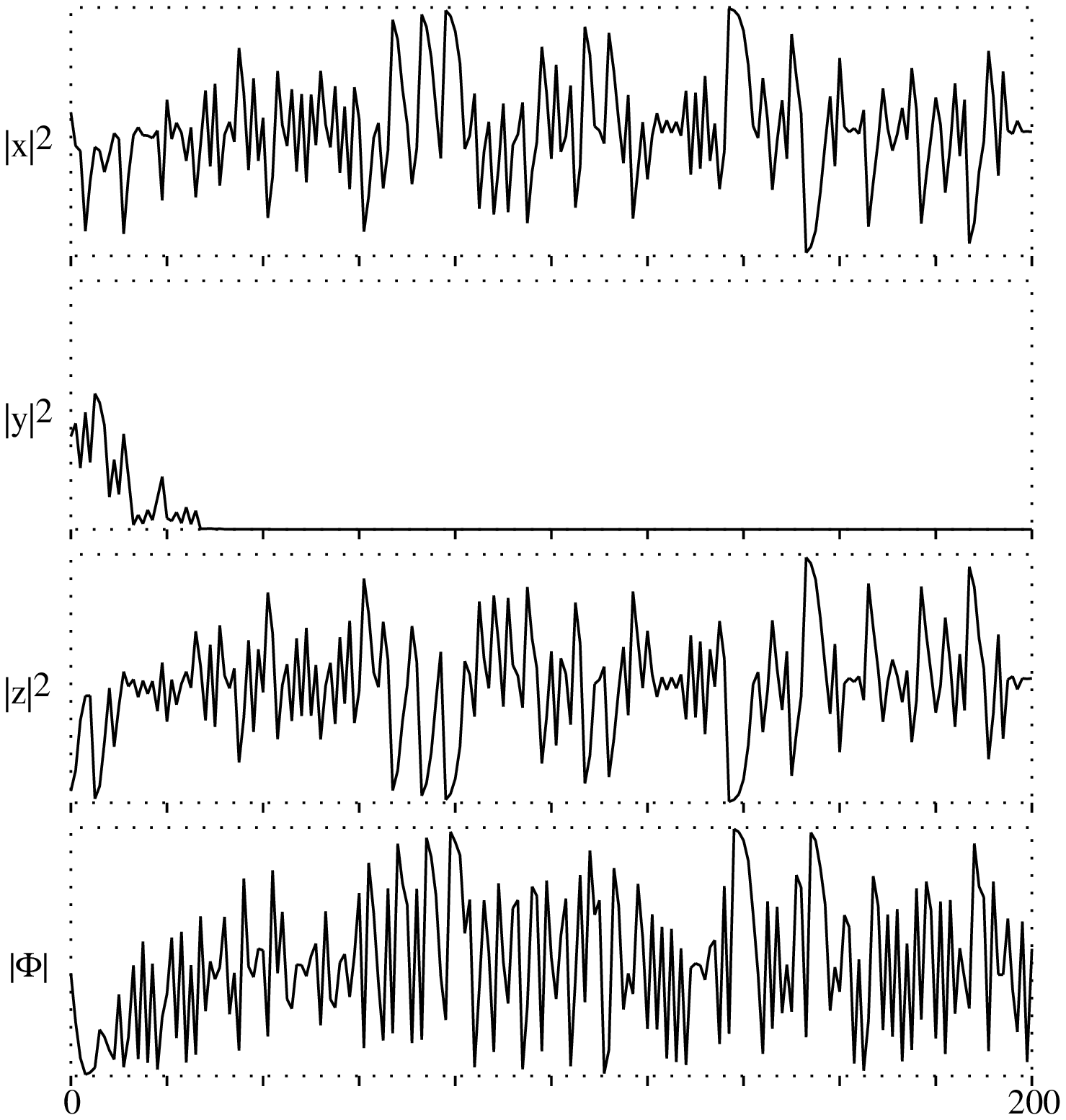,height=3in}}
\mycaption{ \bit \label{e18} $($Example \ref{e19}$\,.)$
Plots for the map with Desboves parameters $\;(-1.4\,,\,-.8\,,\,1.4)\,$.
Here the coloring is as in the previous
figures except that it describes convergence to the
``equator'' $y=0$,  rather than to the invariant Fermat curve. For this map,
the ``north pole'' $\;(0:1:0)\;$ also attracts many orbits.}
\end{figure} 

\begin{figure}
\centerline{\psfig{figure=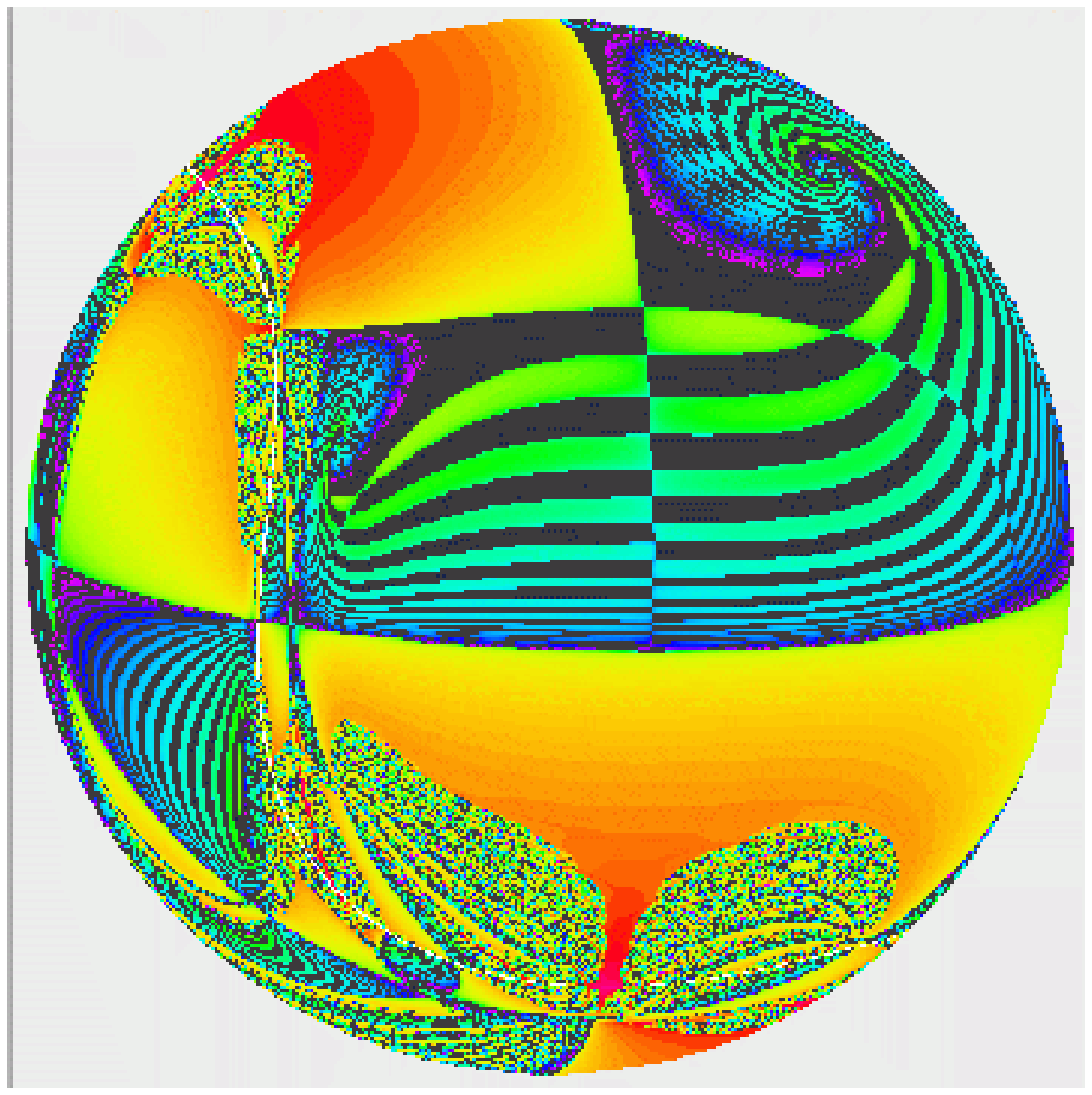,height=3in}
\psfig{figure=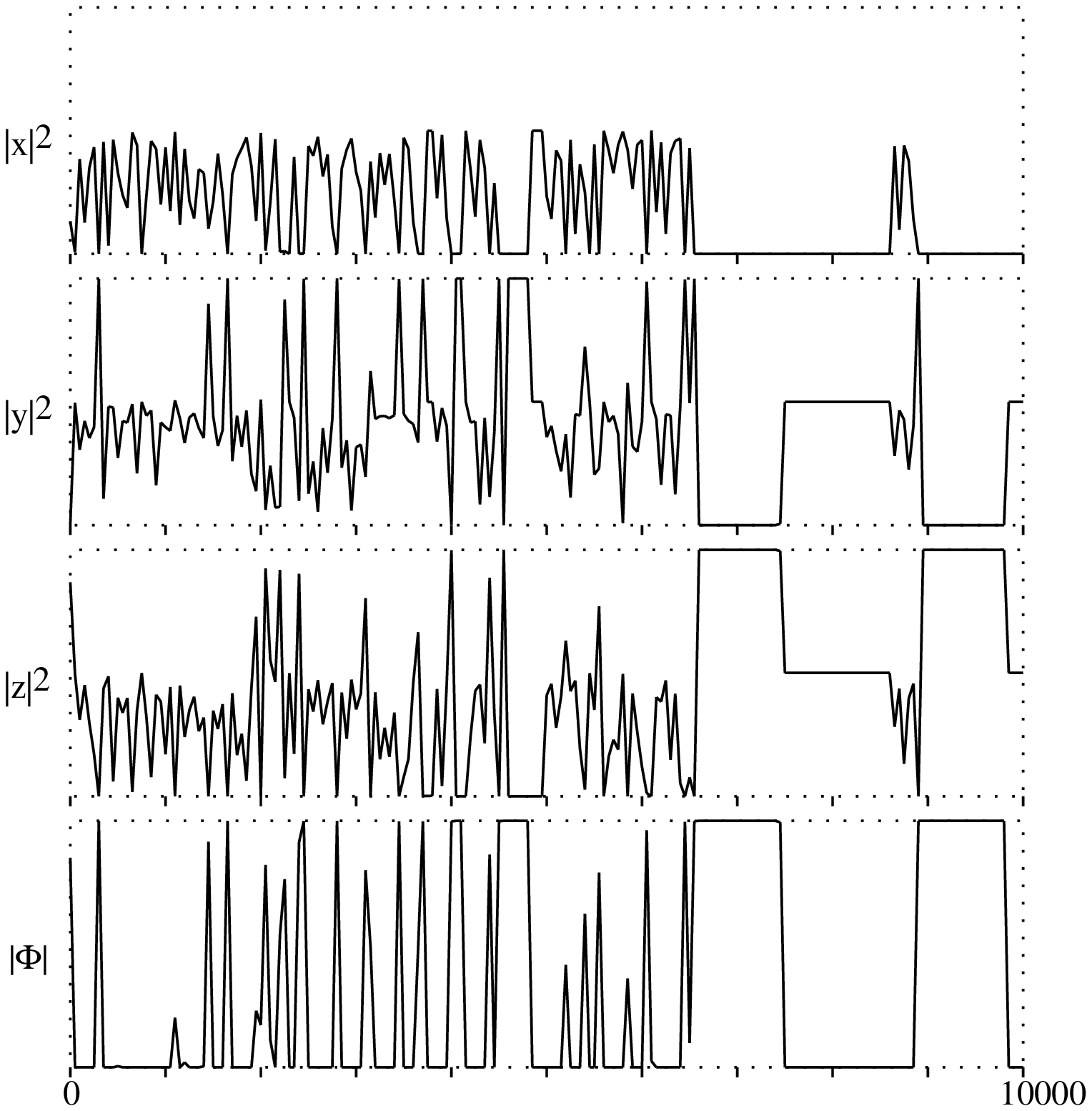,height=3in}}
\mycaption{ \bit\label{e24} $($Example \ref{e8}.$)$
On the left: Corresponding figure for the real
Desboves map with parameters $(\frac{1}{3}\,,\,1\,,\,\frac{5}{3})$,
again describing convergence
to $($or at least coming close to$)$ the Fermat curve. On the right:
One randomly chosen orbit for the complex map through $10000$ iterations.}
\bigskip \bigskip \bigskip

\centerline{\psfig{figure=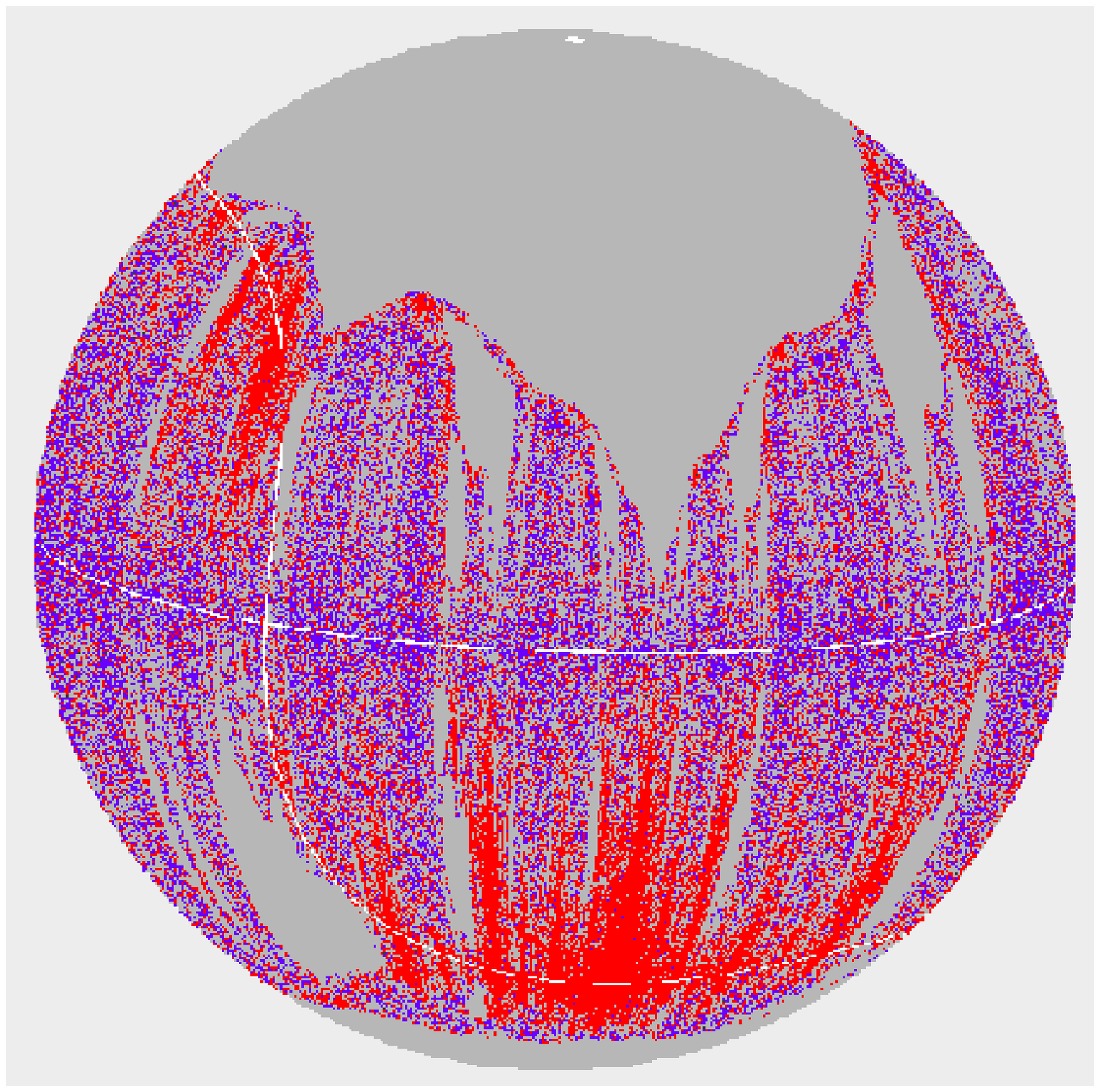,height=3in}\qquad
\psfig{figure=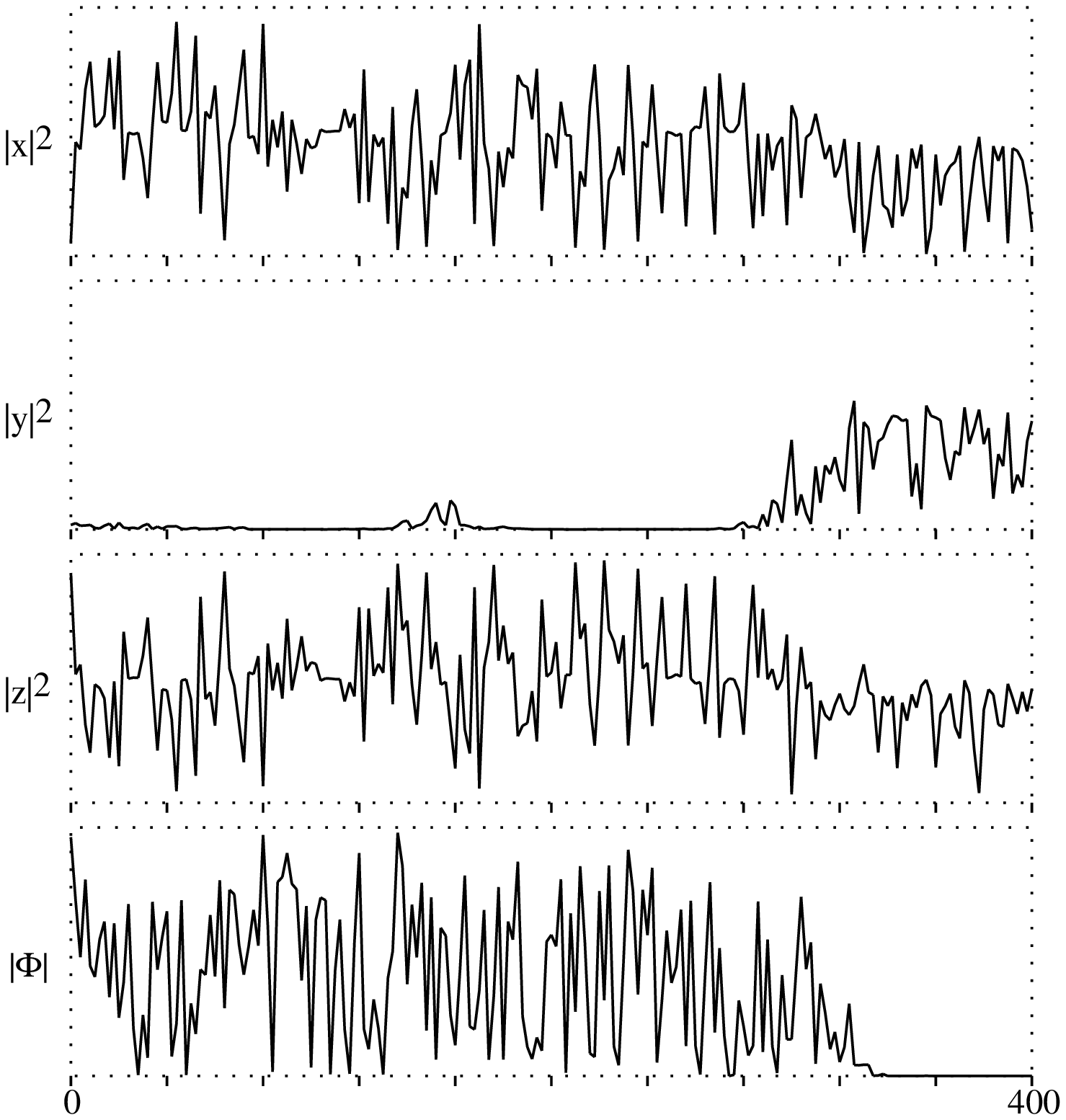,height=3in}}
\mycaption{ \bit \label{e14}  $($Example \ref{3bas}.$)$
Plots for the ``elementary map'' with parameters $(a,b,c)=
(-1\,,\,\frac{1}{3}\,,\,1)$. 
In this case, every great circle through the north pole $(0:1:0)$
maps to a great circle through the north pole.
There are three attractors: the Fermat curve $\;\F\,$, the equator $\;\{y=0\}\,$,
and the north pole, each marked in white. The corresponding attracting basins
are colored red, blue, and grey respectively.
$($However, the closely intermingled blue and red yield a purple
effect.$)$ The graphs on the right show an orbit which nearly converges to
$\;\{y=0\}\;$ but then escapes towards $\;\F\,$.
}
\end{figure}

The graphs to the right of Figs. \ref{fig:min}, \ref{fig:+-}, \ref{e13}, \ref{e18}, \ref{e24},
\ref{e14} illustrate some more or less typical randomly
chosen orbit for the associated complex map. Here each orbit point
 $\,(x:y:z)\,$ has been normalized so that $\;|x|^2+|y|^2+|z|^2=1\,$.
The horizontal
coordinate measures the number of iterations, while the vertical coordinates in
each of the four stacked graphs represent respectively
 $\;|x|^2\,,~|y|^2\,,~|z|^2\,$, and $\;|\Phi(x,y,z)|\,$.

\begin{ex}{\bf The Fermat Curve as a Global Attractor?}\label{ex:1}
If we choose
Desboves parameters $(b-\frac{2}{3}\,,\,b\,,\,b+\frac{2}{3})$ with $|b|$
small, then the transverse Lyapunov exponent is negative in both the
real and complex cases. Numerical computation suggests that nearly all
orbits actually converge to the Fermat curve. (Perhaps even all but a set
of measure zero?)
As an example, consider the case $(a,b,c)=(-\frac{2}{3}\,,\,0\,,\,\frac{2}{3})$.
Using the Gnu multiple precision arithmetic package,
and starting with several thousand randomly chosen points on the
real or complex projective plane, one can check that all orbits land on
the curve, to the specified accuracy, within a few hundred iterations.
Of course, even if we could work with infinite precision arithmetic, such a
computation could not prove that a given orbit converges to the curve, and also
could not rule out the possibility of other attractors with extremely small
basins. In fact it seems possible that periodic attractors with high
period and small basin exist for a dense open set of parameter values. This
case $b=0$ is rather special in one way, since the map $f_{-2/3\,,\,0\,,\,2/3}$
has points of indeterminacy; namely those points where $(x^3:y^3:z^3)$ is equal
to either $(1:7:1)$ or $(0:1:0)$.
However, the
behavior for small non-zero values of $b$ seems qualitatively similar.
\end{ex}

\begin{ex}{\bf An Even Stronger Attractor.}\label{ex:2}
The case $b=\pm\frac{1}{9}$ yields a even more strongly attracting Fermat
curve, as illustrated in Fig. \ref{fig:min}. The transverse
derivative has a simple zero at the point $(-1,1,0)/\sqrt 2$ to the upper
left of the figure, and a double zero at the point $(0,-1,1)/\sqrt 2$ near
the bottom. A numerical search suggests that this is the most attracting
example within the real or complex Desboves family,
in the sense that the transverse exponent takes its most negative
value of $-2.0404\ldots\,$ for the real map or $-0.6801\ldots\,$ for the
complex map. Certainly these are the extreme values for real parameters
within the two-thirds family, as graphed in Fig. \ref{fig:2/3}.
\end{ex}

\begin{ex}{\bf Another Global Attractor?}\label{ex:+-}
 If we take Desboves coordinates  $(\frac{1}{3}\,,\,0\,,\,-\frac{1}{3})$, then again the Fermat
curve seems to attract nearly all orbits. Compare Fig. \ref{fig:+-}.
Here the transverse derivative has a double zero at the fixed point $(-1:0:1)$
in the middle of the large red region.
It is a curious fact that the transverse exponents in this case are precisely
the same as those for Example \ref{ex:1}, namely $-1.456\cdots$ for the
real map, or $-.549\ldots$ for the complex map.
\end{ex}

\begin{ex}{\bf A Cycle of Herman Rings?}\label{e4}
Now suppose that we choose Desboves parameters in the two-thirds
family, with $(a,b,c)$ equal to 
$\,(-\frac{1}{5}\,,\,\frac{7}{15}\,,\,\frac{17}{15})$. Here the
transverse exponent is $-.509\cdots$ for the real map, but $\;+.402\cdots$
for the complex map. Thus
we can expect the
Fermat curve to be an attractor in the real case, but not in the complex case.
The left half of Fig. \ref{e13} illustrates the dynamics in the real case.
Numerical computation suggests that some 83\% of the orbits
converge to the Fermat curve, while the remaining 17\% converge to
a pair of small circle. The attractive basin for this pair of circles
is conjecturally a dense open subset of ${\mathbb P}^2({\mathbb R})$.
The map $f=f_{a,b,c}$ carries each of these circles
to the other, reversing orientation, while $f\circ f$ carries each circle
to itself with rotation number $\pm 0.18587\cdots$.
Of course such a phenomenon can be expected to be highly sensitive
to small changes in the parameters---We cannot really distinguish
between a rotation
circle with irrational rotation number and one with a rational rotation
number which has very large denominator (although the later would
necessarily contain a periodic orbit).

\begin{figure} [t]
\centerline{\psfig{figure=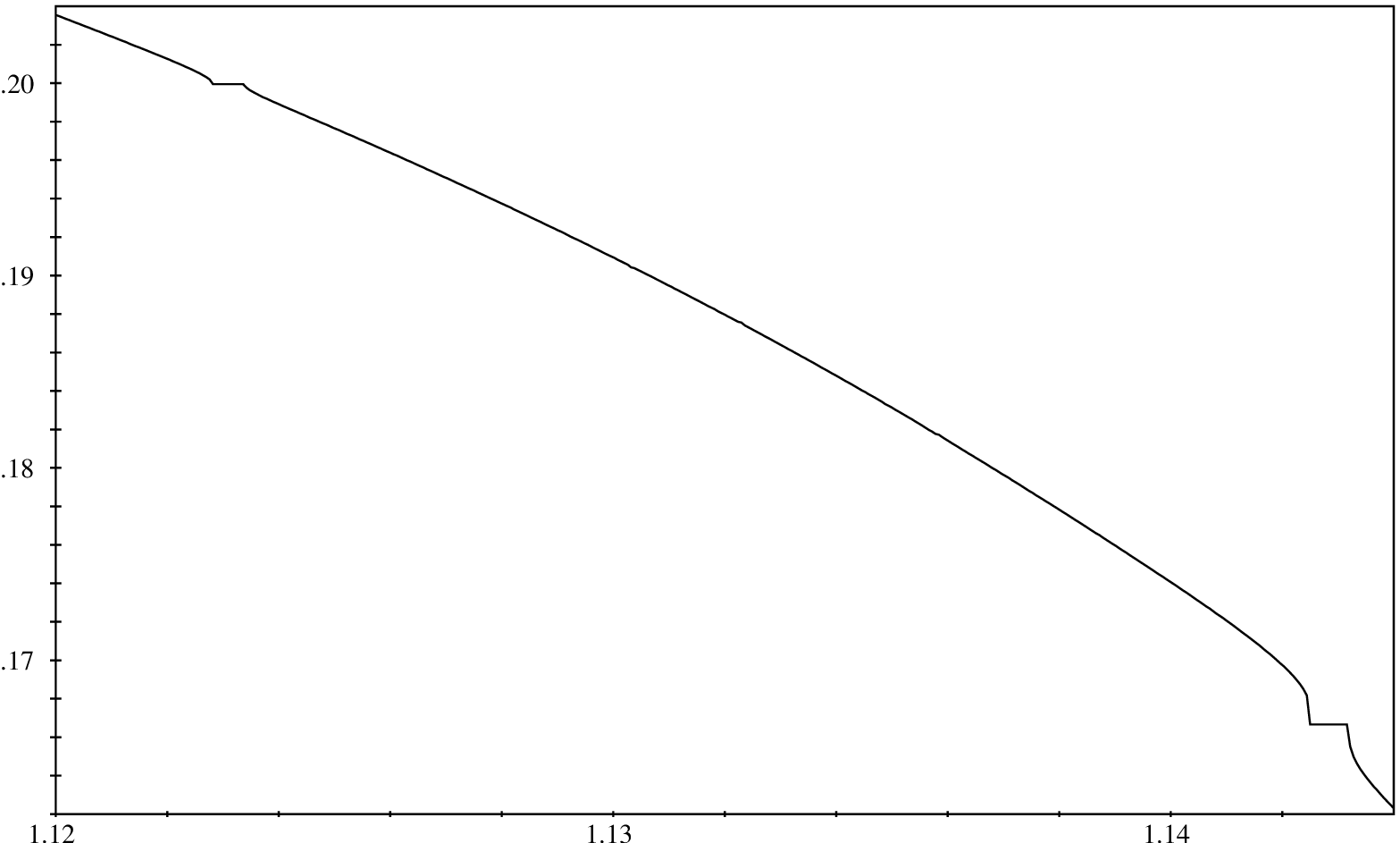,height=3in}}
\mycaption{\label{fig:HRP}\bit
An empirical plot of the rotation number for the pair of attracting
circles in $\bP^2(\R)$ as a function of the parameter $c$
in Example {\rm\ref{e4}}, keeping $a$ and $b$ fixed.
Presumably for each rational value for the rotation number there
corresponds an entire plateau of $c$ values for which the pair of circles
contains an attracting periodic orbit. Only the plateaus of height
$\frac{1}{5}$ and $\frac{1}{6}$ are visible in this figure; but with higher
resolution, tiny blips at height $\frac{3}{16}$ and $\frac{2}{11}$
would also be visible. As $c$ decreases past 1.12 the attracting circles shrink
to points; while as $c$ increases past 1.144 they expand until they break
up upon hitting the boundary of their attracting basin.
It is conjectured that whenever the rotation number is Diophantine, the
corresponding pair of circles in $\bP^2(\R)$ are contained in a pair of Herman
rings in $\bP^2(\C)$.}
\end{figure}

In the complex case, the Fermat curve is no longer an attractor. In fact,
almost all orbits seem to eventually land near this cycle
of circles and then to behave just like an orbit on a pair of nearby circles
with the same rotation number. {\it
This suggests that most orbits converge to
a cycle of two Herman rings in ${\mathbb P}^2({\mathbb C})$, with the pair
of real circles as their central circles.\/}
(For a more detailed discussion of Herman rings in $\bP^2(\C)$, see
\S\ref{s:hr}.) Again we must be cautious, since
such a phenomenon must be highly sensitive to perturbations; but
the empirical evidence certainly suggests the existence of a cycle of two
Fatou components which could only be the immediate basins for attracting
Herman rings.\footnote{For the description of possible Fatou components
in $\bP^2$, compare \cite{FS1}, \cite{U1}.}
(The convergence is very slow, and there may be other much more chaotic
attractors.) These attracting circles persist under small perturbation
of the parameters. (Compare Theorem \ref{th:nearhr}.)
A plot of the rotation number for these circles
as a function of the parameter $c$, keeping $a$
and $b$ fixed, is shown in Fig. \ref{fig:HRP}. It seems empirically
that Herman rings also persist under perturbation, but we have been unable
to justify this statement theoretically.
\end{ex}


\begin{figure}
\centerline{\psfig{figure=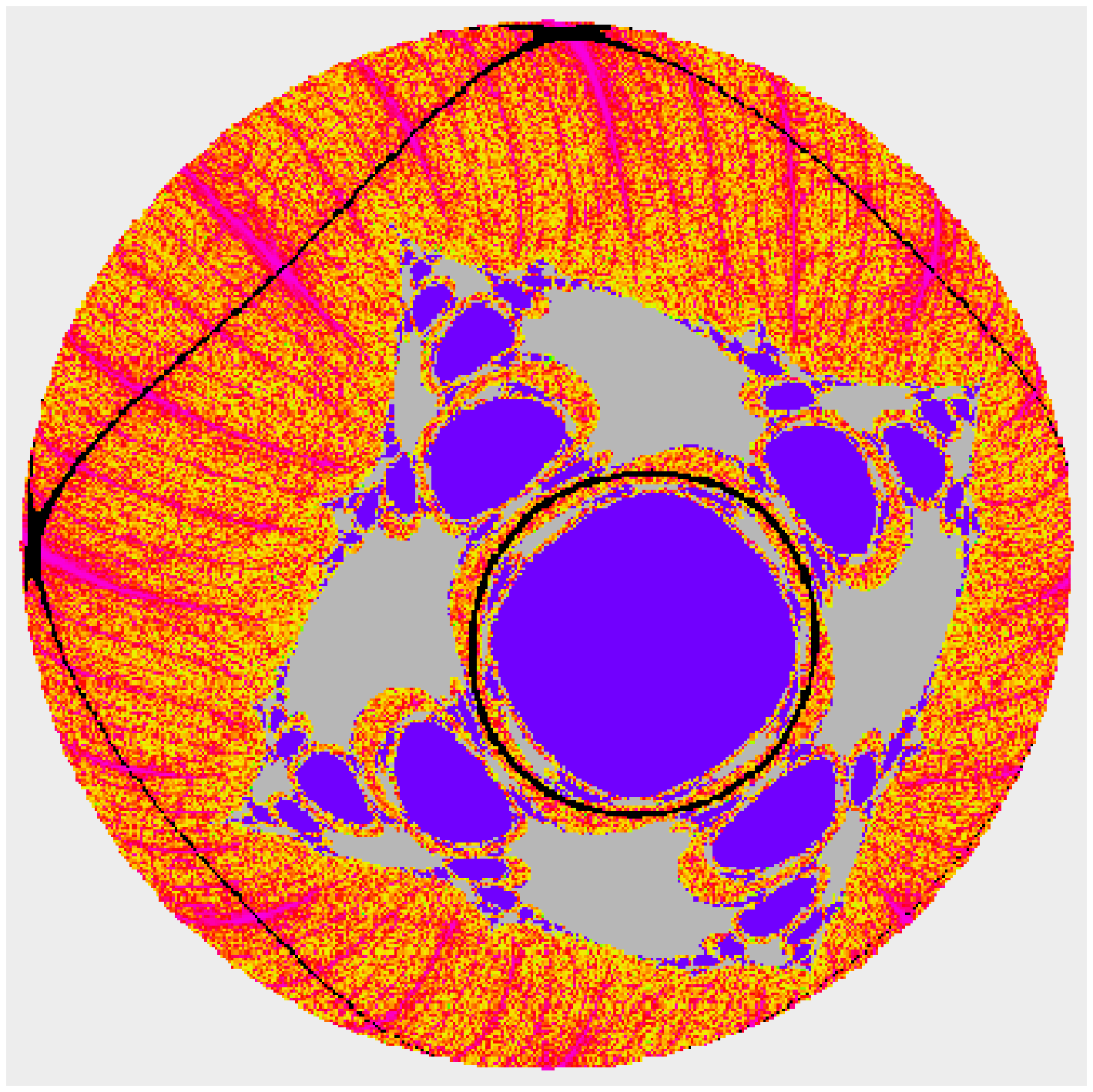,height=3in}}
\mycaption{ \bit \label{e20}$($Example \ref{ex:cas}.$)$
The Cassini quartic with parameter $k=\frac{1}{8}$,
shown in black, consists of an outer circle ${\mathcal C}_{\mathbb R}^0$
with two self-intersections
and a much smaller inner circle ${\mathcal C}_{\mathbb R}^1$.
Here the warmer colors describe convergence towards $\;\cC_\R^0\;$
for the rational map $f_1$ with parameter $a=1$.
The blue region is the basin of
a superattracting fixed point at $(0:0:1)$, while the grey region is the
basin of another attracting fixed point directly above it at $(0:1:2)$.}
\bigskip \bigskip \bigskip

\centerline{\psfig{figure=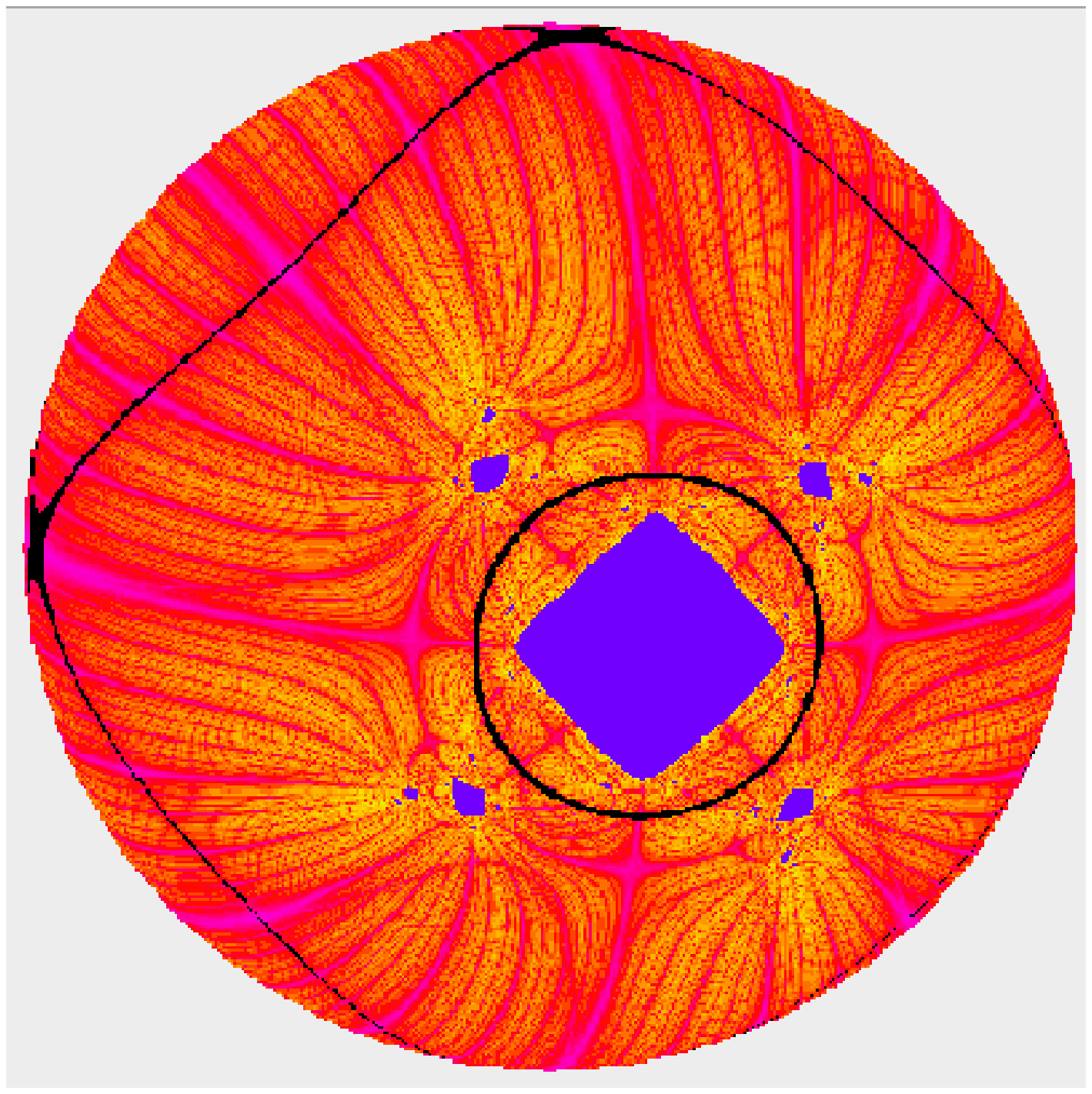,height=3in}}
\mycaption{ \bit \label{e21}
 Corresponding picture for the same Cassini curve
with $k=\frac{1}{8}$, but using the map with parameter $a=\frac{2}{5}$.}
\end{figure}

\begin{ex}{\bf The Line $z=0$ as a Measure Theoretic Attractor?}\label{e19}
(Compare Fig. \ref{e18}.)
For the parameter values $(a,b,c)=(-1.4\,,\,-.8\,,\,1.4)$,
 the Lyapunov exponent turns out to be strictly positive, equal to
$0.247\cdots$ in the real
case, or to $0.352\cdots$ in the complex case. The invariant Fermat curve
does not seem to play any significant dynamical role in this case.
On the other hand, the equator $y=0$ seems to be at least a measure-theoretic
attractor; and there is also an attracting fixed point at the north pole
$(0:1:0)$. In fact many randomly chosen real or complex
orbits converge to the north pole $(0:1:0)$, but even more seem to
converge to the equator.
\end{ex}


\begin{figure} [t]
\centerline{\psfig{figure=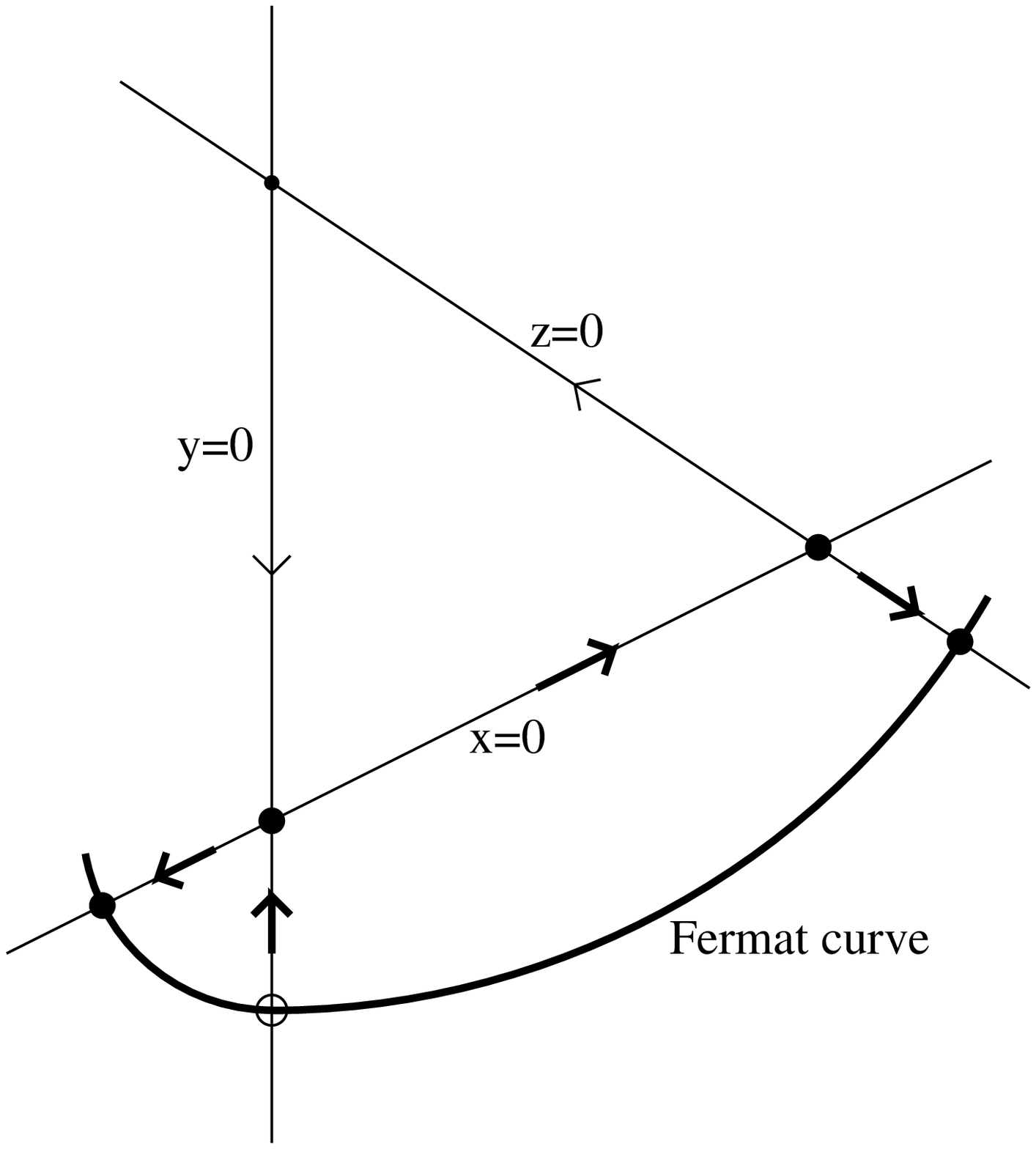,height=3in}}
\mycaption{\bit \label{fig:4.7}
  Schematic diagram illustrating Example \ref{e8}.\hskip 2.5in}
\end{figure}

\begin{ex}{\bf A Composite ``Almost'' Attractor.}\label{e8}
 For the real or complex map with
Desboves parameters $\;(\frac{1}{3}\,,\,1\,,\,\frac{5}{3})\,$
as illustrated in Fig. \ref{e24},
typical orbits seem to spend a great deal of time quite close to the
Fermat curve $\;\F\,$, 
even though the
transverse exponent is strictly positive, equal to
$\;0.081\cdots\;$ in the real case or to $\;1.032\cdots\;$ in the complex case.
This curve is not an attractor by itself, since nearby
orbits eventually get kicked away from it. However, the union
\begin{equation}\label{eq:union}
    A~=~\{x=0\}\cup\{y=0\}\cup\{z=0\}\cup{\cal F}\,,
\end{equation}
or in other words the variety $~x\,y\,z\,\Phi(x,y,z)=0\,$,~
does seem to behave like an attractor, at least in a statistical sense.
(Compare \cite{GI}.)
Typical orbits seem to spend {\it most\/} of the time extremely close to
this variety. However, they do not stay in any one of its four
irreducible components, but sometimes jump from one component to another.
Furthermore, it seems likely that typical orbits will escape completely from a
neighborhood of this variety, very infrequently but infinitely often.

Here is a more detailed description, as illustrated in Fig.~\ref{fig:4.7}.
To fix ideas, we will refer to the real case; but the complex case is not
essentially different.  A randomly chosen orbit seems to spend most of the
time either
wandering chaotically very close to the Fermat curve or else almost stationary,
very close to one of the four saddle fixed points (the black dots in
Fig.~\ref{fig:4.7}).  However, such an orbit does not seem to stay close to any one of
the four components of this variety forever. For example, it is likely to
escape from the neighborhood of the Fermat curve $\;\F\;$ when
it comes very close to the strongly repelling point $\;{\cal F}\cap\{y=0\}\;$
which is circled in Fig.~\ref{fig:4.7}. It will then shadow the coordinate
line $\;y=0\,$, jumping quickly to a small neighborhood of the saddle point
$\;x=y=0\,$, and then slowly coming closer to this point for thousands of
iterates. Again it must eventually escape, now shadowing the
line $\;x=0\;$ and jumping quickly
either towards the saddle point $\;{\cal F}\cap\{x=0\}\;$
or towards the saddle point $\;x=z=0\,$.
In either case it again spends a long time approaching this saddle point, but
then escapes. In the first case, it is now very close to the Fermat curve
and shadows it for a long time with a highly chaotic orbit before
starting the cycle again. In the second case where it escapes near the saddle
point $\;x=z=0\,$, it then shadows the line $\;z=0\,$ as it
quickly converges towards the saddle
point $\;{\cal F}\cap\{z=0\}\,$, where it again remains for a long time before
repeating the cycle.
\end{ex}

\begin{rem}{\bf The Concept of Attractor.}\label{rm:attr}
Such examples have led authors such as \cite{GI}
and \cite{AAN} to suggest modified definitions of attractor, emphasizing not
the omega-limit set of a typical orbit, but rather its asymptotic probability
distribution. To illustrate the effect of such a change, think of a dynamical
system in the plane
in which orbits spiral out towards a limit cycle $\Gamma$ which consists of
a homoclinic loop, begining and ending at a fixed point $p$. Then the
unique ``measure theoretic attracting set'' for the region inside the loop
is the entire loop $\Gamma$. No orbit starting inside actually converges to
the point $p$. However,
every orbit starting inside the loop spends {\it most\/} of its time
apparently converging to $p$, with a statistically insignificant
(but infinite) collection of exceptional times. Thus, from the point of
view of Gorodetski and Ilyashenko or of Ashwin, Aston and Nicol, the
``attractor'' is the single point $p$.
\end{rem}

\begin{figure}
\centerline{\psfig{figure=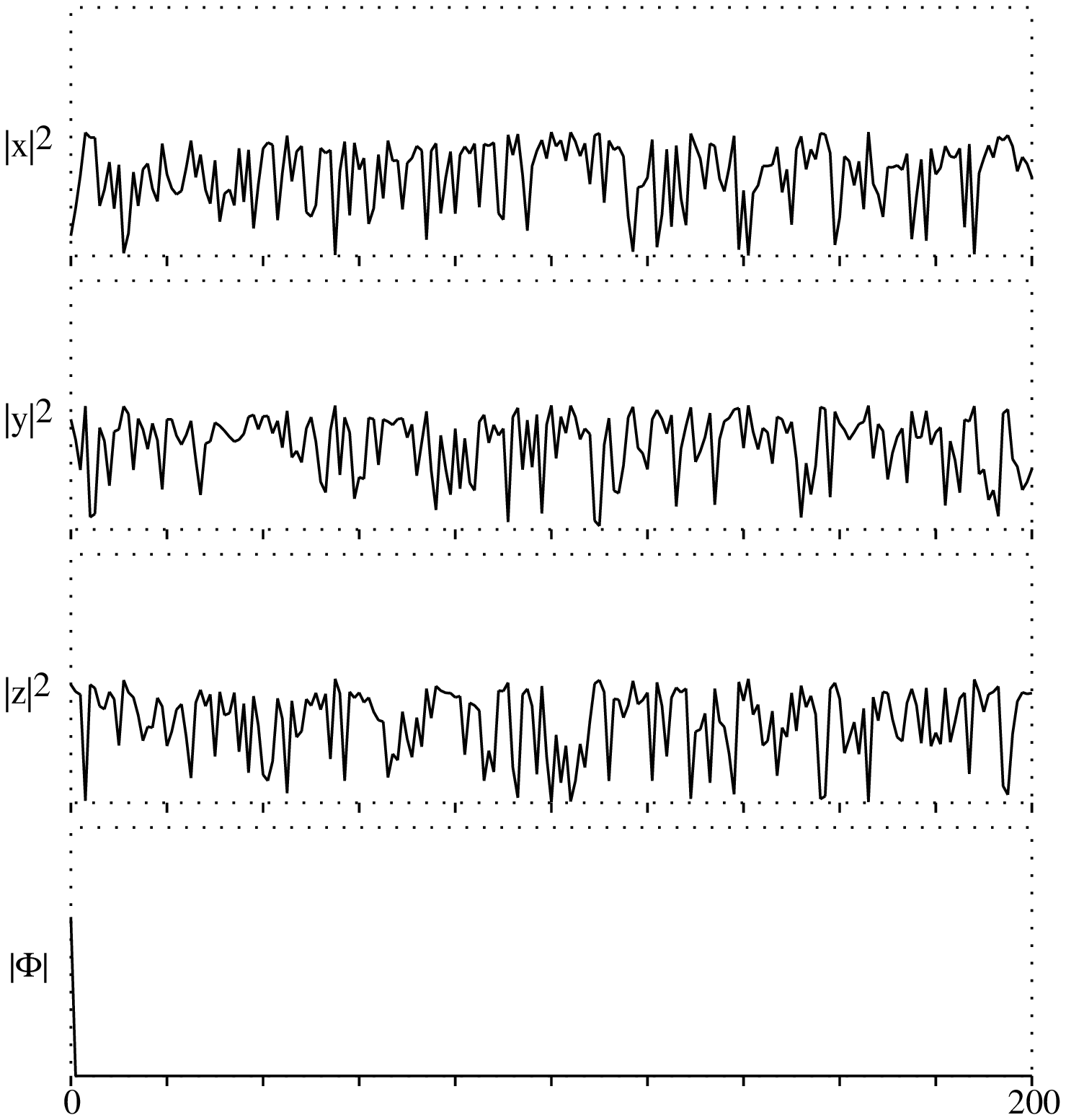,height=3in}}
\mycaption{\bit \label{fig:all-3map} \bit A plot of 200 iterations with
a random start for the degree three map of Example \ref{dg3map}, taking
parameters $\;(a,b,c)=(0,\gamma,-1)/2$,}
\end{figure}

\begin{ex}{\bf A Family of Degree Three Maps.}\label{dg3map}
In the case of a map of degree three  
(or indeed for any
degree which is not a perfect square), the multiplier
on an invariant elliptic curve cannot be a real number. Hence
we cannot describe both the map and the invariant curve by equations with
real coefficients---we can only consider the complex case.

Here is one explicit family of degree three maps which send
the Fermat curve $\F$ into itself.
As in \S\ref{s:De}, we start with a rather degenerate self-map of $\bP^2$.
Let
$$H_0(x, y, z)\;=\; \big(kxyz\,,~~ y^3-\gamma z^3\,,~~ z^3-\gamma y^3\big)$$
where $\gamma$ is the cube root of unity $(-1+i\sqrt 3)/2$
and where $k^3=3(\gamma^2-\gamma)$. (For example $k=\gamma^2-1$
or $\;k=i\sqrt 3$.) Then the associated map
$$h_0(x:y:z)~=
  ~\big(kxyz\; :\; y^3-\gamma z^3\; :\; z^3-\gamma y^3\big)$$
of projective space has a first integral $\Phi/\Psi$, where $\Phi(x,y,z)
=x^3+y^3+z^3$ and $\Psi(x,y,z)=x^3$. In fact a brief computation shows that
$$ {\Phi(H_0(x,y,z))/ \Phi(x,y,z)}~=~
 {\Psi(H_0(x,y,z))/ \Psi(x,y,z)}~=~ 3\,y^3z^3\,,$$
and it follows immediately that the rational function $\Phi/\Psi$ is invariant
under $h_0$. In particular, the Fermat curve, defined by $\Phi=0$, is
$h_0$-invariant. (However, in contrast to the Desboves case of \S\ref{s:De},
the various elliptic curves $~\Phi/\Psi={\rm constant}\in\C\ssm\{1\}~$
are all mutually isomorphic. A similar example will be described in Remark
\ref{rem:degen}.) This map $h_0$ has just one point of indeterminacy,
namely $(1:0:0)$, which is not on ${\mathcal F}$.

Like all maps of $\bP^2$ with first integral, $h_0$ is
not very interesting as a dynamical system, but it does embed in a family of
more interesting maps.
Consider the 3-parameter family of homogeneous polynomials
$$  H=H_{a,b,c}~=~H_0+(a,b,c)\,\Phi~. $$
Each of the associated maps $h_{a,b,c}$ of the projective plane carries
the Fermat curve $\mathcal F$ to itself with degree three and
multiplier $\gamma k$.
There is only one fixed point of $h_{a,b,c}$ on $\F$, namely $(0:-1:1)$.

 First consider the 1-parameter subfamily
of  maps satisfying the conditions that $a=b+\gamma c=0$,
with $x=0$ as invariant line.
(The use of these special parameters simplifies the computation of the
transverse Lyapunov exponent.) 
When $b=\gamma/2$, the transverse Lyapunov
exponent takes its most negative value of -1.647918. Thus $\mathcal F$ appears
to be more strongly attracting under this map than under any of the complex
Desboves maps, where the most negative transverse exponent was
 $-0.6801\cdots$. (See Example 5.2. Such computations will be
explained in Part 2 of this paper.)
In Fig.~\ref{fig:all-3map}, we show the extraordinary attracting properties
of the Fermat curve for this map. Most randomly chosen points seemed to hit
the Fermat curve, up to the resolution of the graph, after only 6 or so
iterations. (Compare with the right hand sides of Figs.: \ref{fig:min},
\ref{fig:+-} and \ref{e14}.)

It is also interesting to consider the subfamily consisting of maps
$h_{a,\,0,\,0}$ with $b=c=0$. These are {\bit elementary maps~}, as described
in  Example \ref{e7}. For this particular
elementary family, it appears that the transverse exponent is
always non-negative, so that the Fermat curve is never an attractor.
\end{ex}

\begin{ex}{\bf The Degree Two Case.}\label{rem:deg2}
According to \cite[Proposition 6.6]{BD}, up to holomorphic conjugacy,
there are exactly 20 distinct examples of holomorphic self-maps of $\bP^2(\C)$
of algebraic degree two with an invariant
smooth elliptic curve. See  Example 6.9 of their paper 
for a detailed study of one of these degree two maps. This example,
with multiplier equal to $i\sqrt 2$, has five
attracting cycles, with common basin boundary equal to the Julia set.
Four of these are attracting fixed points, and the fifth is an attracting
period 2 orbit. Empirically,  randomly chosen orbits for this example
always seem to converge to one of these five cycles.\end{ex}

\setcounter{lem}{0}
\section{Intermingled Basins}\label{s:rb}

\begin{ex}{\bf Elementary Maps.}\label{e7}
Finally we come to examples where we can provide complete proofs.
Following \cite[p.~16]{BD}, a rational map of $\bP^2$ is called {\bit
elementary\/}, with {\bit center\/} $p_0$, if every line through $p_0$
maps to a line through $p_0$. Elementary maps are easier to analyze than more
general rational maps since we can separate the variables to simplify
the discussion.

In particular, consider a Desboves map $f=f_{a,\,b,\,c}$ as in formula
(\ref{e11}) of \S\ref{s:De}, where the parameters $a,\,b,\,c$ satisfy
$a=-1$ and $c=1$.\,\, Then the image $f(x:y:z)=$\break$(x':y':z')$ satisfies
$$  x'~=~x(-x^3-2z^3)\qquad {\rm and}\qquad z'~=~z(2x^3+z^3)~.$$
It follows that each line $\;(x:z)={\rm constant}\;$ through the coordinate
point $\;p_0=(0:1:0)\;$ maps to another line
$\;(x':z')={\rm constant}'\;$ through $p_0$.
If we set $\;X=x/z\;$ and $\;X'=x'/z'\,$, then the correspondence
\begin{equation}\label{eq:Lat}
  \widehat f\;:\;X~~\mapsto~~X'\,=\,-X\frac{X^3+2}{2X^3+1}
\end{equation}
does not depend on the choice of parameter. This rational map  (\ref{eq:Lat})
is described as a {\bit Latt\`es map\/}, since it is the image of a rigid
map on the torus $\;\F\cong\C/\Omega\;$ under the semiconjugacy
$\;(x:y:z)\mapsto(x:z)\,$ of degree three. (In fact $\widehat f$ is conformally
conjugate to the Latt\`es map described in Remark \ref{rem:Lat}.)
It has an ergodic invariant measure
which is smooth except at its critical values, the cube roots of $-1$.
Over the real numbers, $\;\widehat f\;$ is a covering map from the circle
$\;\bP^1(\R)\;$ to itself with topological degree $\;-2\,$.
\end{ex}
\smallskip

\begin{figure}
\centerline{\psfig{figure=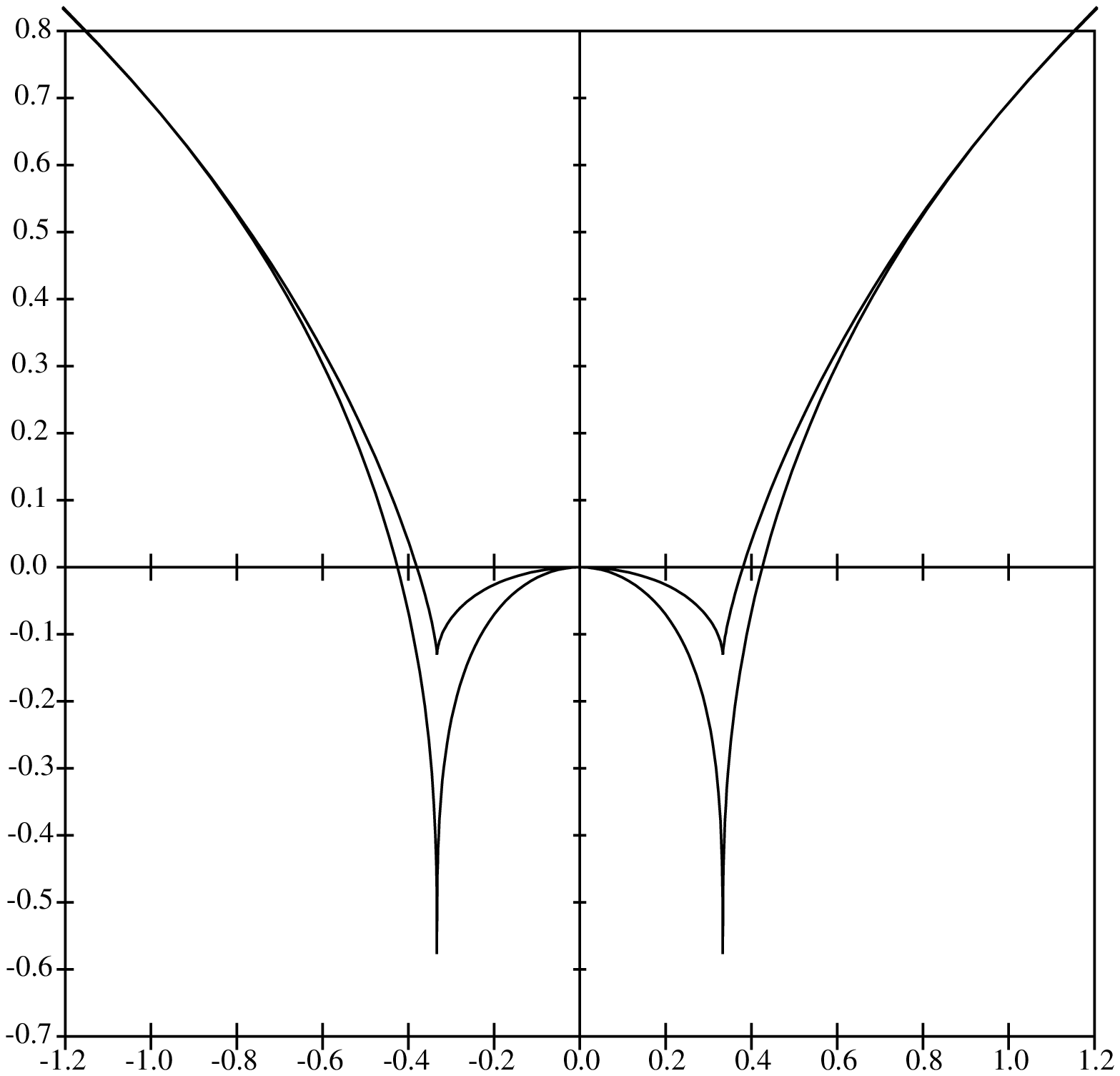,height=2.5in}\qquad
\psfig{figure=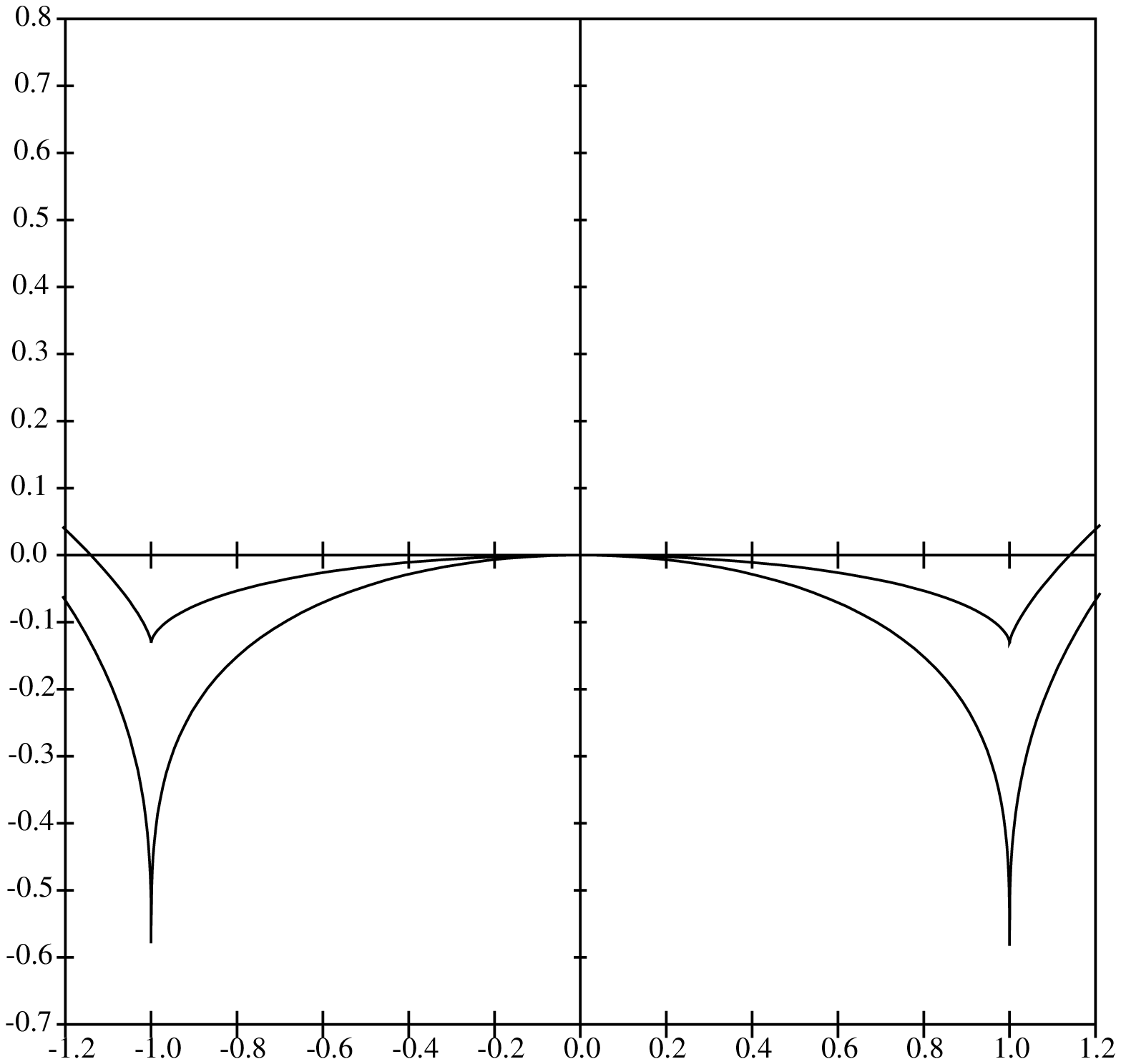,height=2.5in}}
\mycaption{\label{fig:elemg}\bit On the left:
Plot of the transverse Lyapunov exponent along the Fermat curve
as a function of the parameter $b$
for the elementary family of Example \ref{e7},
with Desboves parameters $\;(-1\,,\,b\,,\,1)\,$. On the right:
Corresponding plot for the transverse exponent along the line $\;y=0\,$.
In both cases, the graph for the complex map lies above the graph for the real
map.}
\end{figure}

Over either the real or complex numbers,
if we think of $\;{\mathbb P}^2\ssm\{p_0\}\;$ as a (real or complex)
line-bundle over the projective line $\;{\mathbb P}^1\;$ with
projection $\;\pi:(x:y:z)\mapsto(x:z)\;$, then we have the commutative diagram
\begin{equation}\label{eq:el-diagr}\begin{matrix}
    \;\;\bP^2\ssm\{p_0\} & \buildrel f\over\to& \;\;\bP^2\ssm\{p_0\}\cr
    \pi\downarrow~ &&\pi\downarrow~\cr
    \;\bP^1& \buildrel \widehat f\over\to& \;\bP^1
\end{matrix}\end{equation}
where $\;f\;$ carries each fiber into a
fiber by a polynomial map,
with coefficients which vary with the fiber. As an example, for the two
invariant fibers $x=0$ and $z=0$ we get the maps
$$(0:y:1)~\mapsto ~(0: by^4+(1+b)y:1)\qquad{\rm and}\qquad(-1:y:0)~\mapsto~
(-1: by^4+(1-b)y:0)$$

respectively. If we exclude
the degenerate case  $\;b=0\,$ (compare Remark \ref{rem:degen}), then these
polynomial maps all have degree four. Furthermore, the center point
$\;p_0=(0:1:0)\;$  is superattracting, and
serves as the point at infinity for each one.
In the real case, these polynomial maps are all unimodal, while in the complex
case they all have $120^\circ$ rotational symmetry.

Since the rational map $\;\widehat f\;$ of the base space has no
attracting cycles, it follows that {\it an elementary map
with invariant elliptic curve can have no
attracting cycles other than its center point.}

In this special case of an elementary map, we can give a relatively easy proof
that a negative transverse exponent for any invariant elliptic curve
implies that this curve
is a measure-theoretic attractor. (The most attracting Desboves example,
with $b=\frac{1}{3}$, is illustrated in Fig.~\ref{e14}.)
However, in the case of an elementary Desboves map $f$
we get a surprising bonus. The invariant line $\;\{y=0\}\;$ is also carried
into itself by the Latt\`es map of Eq. (\ref{eq:Lat}). Hence it also has a
canonical ergodic invariant measure, and a well defined
transverse Lyapunov exponent.
According to Fig. \ref{fig:elemg}, for real values of $\;b\,$
{\it both\/} of these transverse exponents
are strictly negative provided that $\;|b|\;$ is fairly small and non-zero.

\begin{Quote}
\begin{theo}{\bf Basins of Positive Measure.}\label{th:el}
Let $f$ be a real or complex elementary map with an invariant elliptic
curve $\cC$. If the transverse Lyapunov exponent $\;{\rm Lyap}_{\,\cC}$
is strictly negative,  then the attracting basin $\;B(\cC)\,$,
consisting of points whose orbit converges to ${\cC}\,$,
has strictly positive measure. In fact any neighborhood of a point of $\;\cC\;$
intersects $\;B(\cC)\;$ in a set of positive measure.
Similarly, if such an $f$ has an invariant line $\mathcal L$ with
strictly negative transverse exponent,
then the attracting
basin for this line has positive measure, and intersects any neighborhood of
a point of this line in a set of positive measure.
\end{theo}
\end{Quote}

In the complex case we can give a much more precise picture.
As usual, define the Fatou set to be the largest open set on which
the sequence of iterates of $\;f\;$ forms a normal family, and define the
Julia set $\,J\,$ to be its complement in $\;\bP^2(\C)\,$. If $\;p\;$
is any point of an invariant elliptic curve, then according to
Bonifant and Dabija  in \cite[Theorem 5.4 and Proposition 6.16]{BD}
the iterated preimages of $\;p\;$ are everywhere dense in the Julia set;
and furthermore:

\begin{Quote}
\begin{prop}{\bf The Fatou Set is a Dense Open Basin.}\label{prop:bd}
If $\;f\;$ is a complex elementary map with smooth invariant elliptic curve,
and if the center $p_0$ is not a point of indeterminacy, then
$p_0$ is a superattracting fixed
point whose basin coincides with the Fatou set. This basin
is connected and everywhere dense in $\bP^2$. Furthermore,
if $\;U\;$ is a small neighborhood of a point of the Julia set,
then the union of the forward images of $\;U\;$ is the entire space
$\;\bP^2\ssm\{ p_0\}\,$.
\end{prop}\end{Quote}

In particular, it follows that the attracting basin for the
the elliptic curve has no interior points.
 Similarly, if there is an invariant line $\mathcal L$ disjoint from
the center $p_0$, then the attracting basin of $\mathcal L$
cannot have any interior point. It also follows that
$f$ is topologically transitive on the Julia set. This means that
the orbit of a ``generic'' point of the Julia set $J$ is everywhere dense
in $J$. Such a generic point of $J$ cannot belong to any of the
attracting basins $B(\cC)\,,~B({\mathcal L})$, or $B(p_0)$.

\begin{Quote}
\begin{coro}{\bf Intermingled Basins.}\label{cor:rb}
If a complex elementary map has both
an invariant line $\mathcal L$ which does not pass through its center,
and an invariant elliptic curve $\;\cC$, then the
two topological closures $\overline{B(\cC)}$ and $\overline{B({\mathcal L})}$ are both
precisely equal to the Julia set. Furthermore, if the transverse
Lyapunov exponent for $\cC$ $($or for ${\mathcal L})$ is negative,
then every neighborhood
$\;U\;$ of a point of the Julia set intersects $B(\cC)$ 
$(\!$or respectively $B({\mathcal L})\,)$ 
in a set of positive Lebesgue measure.
\end{coro}\end{Quote}

\begin{rem}{\bf Three Basins.}\label{3bas}
\quad In the case where the transverse exponents  $\;{\rm Lyap}_\cC(f)\;$ and
$\;{\rm Lyap}_{{\mathcal L}}(f)\;$ are both negative, it follows that the
basins for these two attractors are intimately intermingled. 
For the real Desboves map illustrated in Fig. \ref{e14}, a very rough
estimate suggests that about 66\% of the points in $\;\bP^2\;$ are
attracted to the center $\;(0:1:0)\,$, about 17\% to the line
$\;\{y=0\}\,$, and about
17\% to the Fermat curve. For the associated complex mapping, the figures
are 81\%, 13\%, and 6\%. (However, the computation is highly sensitive, and
these estimates may well be quite inaccurate.) It may be conjectured that
every point outside of a set of
measure zero lies in the union of these three attracting basins.
\end{rem}
\smallskip

\begin{rem}{\bf Terminology.}\label{phys-lit} \quad Such exotic behavior has been
studied extensively, particularly in the applied dynamics literature. The term
``riddled basin'' was introduced in \cite{AKYY} to indicate an
attracting basin whose complement intersects every disk in a set of
positive measure. They define two basins to be {\bit intermingled\/}
if every disk which intersects one basin in a set of positive measure
also intersects the other basin in a set of positive measure. For a
particularly clear example, see \cite{Kn}. Such examples of intermingled
basins seem to be known only in cases where the
attractors themselves are quite smooth---We don't know whether there can be
two fractal attractors whose basins have the same closure.
\end{rem}

{\bf Proof} {\bf of Theorem \ref{th:el}.}
Without loss of generality, we may
assume that the center of the real or complex elementary map is $(0:1:0)$.
Furthermore,  an invariant line with negative transverse exponent
certainly cannot pass through this center; so if there is such
a line $\cal L$, we may assume that it is the line $\{y=0\}$,
as in Fig. \ref{e14}.
(Since we assume that there is an invariant elliptic curve $\cC$, it follows
that any invariant line not passing through the center is mapped to itself by
a Latt\`es map, with an absolutely continuous invariant measure so that the
transverse exponent is well defined.)
Each fiber $\;(x:z)={\rm constant}\;$
of the fibration $\;\pi(x:y:z)= (x:z)\;$ has a canonical flat metric
\begin{equation}\label{e15}
 |dy|/\sqrt{|dx|^2+|dz|^2}\,,
\end{equation}
which gives rise to a norm $\;\|\vec v\,\|_{\rm t}\;$ for vectors tangent to
the fiber. Let
$$ \|f'_{\rm t}(p)\|= \|f'\vec{v}\,\|_{\rm t}/\|\vec{v}\,\|_{\rm t} $$
be the  norm of the partial derivative along the fiber; where $\;\vec v\;$
can be any non-zero vector tangent to the fiber at $\;p\,$. (Note that any
vector tangent to its fiber must map to a vector tangent to the image fiber.)
This norm
is well defined, depending only on the base point $p$ of $\vec{v}\,$. At
points of the curve $\;\cC\,$,
we want to compare $\;\|\vec v\,\|_{\rm t}\;$ with the semi-definite norm
$\;\|\vec v\,\|_{\tr}\;\;$ which is obtained by first projecting $\;\vec v\;$
to the quotient vector space $\;T(\bP^2, p)/T(\cC, p)\;$ and then using a
positive definite norm in this quotient space. Note that most fibers
intersects the degree three curve $\cC$ transversally in three distinct points.
There are only a finite number of exceptional fibers which intersect
tangentially. Therefore the
ratio $\;\|\vec v\,\|_{\tr}\;/\|\vec v\,\|_{\rm t}\ge 0\;$ is a continuous
function on $\;\cC\;$ which vanishes only at these points of tangency.
Furthermore, the logarithm $\;\ell(p)\;$
of this ratio has only logarithmic singularities, and hence is an integrable
function on $\;\cC\,$. Since the measure $d\lambda$ is $f$-invariant,
it follows that the difference
$$ \int_\cC\log\|f'_{\tr}\;\|\;d\lambda\,-\int_\cC\log\|f'_{\rm t}\|\;d
\lambda~~=~~\int_\cC\ell\circ f\;d\lambda\,-\,\int_\cC\,\ell\,d\lambda$$
is zero. In other words, the average value
$~      \int_{\cC}\,\log\|f'_{\rm t}\|\, d\lambda ~$
coincides with the transverse exponent ${\rm Lyap}_{\cC}$ of \S\ref{s:texp}.

For any $\;p\;$ and $\;q\;$ belonging to the same fiber, let $\delta(q,p)\ge 0$
be the distance of $q$ from $p\,$, using the flat metric of Eq.
(\ref{e15}) on this fiber. Then
$$  \delta(f(q), f(p))~=~\|f'_{\rm t}(p)\|\,
 \delta(q,p) \,+\,o(\delta(q,p)) $$
as $\;\delta(q,p)\;$ tends to $\,0\,$. This estimate holds uniformly throughout
a neighborhood
of ${\cC}$. Hence, given any $\epsilon>0\,$, we can choose $\delta_0$ so that
\begin{equation}\label{e16}
\delta(f(q),f(p))\le\big(\|f'_{\rm t}(p)\|+\epsilon\big)\,\delta(q,p)
\end{equation}
when $\quad p\in\cC$  and
$\delta(q,p)\,<\,\delta_0$ with $\pi(p)=\pi(q)~ .$

Choose $\epsilon$ small enough so that
\begin{equation}\label{e17}
\int_{\cC} \log\big(\|f'_{\rm t}(p)\|+\epsilon\big)\,d\lambda(p)~<~0 ~.
\end{equation}
Let $p_0\mapsto p_1\mapsto\cdots$ be the orbit of an arbitrary initial
point $p_0\in\cC$ under $f$.
By the Birkhoff Ergodic Theorem, the averages
$$  {\frac{1}{ n}}\Big(\log(\|f'_{\rm t}(p_0)\|+\epsilon)\,+
 \,\log(\|f'_{\rm t}(p_1)\|+\epsilon)\,+\,\cdots
\,+ \,\log(\|f'_{\rm t}(p_{n-1})\|+\epsilon)\Big) $$
converge to the integral  of Eq. (\ref{e17}) for almost all $\;p_0\in\cC\,$.
In particular, for almost all $p_0$ it follows that the number
$$  \log\big(\|f'_{\rm t}(p_0)\|+\epsilon\big)\;+
\;\log\big(\|f'_{\rm t}(p_1)\|+\epsilon\big)\,+\,\cdots\,+
 \,\log\big(\|f'_{\rm t}(p_{n-1})\|+\epsilon\big)~ \hfill{}$$
\vskip -.2in
$$\hfill{} =~\log\Big((\|f'_{\rm t}(p_0)\|+\epsilon\big)\,\cdots
 \,(\|f'_{\rm t}(p_{n-1})\|+\epsilon\big)\Big)$$
is negative for large $n$. Hence the maximum
$$M(p_0)~=~\max_{n\ge 0}\Big((\|f'_{\rm t}(p_0)\|+\epsilon\big)\,
 (\|f'_{\rm t}(p_1)\|+\epsilon\big)\,\cdots\,(\|f'_{\rm t}(p_{n-1})\|
 +\epsilon\big)\Big)~\ge~ 1 $$
is well defined, measurable, and finite almost everywhere. If
$\delta(q,p_0)\le \delta_0/M(p_0)$, then it follows from the inequality
 of   (\ref{e16})
that $\;\delta\big(f^{\circ n}(q)\,,\,f^{\circ n}(p_0)\big)\le\delta_0\;$
for all
$\;n\,$, and also that this distance converges to zero as $\;n\to\infty\,$.
Now let $\;S\;$
be the set of positive measure
consisting of all $\;q\;$ with $\;\delta(q,p)\le\delta_0/
M(p)\;$ for some $\;p\in{\cC}\,$ with $\pi(p)=\pi(q)$.
(Here we take $\delta_0/\infty$ to be $0$.) Then for all $\;q\in S\;$
it follows 
that the orbit of $\;q\;$
converges to $\;{\cC}\,$. Evidently, \;S\; intersects every neighborhood of
a point of $\;\cC\;$ in a set of positive measure. The proof for the basin
of the line $\;y=0\;$ is completely analogous. \qed

\smallskip

\begin{rem}{\bf What about Positive Transverse Exponent?}\label{rm:Lpos}
 Conversely, it seems natural to conjecture that the basin of $\;\cC\;$ has
measure zero whenever the transverse Lyapunov exponent is positive. However,
this cannot be proved simply by reversing the inequalities in the argument
above. The problem is that $\;\log\big(\|f'_{\rm t}\|-\epsilon\big)\;$
is not a meaningful approximation to $\;\log\|f'_{\rm t}\|\,$, since
$\;\|f'_{\rm t}\|\;$ must sometimes be smaller than any given $\;\epsilon\,$.
In fact,  almost every orbit near $\;\cC\;$ must
pass arbitrarily close to the critical locus of $\;f\,$. But
whenever $\;p\;$ is very close to the critical locus,
there is a real possibility that $\;f(p)\;$ will be much closer
to $\;\cC\;$ than would have been predicted from the differential inequality.
We don't know whether this effect could be strong enough to make $\;\cC\;$
a measure theoretic attractor even in some cases
where the transverse Lyapunov exponent is positive. (Compare Remark
\ref{rm:comp}.)
\end{rem}
\smallskip

{\bf Proof} {\bf  of Corollary \ref{cor:rb}.}
\quad It follows immediately from Proposition
\ref{prop:bd} that the basins of $\;\cC\;$ and $\;{\mathcal L}\;$ are contained in the
Julia set. On the other hand, if $\;p\in\cC\cap{\mathcal L}\;$ then the iterated
preimages of $\;p\;$ are contained in both basins, and are dense in $\;J\,$.
Therefore, the closure of either basin is equal to $\;J\,$.

Now if the open set $\;U\;$ intersects the Julia set, then it contains an
iterated preimage of $\;p\,$. Since $\;f\;$ is an open mapping, it follows
that some forward image $\;f^{\circ n}(U)\;$ is an open neighborhood of
$\;p\,$. If $\cC$ (or $\mathcal L$) has negative transverse exponent, then
by Theorem \ref{th:el} the image  $\;f^{\circ n}(U)\;$ intersects
the corresponding basin in
a set of positive measure. Choosing a regular value of $\;f^{\circ n}\;$
which is a point of density for this intersection, and choosing
a point $\;q\in U\;$ which maps to this regular value, it follows easily that
any neighborhood of $\;q\;$ intersects the corresponding basin in a set of positive measure.\qed

\smallskip

\begin{rem}{\bf The Special Case $b=0$.}\label{rem:degen}
The above discussion  of elementary Desboves maps $f=f_{-1,b,1}$
always assumed that the parameter $\;b\;$ is non-zero.
For the case $\;b=0\,$, we have a much simpler situation.
The center $\;(0:1:0)\;$ then becomes a point of indeterminacy.
If we think of each fiber $V$ as a one-dimensional complex vector space,
taking $V\cap{\cal L}$ to be its
zero vector, then each fiber maps linearly to a fiber.
In fact it is not hard to see that $f$ is well defined as a holomorphic
map from the complement
$\;\bP^2\ssm\{(0:1:0)\}\;$ onto itself, and that
this complement is ``foliated'' by $\;f$-invariant copies of
the Fermat curve, which intersect only along the locus $\;\F\cap{\mathcal L}\,$.
In particular, the map
$f=f_{-1,0,1}$ has a first integral. (Compare Lemma \ref{lem:fimaps}.)
\end{rem}

\setcounter{lem}{0}
\section{Trapped Attractors: Existence and Nonexistence.}\label{s:sa}
The first half of this section will provide explicit
examples of everywhere defined rational maps of $\bP^2(\R)$ which have a
smoothly immersed real elliptic curve as trapped attractor. The second half
will prove that a complex elliptic curve, with or without singular points,
can never be a trapped attractor.

\begin{ex}{\bf A Real Elliptic Curve as Trapped Attractor.}\label{ex:cas}
This last example will study the case of a singular elliptic curve. As in
\cite[\S8.6]{BD}, consider the Cassini quartic curve \;${\mathcal C}\;$ with
homogeneous equation \;$\Phi(x,y,z)=0\,$, where\thinspace\footnote
{This expression yields curves which are equivalent, under a complex linear
change of coordinates, to quartic curves introduced in 1680 by the
French-Italian
astronomer Giovanni Domenico Cassini, in connection with planetary orbits.}
\begin{equation}\label{eq:CQ}
    \Phi(x,y,z)~=~\Phi_k(x,y,z)~=~x^2y^2\,-\,(x^2+y^2)z^2\,+\,kz^4
\end{equation}
depends on a single parameter \;$k\ne 0,1\;$. Over the complex
numbers, this is an elliptic curve with nodes at the two points \;$(1:0:0)\;$
and \;$(0:1:0)\,$.
That is, the uniformizing map
\;${\mathbb C}/\Omega\to{\mathcal C}\subset{\mathbb P}^2({\mathbb C})\;$ has
transverse self-intersections at these two points. Define a one-parameter
family of homogeneous polynomial maps from $\;\C^3\;$ to itself
by the formula $\;F(x,y,z)=F_a(x,y,z)=(X,Y,Z)\,$, where
\begin{equation}\label{eq:cfunc}
    X=-2xy(x^2+y^2-2kz^2),\,
    Y=y^4-x^4,\, Z=-a\,\Phi(x,y,z)+2xy(x^2-y^2)~.
\end{equation}
According to \cite{BD}, the curve $\;\cC\;$ is invariant under the induced
rational map $\;f=f_a:\bP^2\to\bP^2\,$. 
It is not hard to check that the singular point $\;(0:1:0)\in\cC\;$
(the ``north pole'' in Figs. \ref{e20} and \ref{e21})
is a saddle fixed point of $\;f_a\,$, with eigenvalues $-2$
and $0$, and that the point $\;(0:0:1)\;$ (near the center of
these figures) is a superattracting fixed point whenever $a\ne 0$.

If the parameters $\;k\;$ and $\;a\;$ are real, then the corresponding real
curve \;${\mathcal C}_{\mathbb R}={\mathcal C}\cap{\mathbb P}^2({\mathbb R})\;$
is connected when \;$k<0\;$, but has two connected
components when \;$k>0$.
These maps are illustrated in Figs. \ref{e20} and \ref{e21} for the case
$\;k=1/8>0$, with $\cC_\R$ in black. Here the smaller nonsingular component
maps to the larger singular component.
The two branches of $\;\cC_\R\;$
through the singular point $(1:0:0)$ map to one of the two branches through
the fixed point
$(0:1:0)$ while both branches through $(0:1:0)$ map to the other
branch through $(0:1:0)$. (All four branches lie within a single
immersed circle which crosses itself twice within the nonorientable manifold
$\bP^2(\R)$.) Note the identities
$$  F(-x,-y,z)~=~F(x,y,z)\qquad{\rm and}\qquad F^{\circ 2}(-y,x,z)~=~
 F^{\circ 2}(x,y,z)\,, $$
which imply that the Julia set of $\;f_a\;$ has $90^\circ$-rotational
symmetry about the point $\;(0:0:1)\,$.\smallskip

Let $\;\cC^0_\R\;$ be the connected component of the real curve
which contains this singular fixed point, and hence maps onto itself.
We will 
prove the following result.
\end{ex}


\begin{Quote}
\begin{theo}{\bf A Trapped Attractor.}\label{th:cas}
If $0<|k|< 1/4$, and if $a$ is sufficiently
small, then the curve ${\mathcal C}^0_{\mathbb R}$ is a trapped attractor
for the map $f_a$ on the real projective plane.
\end{theo}
\end{Quote}
\smallskip

{\bf Proof.}
Let $\;(X,Y,Z)=F(x,y,z)\,$. The quotient
$$ \Phi_F(x,y,z)~=~\Phi(X,Y,Z)/\Phi(x,y,z) $$
is a polynomial of degree 12 in $x,y,z$, depending on the parameter $a$.
In general this polynomial seems rather complicated, but in the special case
$\;a=0\;$ computation shows that it takes the simple form
\begin{equation}\label{e22}
\Phi_F(x,y,z)= 16\,k\,x^2y^2(x^2-y^2)^4~.
\end{equation}
As a convenient measure of the distance of a point of ${\mathbb P}^2$ from the
curve $\Phi=0$ we take the ratio
$$ r(x:y:z)~=~ |\Phi(x,y,z)|/(x^2+y^2)^2\,,$$
with $\Phi$ as in (\ref{eq:CQ}).
This ratio is finite except at the value
$r(0:0:1)=+\infty$, and it vanishes only on the Cassini curve.
We want to prove an inequality of the form
\begin{equation}\label{e23}
r(X:Y:Z)~\le ~\lambda\; r(x:y:z)
\end{equation}
whenever $r(x:y:z)$ is sufficiently small, where $\lambda<1$ is constant.
To do this, we consider the ratio
$$ r_f(x,y,z)~=~{\frac{r(X:Y:Z)}{ r(x:y:z)}}~=~
 \left|{\frac{\Phi(X,Y,Z)}{\Phi(x,y,z)}}\right|\;
 {\frac{(x^2+y^2)^2}{(X^2+Y^2)^2}}~. $$
In the special case $a=0$, using the identity (\ref{e22}) and the inequality
\begin{equation}\label{ineq:X}
 X^2+Y^2~ \ge~ Y^2=(y^4-x^4)^2=(x^2+y^2)^2(x^2-y^2)^2\,,
\end{equation}
together with $|2xy|\le x^2+y^2$, we see that
$$ {\frac{r(X:Y:Z)}{ r(x:y:z)}}~\le~
\left|{\frac{\Phi(X,Y,Z)}{ \Phi(x,y,z)}}\right|\,{\frac{(x^2+y^2)^2}{
(x^2+y^2)^4(x^2-y^2)^4}}= {\frac{16|k|x^2y^2}{(x^2+y^2)^2}}~\le~4|k|~.$$
If $\;0<|k|<1/4\,$, then we can choose $\;\lambda\;$ so that
$\;4|k|<\lambda<1\,$. If $\;N\;$ is any compact
subset of ${\mathbb P}^2({\mathbb R})\ssm\{(0:0:1)\}$, then for any $a$ which
is sufficiently close to zero it then follows by continuity
that the required inequality (\ref{e23}) holds uniformly throughout $N$.
Thus all orbits of $f_a$
in $N$ converge uniformly towards the subset ${\mathcal C}_{\mathbb R}$. In
the case where there are two components, the image
$f_a({\mathcal C}_{\mathbb R})$ is necessarily equal to the component
${\mathcal C}_{\mathbb R}^0\subset {\mathcal C}_{\mathbb R}$
which contains the fixed point, so it follows
that all orbits in $N$ converge uniformly to ${\mathcal C}_{\mathbb R}^0$.
 \qed

\begin{rem}{\bf The Case $a=0$.}\label{rem:0case}
In the limiting case $\;a=0\,$, there is no superattracting
point, and in fact $\;(0:0:1)\;$ becomes a
point of indeterminacy. The proof shows that the basin of $\;\cC^0_\R\,$
under this limiting map $\;f_0\;$ is its entire domain of definition
$\;\bP^2(\R)\ssm\{(0:0:1)\}\,$,
provided that $\;0<|k|<1/4\,$.
\end{rem}

Now let us work over the complex numbers. 
We will show that a
complex elliptic curve can never be a trapped attractor. The proof will
occupy the rest of this section.

\begin{Quote}
\begin{theo}{\bf  No Complex Trapping.}\label{th:notrap}
Let $\,\cC\subset\bP^2=\bP^2(\C)\,$ be an elliptic curve
and let $N$ be a neighborhood of $\;\cC$. Then there cannot exist any
holomorphic map $~f:N\to N~$ such that $\bigcap_n f^{\circ n}(N)=\cC$.
\end{theo}
\end{Quote}

We first carry out the argument for a smooth elliptic curve (necessarily
of degree three), and then show how to modify the proof for a singular
curve (necessarily of degree greater than three).
The proof for a smooth $\;\cC$
will be based on the following construction.

Let $N_\epsilon$
be the $\epsilon$-neighborhood consisting of all points with distance
less than $\epsilon$ from  $\;\cC$, using the standard Fubini-Study
metric\footnote{In terms of homogeneous coordinates
$(x_0:x_1:\cdots:x_n)$ normalized so that $\sum|x_k|^2=1$, this metric
takes the form $dt^2=\sum|dx_k|^2-|\sum\overline x_k\, dx_k)|^2$.
With this normalization,
the Riemannian distance  $0\le\theta\le\pi/2$ between two points
$(x:y:z)$ and $\;(u:v:w)\;$ of $\;\bP^2\;$ can be computed
by the formula
$\;\cos(\theta)=|\overline xu+\overline yv+\overline zw|\,$.}
on $\bP^2$. If $\cC$ is smooth and $\epsilon$ is
sufficiently small, then $N_\epsilon$ is
the total space of a real analytic fibration $\pi:N_\epsilon\to\cC$, where
$\pi(p)$ is defined to be that point $q\in\cC$ which is closest to $p$.
Furthermore, although this projection map $\pi$ is not holomorphic,
each fiber $F_q=\pi^{-1}(q)$ is a holomorphically embedded complex disk
which is contained in the complex line orthogonal to $\,\cC\,$ at $q$.

\begin{Quote}
\begin{lem}{\bf Curves in a Neighborhood.}\label{lem:transv}
With $\,\cC\subset N_\epsilon\,$ 
as above, any non-constant
holomorphic map $\phi:\cC_1\to N_\epsilon$ from an elliptic curve into
$N_\epsilon$ must be an immersion, and must intersect each fiber $F_q$
transversally, so that the composition $\pi\circ\phi$ is a real analytic
immersion of $\;\,\cC_1\,$ onto $\;\cC\,$ of degree $\ge 1$.
\end{lem}
\end{Quote}

{\bf Proof.}
 Suppose to the contrary that there exists a critical point
for the composition $\pi\circ\phi:\cC_1\to\cC$. It will be convenient
to rotate the coordinates for $\bP^2$ as follows. Using $(x,y)$ as an
abbreviation for the point with coordinates $(x:y:1)$, we may assume that
the critical value in $\;\cC\subset \bP^2$ has coordinates $(0,0)$ and that
the tangent line to $\;\cC\,$ at this critical value has equation $y=0$.
The fiber through this point is then a disk in the line $x=0$.
Choose a parametrization $t\mapsto\big(x(t),y(t)\big)$ for $\;\cC\,$ so that
$x(0)=y(0)=y'(0)=0$, and choose a parametrization
$s\mapsto\big( x_1(s),y_1(s)\big)$ for $\phi(\cC_1)$ so that the critical
point in $\cC_1$ is $s=0$, mapping to a point
$\big(x_1(0),y_1(0)\big)=(0, y_0)$ which
lies in the fiber $x=0$. Now expand the function $x_1(s)$ as a power series
$$  x_1(s)~=~c\,s^n\,+\,{\rm(higher~order~terms)}\,,$$
with $c\ne 0$. Here $n\ge 2$, since otherwise $\phi(\cC_1)$ would cross
the fiber $x=0$ transversally.

Using coordinates $(x,y)$ for $N_\epsilon$ and $t$ for $\;\cC$, we can think
of  the real analytic projection $\pi:N_\epsilon \to \cC$ as a correspondence
$\pi:(x,y)\mapsto t=\widehat t(x,y)$.  Setting $y=y_0+\eta$,
we can write the power series expansion for $\widehat t(x,y)$
at the point $(0, y_0)=\big(x_1(0),\,y_1(0)\big)$ as
$$  \widehat t= x\big(a_1+(a_2x+a_3\overline x+a_4\eta+a_5\overline\eta)
+\cdots\big)+\overline x\big(b_1+(b_2x+b_3\overline x+
b_4\eta+b_5\overline\eta)+\cdots\big),$$
where the dots stands for terms of degree $\ge 2\;$ in $\;x,\;\overline x,\;
\eta,$ and $\overline\eta$, and where the $a_j$ and $b_j$ are
complex numbers with $|a_1|>|b_1|$ since the projection
from $(x,\,y_0)$ to its image in $\;\cC$ must preserve orientation for $x$
near $0$. (Here we can assume that $b_2=a_3$.)
Therefore the composition
$\;s\mapsto(x_1(s),\,y_1(s))\mapsto\widehat t\;$ has power series
$$\widehat t~=~a_1\,c\,s^n\;+\;b_1\,\overline c\;\overline s^n\;+\;{\rm(higher
~order~terms)}\,.$$
This proves that the composition $\pi\circ\phi:\cC_1\to\cC$ has an isolated
critical point of local degree $n\ge 2$ (and hence multiplicity $n-1\ge 1$)
at the point $s=0$.
Thus every critical point is isolated, and it follows that
$\pi\circ\phi$ is a branched covering, with only finitely many critical points.

As in Remark \ref{rm:highg}, we apply
the Riemann-Hurwitz Theorem, which asserts that the Euler characteristic
$\chi(\cC_1)$ is equal to $\delta\,\chi(\cC)-\nu$, where $\delta$ is the degree
of $\pi\circ\phi$ and where $\nu$ is the number
of critical points counted with multiplicity. Since $\chi(\cC_1)=\chi(\cC)
=0$, this proves that $\nu=0$, as required. \qed
\bigskip

In particular, it follows that the composition $\pi\circ\phi:\cC'\to\cC$
is quasiconformal.
 We will need the following.\medskip

{\bf Definitions.}
 The  {\bit complex dilatation\/} of
a quasi-conformal map $g(z)$ is the ratio
$$\mu_g(z)~=~(\partial g/\partial\overline z)/(\partial g/\partial z)\,.$$
(The absolute value $|\mu_g|<1$ is sometimes known as the ``small dilatation'',
while the ratio
 $$ K(z)~=~(1+|\mu_g|)/(1-|\mu_g|)~\ge~ 1$$ is called the {\bit dilatation\/}.
The associated {\bit Teichm\"uller distance\/} can be defined
by first taking the maximum $\max K(z)$ as $z$ ranges over the Riemann
surface, and then
taking the infimum of $\;\log\!\big(\!\max K(z)\big)\;$ over an isotopy class
of maps.)

Now consider an infinite sequence of holomorphic immersions
$\phi_j:\cC_j\to N_\epsilon$, where each $\;\cC_j$ is a compact
Riemann surface of genus one.

\begin{Quote}
\begin{lem}{\bf Converging Quasiconformal Maps.}\label{lem:dil20}
If the successive images $\phi_j(\cC_j)$ converge to $\;\cC$ in the Hausdorff
topology $($or in other words if the distance of each point of $\phi_j(\cC_j)$
from $\;\cC$ converges uniformly to zero$)$, then the
complex dilatation of the quasiconformal map $\pi\circ\phi_j:\cC_j\to\cC$
converges uniformly to zero as $j\to\infty$.
\end{lem}
\end{Quote}

The proof of Lemma \ref{lem:dil20} will be based on the following preliminary
statement, which is needed in order to control first derivatives.

\begin{Quote}\begin{lem}{\bf Converging Tangent Spaces.}\label{lem:tanlim}
With the $\phi_j$ as in Lemma \ref{lem:dil20}, consider a
sequence of points $\phi_j(p_j)\in\phi_j(\cC_j)$ converging
to a point $q\in\cC$. Then the tangent space to $\phi_j(\cC_j)$ at
$\phi_j(p_j)$ converges to the tangent space to $\;\cC$ at $q$. Furthermore, this
convergence is uniform as we vary the points  $p_j$ and $q$.
\end{lem}\end{Quote}

It will be convenient to describe a tangent vector at a point $p\in N_\epsilon$
as being {\bit parallel\/} to $\cC$ if it it orthogonal to the fiber through
$p$, as described in the proof of Lemma \ref{lem:transv}. Thus to prove
Lemma \ref{lem:tanlim} we must show that the angle $\theta_j$ between\footnote{
This angle is measured in the tangent space of $\bP^2$ at the
given point. By definition, the angle $\;0\le\theta\le\pi/2\;$
between two complex lines through the origin in a Hermitian vector space
is defined by the formula
$\;\cos(\theta)=|\langle{\bf v}_1,\,{\bf v}_2\rangle|$,
where ${\bf v}_1$ and ${\bf v}_2$ are representative unit vectors. (This
can be identified with the geodesic distance between the corresponding points
in the projective space of lines through the origin.)}
the tangent space at $\phi_j(p_j)$ and the directions parallel to $\cC$ at
this point must converge to zero as $j\to\infty$. The proof will use the
methods of Lemma \ref{lem:transv}, but with a more flexible construction.
\medskip

{\bf Proof} {\bf of Lemma \ref{lem:tanlim}.}  Let $\sigma$
be a $C^\infty$-function which assigns to each point $q\in\cC$ a complex
line $\sigma(q)$ which cuts the curve $\;\cC$ transversally at $q$, and let
$D_\epsilon(q)$ be the disk of radius $\epsilon$ centered at $q$ within this
transverse line. Then the union  $N(\sigma,\,\epsilon)$ of all of these
transverse disks forms a neighborhood of $\cC$. Furthermore,
there exists a number $\epsilon_\sigma>0$ which is small enough
so that the projection $\pi_\sigma:N(\sigma,\,\epsilon_\sigma)\to\cC$
which sends each such disk to its center $q$ is well defined and smooth.
Just as in the proof of Lemma \ref{lem:transv}, any immersion of an
elliptic curve into $N(\sigma,\,\epsilon_\sigma)$
must cut every fiber transversally.

We will construct a compact family $K$ of such transverse line fields $\sigma$.
Since a radius $\epsilon_\sigma$ which is small enough for a given $\sigma$
will also work for any nearby $\sigma$, it follows that, for any such $K$,
we can choose a number
$\epsilon_K>0$ which is small enough so that $\pi_\sigma$ fibers $N(\sigma, \,
\epsilon_K)$ over $\;\cC$ for {\it every\/} $\sigma\in K$.
Furthermore, since the distance of the boundary
of $N(\sigma, \epsilon_K)$ from $\;\cC$ depends continuously on $\sigma$,
it follows that the intersection $N(K,\,\epsilon_K)
=\bigcap_{\,\sigma\in K}N(\sigma,\,\epsilon_K)$ is itself a
neighborhood of $\;\cC$.

Now suppose that Lemma \ref{lem:tanlim}
were false, so that we could find a sequence
of points $p_j\in\cC_j$ such that the corresponding angles $\theta_j$
did not converge to zero. After passing to an infinite subsequence,
we may assume that the images $\phi_j(p_j)$ converge to some point $q_0\in\cC$,
and that the corresponding angles $\theta_j$ are bounded away from zero,
say $\theta_j\ge\alpha>0$.
Choose the compact family $K$ of sections
$\sigma$ such that, for any $q$ in some neighborhood $U_0$ of $q_0$ in $\;\cC$,
every possible transverse line which intersects $\;\cC$ at $q$ at
an angle greater than $\alpha/2$ will occur as the value $\sigma(q)$ for some
$\sigma\in K$. For example we can construct such a $K$ as follows:
First choose a nonzero tangent vector field $\vec v(q)$ and a nonzero
orthogonal vector field $\vec w(q)$, where $q$ varies over some neighborhood
$U$ of $q_0$ within $\cC$. Then choose a
$C^\infty$ function $\rho:\cC\to\R$ which has compact support contained in $U$,
and which takes the value $+1$ throughout a smaller neighborhood $U_0$
of $q_0$. For any complex constant $a$, define a transverse line field
$q\mapsto \sigma_a(q)$ as follows. For $q\in U$ let $\sigma_a(q)$ be the
transverse line spanned by the vector
$$\vec w_a(q)~=~ a\,\rho(q)\vec v(q)+\vec w(q)\,,$$
while for $q\not\in U$ let $\sigma_a(q)$ be the line orthogonal to $\cC$.
Now we can choose the compact set $K$ to be the set of
all $\sigma_a$ with $|a|\le c$, where $c$ is some large constant.

Then for $j$ sufficiently large, the point
$\phi_j(\cC_j)$ will be contained in the intersection $N(K,\,\epsilon_K)$.
Furthermore, the tangent line to $\phi_j(\cC_j)$ at $p_j$ will intersect
$\;\cC$ at a point of the neighborhood $U_0$ and at an angle greater than
$\alpha/2$. Hence, for a corresponding choice of $\sigma$, this tangent line
will contain the fiber of $\pi_\sigma$ through $p_j$.
Since this is impossible, we obtain a contradiction
which completes the proof of Lemma \ref{lem:tanlim}. \qed

\smallskip

\begin{rem}{\bf Convergence of Nearby Curves.}\label{rm:curvconv}
 Lemma \ref{lem:tanlim}  actually implies the stronger
statement that, as $j\to\infty$, each branch of the image curve
$\phi_j(\cC_j)$ converges, with all of its derivatives, to $\;\cC$.
(The number of such
 branches is equal to the degree of the map $\pi\circ\phi_j:\cC_j\to\cC$.)
In fact, given any point $q\in\cC$, we
can rotate coordinates, as in the proof of Lemma \ref{lem:transv}, so that $q$
has coordinates $(x,y)=(0,0)$, and so that the fiber $\pi^{-1}(q)$ is contained
in the line $x=0$. Then $\;\cC$ can be described locally as the graph of
a holomorphic function $y=y(x)$, and it follows from Lemma \ref{lem:tanlim}
that, for large $j$, each branch of the image $\phi_j(\cC_j)$ can be described
as the graph of a holomorphic function $y=y_j(x)$ which converges uniformly
to $y(x)$ throughout some neighborhood as $j\to\infty$. It then follows
from a theorem of
Weierstrass that every iterated derivative $d^n y_j(x)/dx^n$ converges
uniformly to $d^ny(x)/dx^n$ throughout a smaller neighborhood.
\end{rem}

{\bf Proof  of Lemma \ref{lem:dil20}.} With coordinates $(x,y)$ as in the
preceding remark, let $t\mapsto\big(x(t),\,y(t)\big)$ be a local
parametrization of the curve $\;\cC$ with $x(0)=y(0)=y'(0)=0$.
It will be convenient to construct
new local holomorphic coordinates $(u,v)$
by the formula $$(x,y)~=~\big(x(u),\, y(u)+v\big)\,.$$
In these new
coordinates, the curve $\;\cC$ has equation $v=0$.
The projection $\pi$ is then represented by a real
analytic map $(u,v)\mapsto \hat t(u,v)$, where $\hat t(u,0)=u$, so that 
$\partial \hat t/\partial u=1$ and $\partial \hat t/\partial\overline u=0$,
when $v=0$. For $j$ large, we can choose $s=u$ as parameter for the nearby curve
$\phi_j(\cC_j)$, so that the parametrization takes the form $s\mapsto
\big(u_j(s),\,v_j(s)\big)$ with $u_j(s)=s$.
For the composition
$s\mapsto \hat t(s,\,v_j(s))$, it follows that
$$\frac{\partial \hat t} 
{\partial\overline s}~=~\frac{\partial\hat t}{\partial\overline u}\;+\;
 \frac{\partial \hat t}{\partial\overline v}\;\overline{dv_j/ds}\,,$$
where $ {\partial \hat t}/{\partial\overline u}$ 
tends to zero as $v_j\to 0$ by the remarks above, and where $dv_j/ds$
tends to zero as $v_j\to 0$ by Lemma \ref{lem:tanlim}.
It follows easily that the complex dilatation
 $$ \mu_{\pi\circ\phi}~=
(\partial\hat t/\partial\overline s)\, / \, (\partial\hat t/\partial s)$$
tends to zero as $j\to\infty$, as required. \qed

\smallskip

{\bf Proof of Theorem \ref{th:notrap} for an embedded curve.}
The proof will be based on the fact that any smooth
elliptic curve  $\;\cC\subset{\mathbb P}^2({\mathbb C})$ can be approximated
arbitrarily closely by other elliptic
curves which are not conformally isomorphic to it. For example,
after a linear change of coordinates, any such $\;\cC$ is defined by an equation
of the form $x^3+y^3+z^3=3k\, xyz$, and by varying the parameter $k$
we can then find nearby curves which are not conformally isomorphic to $\;\cC$.

\begin{figure}[t]
\centerline{\psfig{figure=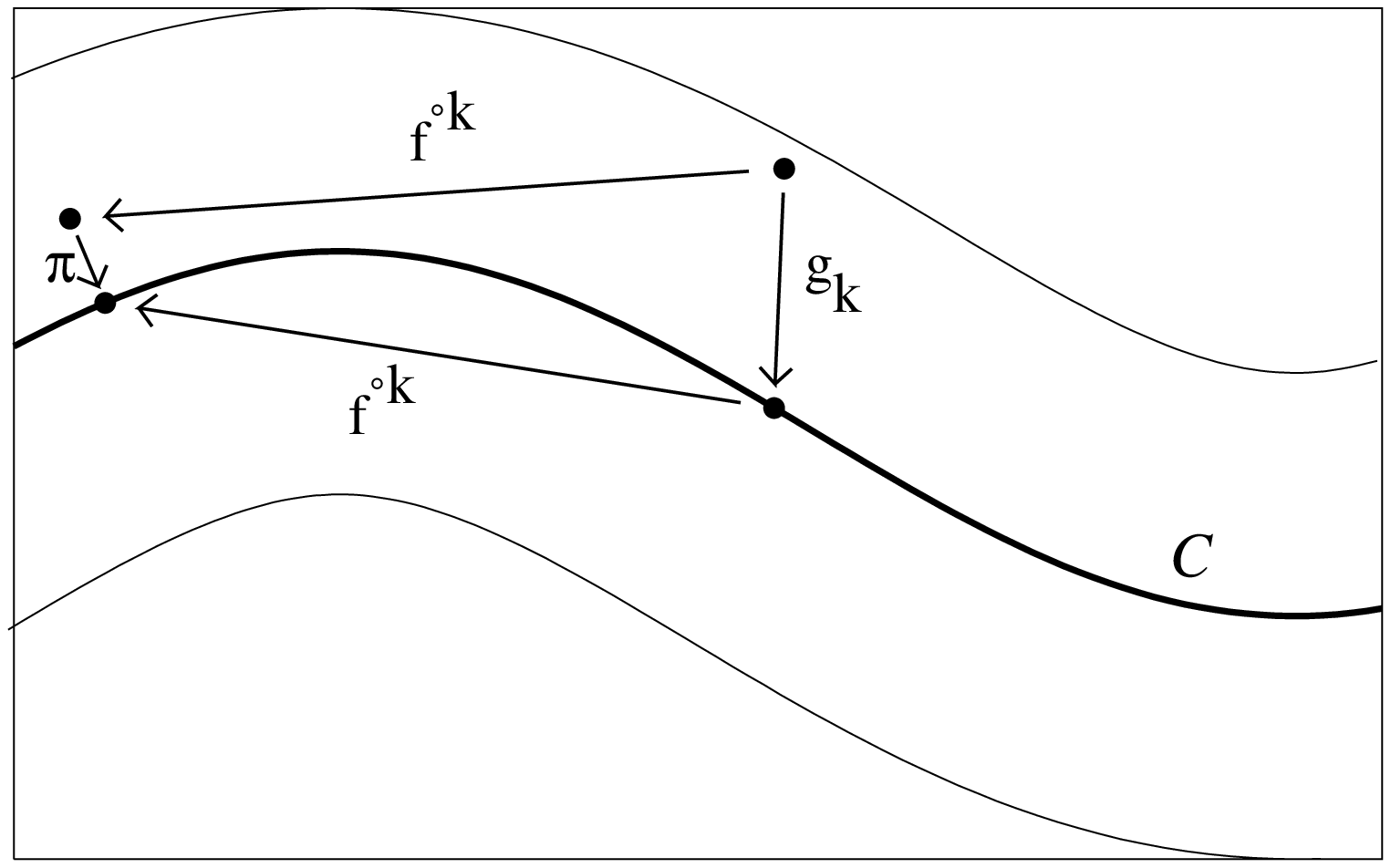,height=2.2in}}
\mycaption{ \label{e26} \bit  Construction of the retraction $g_h$
from $N_\delta$ to $\cC$.\hskip 1.5in}
\end{figure}

Now assume that $\;\cC$ is a trapped attractor. Let $N_\epsilon$ be a tubular
neighborhood, as in Lemma \ref{lem:transv}. Then we can chose
a trapping neighborhood $N_{\rm trap}\subset N_\epsilon$ and then a
smaller tubular neighborhood  $N_\delta\subset N_{\rm trap}$, so that
$$\cC~\subset~N_\delta~\subset~ N_{\rm trap}~\subset N_\epsilon\,.$$
We will then construct a sequence of real
analytic retractions $\,g_h:N_{\delta}\to\cC\,$ of the form
$$  g_h=f^{-h}\circ\pi\circ f^{\circ h}\,. $$
(Compare Fig. \ref{e26}.)
More precisely, since $f^{\circ h}$ is many-to-one, we will construct
 $\,g_h:N\to\cC\,$ so that
$$  f^{\circ h}\circ g_h=\pi\circ f^{\circ h}\,, $$
with $g_h$ equal to the identity map on ${\mathcal C}$.
(Here $\pi:N_\epsilon\to {\mathcal C}$ is again the
orthogonal projection which carries each $p\in N_\epsilon$ to the closest
point of ${\mathcal C}$.)
To do this, let us first pass to the universal covering spaces
$\widetilde\cC\subset\widetilde N_\delta\subset\widetilde N_\epsilon$.
(Since $\cC\cong\C/\Omega$, it follows that $\widetilde\cC\cong\C$.)
Then $f$ and $\pi$ lift to smooth maps
$$  \widetilde {\mathcal C}~\subset~\widetilde N_\delta~\buildrel
 \widetilde f\over
\longrightarrow~\widetilde N_\epsilon~\buildrel \widetilde\pi\over
\longrightarrow~\widetilde {\mathcal C}\,, $$
where we can choose the lift so that $\tilde\pi$ reduces to the identity
map on $\widetilde\cC$.
Since $\widetilde f$ is a linear map of $\widetilde{\mathcal C}\cong\C$,
it follows that $\widetilde f^{-1}$ is well defined.
Therefore the map
$$ \widetilde g_k=\widetilde f^{-k}\circ\widetilde\pi\circ
    \widetilde f^{\circ k}~:~\widetilde N_\delta
~\to~\widetilde {\mathcal C} $$
is well defined; and reduces to the identity
map on $\widetilde {\mathcal C}$.
Finally, since $\widetilde g_k$ commutes with the group of deck
transformations\footnote{On the other hand, the lifted map $\widetilde f$ does
not commute with deck transformations. In fact, if $\Pi$ is the group of
deck transformations, then $f$ induces an
embedding $f_*:\Pi\to\Pi$,
with $\widetilde f(\sigma \widetilde p)=f_*(\sigma)\widetilde f(\widetilde
p)$ for each $\sigma\in\Pi$.}
of $\widetilde N_\delta$ over $N_\delta$, it follows that $\widetilde g_k$
gives rise to a corresponding retraction $g_k:N_\delta\to {\mathcal C}$.

Let $\;\cC'\subset N_{\delta}$ be a smoothly embedded
elliptic curve which is not
conformally equivalent to $\;\cC$. Then it follows using Lemma \ref{lem:transv}
that each $g_h$ maps $\;\cC'$ diffeomorphically onto $\;\cC$. Furthermore, since
the successive images
$f^{\circ h}(\cC')$ must converge towards $\;\cC$, it follows from Lemma
\ref{lem:dil20} that the complex dilatation of the immersion
$\pi\circ f^{\circ h}:\cC'\to \cC$ tends to zero as $h\to\infty$.
Since $f^{\circ h}$ is biholomorphic on both $\;\cC'$ and $\;\cC$,
this implies that the complex dilatation of $g_h:\cC'\to \cC$
also tends to zero as
$h\to\infty$. Thus, to complete the proof of Theorem \ref{th:notrap}, we need
only note the following well known statement from Teichm\"uller theory.

\begin{Quote}
\begin{lem}{\bf A Conformal Isomorphism Criterion.}\label{e30}
 Suppose that there exist quasiconformal homeomorphisms
from the elliptic  curve $\; {\mathcal C}_1$ to ${\mathcal C}_2$ with complex dilatation
arbitrarily close to zero. Then ${\mathcal C}_1$ must be conformally
isomorphic to ${\mathcal C}_2$.
\end{lem}
\end{Quote}

{\bf Proof.} This is an immediate consequence of compactness of the space
of quasiconformal homeomorphisms with bounded complex dilatation. On a more elementary
level, if ${\mathcal C}_1\cong{\mathbb C}/\Omega_1$ and ${\mathcal C}_2\cong{\mathbb C}/\Omega_2$ where
$\Omega_1$ and $\Omega_2$ are unimodular lattices, then the optimal
quasiconformal map in any homotopy class is given by a real-linear map,
corresponding to a linear transformation $L\in{\rm SL}(2,{\mathbb R})$ with
$L(\Omega_1)=\Omega_2$. (Compare \cite[p. 101]{K}.)
Such a linear transformation has complex dilatation zero only if $L$ is a rotation.
Similarly if a sequence of elements of ${\rm SL}(2,{\mathbb R})$
has complex dilatation converging to zero, then some subsequence must converges
to a rotation. The conclusion of Lemma \ref{e30} follows easily.
This completes the proof of Theorem \ref{th:notrap}
for the case of an embedded curve. \qed

\smallskip

{\bf Proof} {\bf  in the Singular Case.} Now consider an elliptic curve $\;\cC\subset
\bP^2$ (necessarily of degree four or more) which has singular points,
so that the uniformizing map $\iota:\C/\Omega\to\cC$ has either critical
points or self-intersections or both.
We will first prove the following preliminary statement.

\smallskip

\begin{Quote}
\begin{lem}{\bf The Branches Fold Together.}\label{branchtogether}
If $\;\cC$ is
a trapped attractor under some rational map $f$ of $\bP^2$, then the
uniformizing map $\iota:\C/\Omega\to\cC$ is necessarily an immersion.
In particular, $\;\cC$ cannot have any cusps. Furthermore, 
some iterate $f^{\circ n}$ must map all of the
branches of $\;\cC$ through any singular point $p$ into a single branch
through $f^{\circ n}(p)$.\end{lem}\end{Quote}

{\bf Proof.} 
Recall that the map $f$ restricted to $\;\cC\,$ lifts to a linear map, which we
will denote by $f^\#$, from $\C/\Omega$ to itself.
First suppose that the uniformizing map $\iota$
has critical points in $\C/\Omega$.
Since there can be only finitely many, and since the lifted map $f^\#$
must send critical points to  critical points,
it follows that there must be a periodic critical point. Thus replacing $f$
and $f^\#$ by some iterate, we may assume that there is a fixed
critical point. In terms of suitable local coordinates around this
point and its image in $\bP^2$, the map $\iota$ will have
power series expansion of the form
$t\mapsto (x,y)=(t^m+\cdots,\;t^n+\cdots)$ with $m> n\ge 2$.
Let $\lambda$ be the multiplier of
$f^\#$. Then the equation $\iota(\lambda\,t)=f\big(\iota(t)\big)$
implies that the eigenvalues of the derivative $f'$ at the critical value
are $\lambda^m$ and $\lambda^n$. Since $|\lambda|>1$, it follows that
this critical value is a repelling fixed point for $f$; which
contradicts the hypothesis that $\;\cC$ has a trapping neighborhood.

Now suppose that we could find two branches through $\iota(t_1)=\iota(t_2)=p$
which map to distinct branches under all iterates of $f$. Since there are
only finitely many singular points, these images must eventually cycle
periodically. Thus, after replacing $f$ by an iterate, we could find two
distinct branches through some singular point $q$ which both map to themselves.
Since $f^\#$ has multiplier $\lambda$, it would follow easily
that the eigenvalues of $f'$ at $q$ have the form $\lambda$ and $\lambda^m$,
where $m\ge 1$ is the intersection number between these two branches.
Again this shows that $q$ is a
repelling point, contradicting the hypothesis that $\;\cC$ has
a trapping neighborhood. \qed
        \smallskip

{\bf The Pulled Back Neighborhood.} For any small $\epsilon>0$ we
can ``pull back'' the $\epsilon$-neighborhood $\;N_\epsilon=N_\epsilon(\cC)$
under the immersion $\iota$
to construct a formal neighborhood of $N_\epsilon^\#\supset\C/\Omega$. For
each $t\in \C/\Omega$, let $D_\epsilon(t,~\iota)\subset\bP^2$ be the open unit
disk in the line normal to $\iota(\C/\Omega)$ at $\iota(t)$, and let
$$ N^\#_\epsilon~=~N^\#_\epsilon(\iota)~\subset~
 (\C/\Omega)\times\bP^2 $$
be the set of all pairs $(t,p)\in\C/\Omega\times\bP^2$ with $t\in\C/\Omega$ and
$p\in D_\epsilon(t,\iota)$. Then $N^\#_\epsilon$ is a
real analytic
manifold, and the projection $\;\pi(t,p)=t\;$ is a real analytic fibration of
$N^\#_\epsilon$ over $\C/\Omega$. Furthermore, if $\epsilon$ is
sufficiently small, then the projection $\;\widetilde\iota(t,p)=p\;$
will be a local diffeomorphism from $N^\#_\epsilon$ onto the open neighborhood
$N_\epsilon\subset\bP^2$. Using this local diffeomorphism
$\widetilde\iota$, we can pull back the complex structure and make
$N^\#_\epsilon$ into a complex manifold.\smallskip

For the next lemma, we assume that $f$ has been replaced by a suitable iterate
$f^{\circ n}$, as in Lemma \ref{branchtogether}.

\begin{Quote}\begin{lem}{\bf Lifting the Trapping
Neighborhood.}\label{stilltrapped}\quad  If the singular curve
$\;\cC$ is a trapped attractor under some rational map $f$
of $\bP^2$ which folds branches together as in Lemma \ref{branchtogether},
then $f$ lifts to a holomorphic map $f^\#$
from a neighborhood of $\C/\Omega$ in $N^\#_\epsilon$ into
 $N^\#_\epsilon$, with  $\C/\Omega$ as trapped attractor.
\end{lem}\end{Quote}

{\bf Proof.} Given a neighborhood $N^\#_\epsilon
=N^\#_\epsilon(\iota)$ as above, by
the uniform continuity of $f$ on the compact set $\overline N_\epsilon(\cC)$,
we can choose $\delta<\epsilon$ so that any curve of length $<\delta$ in
$N_\epsilon(\cC)$ maps to a curve of length $<\epsilon$ in $\bP^2$. We
can then form the neighborhood $\; N^\#_\delta(\iota)
\subset N^\#_\epsilon(\iota)$ of $\;\C/\Omega\,$, with image
$\;\cC\subset N_\delta(\cC) \subset N_\epsilon(\cC)$, and with
$f(N_\delta)\subset N_\epsilon$. We may also assume that $\delta$ is
small enough so that the projection which sends each point of $N_\delta$
to the closest
point of $\;\cC$ is uniquely defined, except within the $\epsilon$-neighborhood
of a branch point.

Now let $T\subset N_\delta$ be a trapping neighborhood for $\;\cC$ and let
$T^\#$ be the full preimage of $T$ in $N^\#_\epsilon$.
Then a lifted map $f^\#:T^\#\to T^\#$ can be
constructed as follows. For each point
$(t, p)\in T^\#\subset N^\#_\epsilon$ we can
drop a perpendicular of length $<\delta$ from $p$ to some point $q\in\cC$.
The image under $f$ will then be a curve of length $<\epsilon$ joining
$f(p)$ to $f(q)\in\cC$. Deforming this curve to a minimal geodesic from
$f(p)$ which meets $\;\cC$ orthogonally, say at $\hat q=\iota(\hat t)$, it follows that
$(\hat t,\,f(p))\in T^\#$, and we will set $f^\#(t,p)=(\hat t, f(p))$.

This construction does not appear to be well defined in the neighborhood of
a singular point, since there may be perpendiculars of length $<\delta$ from
$p$ to points on two or more branches of $\;\cC$. However, by hypothesis
these branches all map to a single
branch of $\;\cC$, so that the minimal geodesic from $f(p)$ to
that branch of $\;\cC$ is unique.

Finally, we must show that the intersection $A^\#$ of the iterated
forward images of $T^\#$ 
is equal to $\C/\Omega$.
Clearly the projection from $T^\#$ onto $T$ maps $A^\#$ onto
$\;\cC$. Therefore $A^\#$ is contained in the preimage of $\;\cC$ in
$T^\#$.
This consists of $\;\cC$, together with preimages of the $\epsilon$-neighborhoods
of the various branch points. (If $\iota(t_1)=\cdots=\iota(t_k)$, then there
are $k-1$ extra preimage branches
 through each of the points $t_1\,,\,\ldots\,,\,t_k$.)
But $f^\#$ maps each of these extra branches back to $\C/\Omega$,
so the attractor $A^\#$ is precisely $\C/\Omega$. \qed
\medskip

 The proof of Theorem \ref{th:notrap}
now proceeds just as in Lemmas \ref{lem:transv} through
\ref{e30} above. However, to carry out the argument in this new context, we
must show that there are nearby curves which are not conformally isomorphic to
$\C/\Omega$. In fact, we will prove the following, which will complete the
proof.

\begin{Quote}\begin{lem}{\bf Deforming Immersed Curves.}\label{curvapprox}
Consider a Riemann surface of genus one of the form $\C/\Omega$.
and let $F_\Omega:\C/\Omega\to\bP^2$ be an immersion.
Then for any small deformation $\Omega_t$ of the lattice $\Omega$, we can construct
a corresponding deformation $F_{\Omega_t}$ of the immersion $F_\Omega$.
\end{lem}\end{Quote}

{\bf Proof.}
 Let  $\widehat\bP^2$  be the projective plane blown up at the three coordinate
 points  $(1:0:0)$, $(0:1:0)$, and $(0:0:1)$. The following two statements are easily verified.

\begin{quote} 1. {\it
Any holomorphic immersion of a Riemann surface into  $\bP^2$  lifts
uniquely  to an immersion into  $\widehat\bP^2$ ,  and any immersion into
 $\widehat\bP^2$    projects to a map into  $\bP^2$ which is an
  immersion, except possibly over the three coordinate points.}
\end{quote}

\noindent(The qualification is necessary since, for example, the
non-immersion\; $\;t\mapsto$  $(t^3:t^2:1)$ from $\bP^1$
into $\bP^2$ lifts to an immersion into $\widehat\bP^2$.)

\begin{quote}2. {\it  We can construct a smooth embedding of $\widehat\bP^2$
as a hypersurface in $\bP^1 \times \bP^1 \times  \bP^1$
  by sending each  $(x:y:z)\in\widehat\bP^2$  to the triple  $(f, g, h)$  where
$$      f = x/y,~~   g = y/z, ~~  h = z/x . $$
  with product   $f g h = 1$.}
\end{quote}

\noindent(The blowup guarantees, for example, that  $x/y$ makes sense, even
at the point $(0:0:1)$.
Interpreted in terms of local coordinates for  $\bP^1$,  the equation
$f g h = 1$  is well behaved,  even when one or two of these functions take
the value $\infty$ . For example, near a point where  $h = \infty$  but  $f$
and  $g$  are finite,  we use   $h^{-1}$  as local coordinate, so that the
 equation takes the form
            $h^{-1} = f g .)$\medskip

Thus to immerse a Riemann surface  $S$  into  $\bP^2$  we need only find three
holomorphic functions  $f, g, h$   from  $S$  to  $\bP^1$   which
yield an immersion of  $S$  into  $\bP^1 \times \bP^1 \times \bP^1,$
and which have product equal to 1, taking care that nothing goes wrong over
the three coordinate points.
(The maps  $f,~ g,~ h$  need not have the same degree. For example the
functions
$       f(t) = g(t) = t ,~~~  h(t)= 1/t^2 $
from  $\bP^1$  to $\bP^1$  yield the smooth quadratic variety   $x z = y^2$.)
\medskip

 As noted above, $F_\Omega$  lifts to an embedding
$~  t~\mapsto~\big(f(t), ~g(t),~h(t)\big)~ $
of $\;\C/\Omega\;$ into the subset
  $\widehat\bP^2\subset\bP^1\times\bP^1\times\bP^1$. Furthermore,
each of the functions $f$ and $g$ can be expressed as a rational function
of the Weierstrass function $\wp_\Omega(t)$ and its derivative
$\wp'_\Omega(t)$.
Choosing some explicit expressions for these rational functions
and setting $h=1/(fg)$, it follows
that the map $t\mapsto\big(f(t),~g(t),~h(t)\big)$ deforms smoothly as
we modify the lattice $\Omega$. Evidently the requirement that this
map project to an immersion into $\bP^2$ will remain satisfied for all
sufficiently small deformations. This completes the proof of Lemma
\ref{curvapprox} and hence of Theorem \ref{th:notrap}. \qed

\setcounter{lem}{0}
\section{Herman Rings in $\bP^2$.}\label{s:hr}

By a {\bit Herman ring\/} for a rational map $f$ of $\bP^2=\bP^2(\C)$ we mean
a (complex  one-dimensional) annulus $H$ which is smoothly embedded in $\bP^2$
and which maps to itself by an irrational rotation under $f$, or under some
iterate of $f$.
In Example \ref{e4} we presented empirical evidence for the existence of
attracting Herman rings, invariant under $f\circ f$, for a substantial
collection of maps in the Desboves family  with real parameters.
This section will explore what we can say more generally
about Herman rings in $\bP^2$.
Note first that it is easy to construct special examples.

\begin{ex}{\bf Rings Contained in a Complex Line.}\label{ex:hr1}  Let
$\,(z_0:z_1)\, \longmapsto\,$  $ \big(p(z_0,z_1)\,:\,q(z_0,z_1)\big) $
be any degree $d$ rational map of ${\bP}^1(\C)$
which possesses a Herman ring. (Compare Remark \ref{rem:1var}, and
see for example \cite{Sh}.) Let $r(z_0,z_1,z_2)$ be any nonzero
homogeneous polynomial of degree $d-1$. Then the map
$$  f(z_0:z_1:z_2)~=~\big(p(z_0,z_1):q(z_0,z_1):z_2\,r(z_0,z_1,z_2)\big) $$
of ${\bP}^2(\C)$
clearly has a Herman ring $H$ which lies in the invariant line $z_2=0$. We can
measure the extent to which this ring is attracting or repelling by
using the ratio $\rho(z_0:z_1:z_2)=|z_2|/\sqrt{|z_0|^2+|z_1|^2}$ as a measure of
distance from the line $z_2=0$. If the ratio
$$  \frac{\rho(f(z_0:z_1:z_2))}{\rho(z_0:z_1:z_2)}~=~|r|\sqrt{\frac{|z_0|^2+|z_1|^2}
{|p|^2+|q|^2}}$$
is less than $1$ everywhere on $H$ then this ring will be locally uniformly
attracting, while if it is greater than $1$ everywhere on $H$ then it will be
locally uniformly repelling.
Thus we can always obtain an attracting $H$ simply by multiplying any given
$r(z_0,z_1,z_2)$ by a constant which is close to zero.
Similarly, if $r(z_0:z_1:0)$ is bounded away from zero
throughout $H$, then we can obtain a repelling $H$ by multiplying $r(z_0,z_1,z_2)$
by a large constant. In the first case, note that $H$ will be contained in
the Fatou set, while in the second case it will be contained in the Julia set.
\end{ex}
\smallskip

Note however that a Herman ring can never be a trapped attractor, since
its boundary points will never be attracted to it. Thus
it can never attract a {\it closed\/} neighborhood of itself.

\begin{ex}{\bf The Ueda Construction.}\label{ex:ahr2}
Here is a quite different procedure.
(Compare \cite[p.~13]{F}.) Recall that the $n$-fold symmetric product
of $\bP^1(\C)$ with itself, that is the quotient
$(\bP^1\times\cdots\times\bP^1)/S_n$ of the $n$-fold product
by the symmetric group $S_n$ of permutations of
the $n$ coordinates, can be naturally identified with $\bP^n(\C)$.
Hence any rational map $f$ of $\bP^1$ gives rise to an everywhere defined
rational map
$(f\times\cdots\times f)/S_n$ of $\bP^n$. In particular, it gives rise to a
map\footnote{More explicitly, if $f(x:y)=\big(p(x,y):q(x:y)\big)$ then the map
$(f\times f)/S_2$ can be described by the formula
$(x_1x_2:x_1y_2+y_1x_2:y_1y_2)~\mapsto~(p_1p_2:p_1q_2+q_1p_2:q_1q_2)~$,
where $p_j=p(x_j,y_j)$ and $q_j=q(x_j,y_j)$.}
$(f\times f)/S_2$ of $\bP^2$. Now if $U_1$ and $U_2$ are disjoint invariant
Fatou components in $\bP^1$, then the product
$U_1\times U_2\subset\bP^1\times\bP^1$ projects diffeomorphically to an
invariant Fatou component in $\bP^2$. In particular, if $U_1$ is a Herman
ring and $U_2$ is the immediate basin of an attracting fixed point, then
the image of $U_1\times U_2$ in $\bP^2$ is
the immediate basin of an attracting Herman ring. Similarly, given
two disjoint Herman rings, we obtain their cartesian product as a Fatou
component,\footnote{This construction can even be used to construct a
Herman ring (or Siegel disk) in the {\it real\/} projective plane: Start with
a rational map of $\bP^1(\C)$ with real coefficients which has
two complex conjugate Herman rings (or Siegel disks). Carry out the
construction described above, and then intersect with $\bP^2(\R)$.}
but given a Herman ring and a parabolic or repelling point
we obtain a Herman ring which is contained in the Julia set.
\end{ex}
\smallskip

{\bf The Transverse Lyapunov Exponent.} For any $f$-invariant Herman ring
or Siegel disk
the transverse Lyapunov exponent is defined much as in the case of an
 $f$-invariant
elliptic curve, but is even more decisive as a test for attraction or
repulsion. However, in these rotation domain cases, this transverse
exponent is no longer a constant, but is rather a convex
piecewise linear real valued function (constant on each invariant circle).
To fix our ideas, we will concentrate on the Herman ring case.

To begin the discussion, note
that for any Herman ring $H\subset\bP^2(\C)$
(and more generally any annulus) there is an essentially unique conformal
embedding $t:H\to\C/\Z$ which maps $H$ isomorphically onto a cylinder of the
form $h_0<\Im(t)<h_1$ in $\C/\Z$. (More precisely, this embedding is unique
up to a translation or change of sign $t\mapsto\pm t+{\rm constant}$.
The difference $h_1-h_0$ is called the {\bit modulus\/}  of $H$.)
We will call $t$ a {\bit canonical parameter\/} on $H$.

\begin{Quote}\begin{lem}{\bf Holomorphic Tubular Neighborhoods.}\label{lem:tubu}
Let  $\;\Gamma_h $
$\,\subset\, H\, \subset\, \bP^2$ be the invariant circle $\Im(t)=h$
contained in a Herman ring in the complex projective plane.
Then  we can parametrize some neighborhood $N=N(\Gamma_h)$ in $\bP^2$
by holomorphic coordinates $(t,u)$, where $u$ ranges over a neighborhood of
zero in $\,\C$ and $t$ ranges over a neighborhood of the circle $\Im(t)=h$ in
$\C/\Z$. Furthermore these coordinates can be chosen so that we have $u=0$
on the intersection $N\cap H$, with $t$ the canonical parameter for $H$.
\end{lem}\end{Quote}

{\bf Proof.} In the dual projective space consisting of all lines in $\bP^2$,
those lines which intersect $\Gamma_h$ form a real 3-dimensional subset.
Hence we can choose a line which misses $\Gamma_h$. After rotating the
coordinates $(x:y:z)$, we may assume that this is the line $z=0$.
In other words, setting $X=x/z\,,~Y=y/z$, we can
introduce the affine coordinates $(X:Y:1)=(x:y:z)$
throughout some neighborhood of $\Gamma_h$. Let $X=X(t)\,,~Y=Y(t)$ be the
canonical parametrization of $H$ throughout this neighborhood.
Then the space of all unit vectors in $\C^2$ which are multiples of
the tangent vector
$\big(\dot X(t)\,,\,\dot Y(t)\big)$ for some
$\big(X(t)\,,\,Y(t)\big)\in\Gamma_h$ has
real dimension 2.  Hence we can choose a fixed unit vector $\vec V\in\C^2$
which is not tangent to $H$ at any point of $\Gamma_h$. The required
coordinates $(t,u)$ are now defined by the formula
$$(t\,,\,u)~\mapsto~\big(X(t)\,,\,Y(t)\big)\;+\;u\,\vec V $$
for all $(t,u)\in(\C/\Z)\times\C$ with both $|\Im(t)-h|$ and $|u|$ sufficiently
small. \qed
        \smallskip

The map $f$, expressed in terms of these tubular coordinates
in a neighborhood of $\Gamma_h$, has the form $(t,u)\mapsto(T,\,U)$,
where $(t,0)$ maps to ($t+\alpha,\,0)$ for some irrational rotation number
$\alpha$. Evidently we can identify the transverse derivative along $H_0$
with the partial derivative
$$ 
 \frac{\partial U}{\partial u}\big(t,0)\,.$$
The transverse Lyapunov exponent is then defined as the horizontal average
$$ {\rm Lyap}(t)~=~ \int_0^1 \log\left|\frac{\partial U}{\partial u}(t+\tau)
\right|\,d\tau\,,$$
where we integrate over the interval $0\le\tau\le 1$. If there are no zeros of
$\partial U/\partial u$ in the strip  $h_0<\Im(t)<h_1\,,~u=0$,
then ${\rm Lyap}(t)$ is an average of harmonic
functions, and hence is harmonic. Since this harmonic function is constant
on horizontal lines, it must be a {\it linear\/} function of the imaginary
part $\Im(t)$. (We will sharpen this statement in Lemma \ref{lem:Jensen}.)
The dynamical implications of this Lyapunov exponent
can be described as follows.

\begin{Quote}\begin{lem}{\bf Attraction or Repulsion.}\label{lm:Weyl}
If ${\;\rm Lyap}(t)$ is negative along the invariant
circle $\,\Gamma_h\subset H$, then a neighborhood of
$\,\Gamma_h$ in the Herman ring
$H$ is uniformly attracting, and hence is contained in the Fatou set of $f$.
On the other hand, if ${\;\rm Lyap}(t)$ is positive, then $\,\Gamma_h$ is contained
in the Julia set.
\end{lem}\end{Quote}

However, in the intermediate case where ${\rm Lyap}(t)$ is identically zero
near $\Gamma_h$, we do not have enough information to decide.
In fact, using Ueda's construction (Example \ref{ex:ahr2}), we can find
examples illustrating both possibilities. We can choose $f$ so that a
neighborhood of $H$, with its
dynamics, is isomorphic to ~(Herman ring)$\times$(Siegel
disk)~ and hence belong to the Fatou set. On the other hand, we can choose $f$
so that $H$ corresponds to ~(Herman ring)$\times$(parabolic point),
and hence belong to the Julia set.\smallskip

{\bf Proof  of Lemma \ref{lm:Weyl}.} Recall that the map $f$ restricted to the
ring $H$ has the form $t\mapsto t+\alpha$ where the rotation number $\alpha$
is irrational. To simplify the notation, let us
translate the canonical parameter $t$ so that it takes real values (modulo one)
on our invariant circle, with $h=0$. Following
Hermann Weyl, for any continuous function $g:\R/\Z\to\C$, the
successive averages
$$A_ng(t)\,=\,(1/n)\sum_{0\le j<n} g(t+j\alpha) $$
converge uniformly to the integral $\,\int_{\R/\Z} g(t)\,dt\,. $
To prove this, note that any continuous $g$ can be uniformly approximated
by trigonometric polynomials 
of the form $\;\sum_{|k|\le N}\, a_k\,e^{2\pi i kt}$. But this
statement is easily verified for each such trigonometric polynomial.

First suppose that the transverse derivative
$\frac{\partial U}{\partial u}(t,0)$ has no zeros
on $\R/\Z$, so that the function
 $t\mapsto g(t)=\log\big|\frac{\partial U}{\partial u}(t,0)\big|$
is finite valued near $\R/\Z$. If the average ${\rm Lyap}(0)$ of $g$ on
$\R/\Z$ is strictly positive (or strictly
negative), then we can choose an integer $n>0$
so that $A_ng(t)$ is strictly positive (or negative) and bounded away from
zero on $\R/\Z$, and hence
throughout some neighborhood of $\R/\Z$ in $\C/\Z$. Now consider an
orbit $(t_0,u_0)\mapsto (t_1,u_1)\mapsto\cdots\mapsto(t_n,u_n)$ near the given
circle $\Gamma_0$. Note that
$$ \lim_{u_0\to 0}\,
 \frac{u_n}{u_0}~=~\frac{\partial u_n}{\partial u_0}\big(t_0,0)~=~
 \prod_{0\le j<n} \frac{\partial u_{j+1}}{\partial u_j}\big(t_0+j\alpha,\,0)
~=~\exp\big(n\,A_ng(t)\big) \,.$$
Taking the log absolute value of both sides, if $A_ng(t)<\log(c)<0$ on
$\R/\Z$, then $|\partial u_n/\partial u_0|<c^n$ when $u_0=0$, and
it follows easily that $ |u_n|\le c^n\,|u_0|$ uniformly
throughout a neighborhood of our invariant circle. This proves
that some neighborhood $H_0$
in $H$ is uniformly attracted to $\Gamma_0$.
 Similarly, if $A_ng(t)>\log(c)>0$, then a
neighborhood is uniformly repelled, so that $\Gamma_0$
is contained in the Julia set.

Now suppose that the holomorphic function
 $t\mapsto\frac{\partial U}{\partial u}(t,0)$ has zeros along the real axis.
The superattracting case where this derivative is identically zero on $H$
is easily dealt with, so we will assume that these zeros are isolated.
If ${\rm Lyap}(0)<0$, then we can replace $g(t)$ by the truncated function
$g_\nu(t)=\max(g(t)\,,\,\nu)$, where $\nu$ is some negative real constant.
If $\nu$ is sufficiently negative, then the
integral $\int_0^1 g_\nu(t)\,dt$ will still be negative. Hence we can choose
$n$ so that $A_ng(t)\le A_ng_\nu(t)$ \break $<\log(c)<0$, and it again follows that
a neighborhood $H_0$ of $\Gamma_0$  is uniformly attracting. In the case
${\rm Lyap}(0)>0$, we cannot assert that a neighborhood is uniformly
repelling when it contains zeros of $\partial U/\partial u$.
However, it is not hard to check that the function
$t\mapsto{\rm Lyap}(\Im(t))$ is continuous, and hence is positive throughout
a neighborhood. Since the Julia set of $f$ is closed, we can at least conclude
that $\Gamma_0$ is contained in the Julia set. \qed

\begin{Quote}\begin{lem}{\bf Piecewise Linearity.}\label{lem:Jensen}
Let $H\subset\bP^2$ be a Herman ring with canonical parameter $t\in\C/\Z$,
where $\;h_0<\Im(t)<h_1$. Then the function ${\rm Lyap}:(h_0\,,\,h_1)\to\R$
is convex and piecewise linear, with a jump in derivative at $h$ if and only
if the transverse derivative $\partial U/\partial u$ has a zero on the
circle $\Im(t)=h$. In fact the change in derivative at $h$ is equal to $2\pi$
times the number of zeros of $\partial U/\partial u$
on this circle, counted with multiplicity.
\end{lem}\end{Quote}

{\bf Proof.} We will adapt a classical argument due to Jensen. (See
for example \cite[Appendix A]{M3}.)  It will be convenient to use the
abbreviated notation $\phi(t)
=\partial U/\partial u$ for the transverse derivative evaluated at $(t,0)$,
 and $\phi'(t)$ for its derivative. If $\phi$ has
no zeros on the circle $\Im(t)=h$, then we can compute the
derivative of the transverse Lyapunov exponent by differentiating under the
integral sign. Setting $t=\tau+ih$, we have

$$ \frac{\partial }{\partial h}\log|\phi(t)|
~=~\frac{\partial}{\partial h}\Re\big(\log\phi(t)\big)~=~
\Re\Big(\frac{d\log\phi(t)}{dt}\,\frac{\partial t}{\partial h}
\Big)~=~\Re\Big(\frac{\phi'}{\phi}i\Big) $$
and therefore
$$ {\rm Lyap}'(h)
~=~\frac{d}{dh}\int_0^1\log|\phi(\tau+ih)|\,d\tau
~=~\Re\int_0^1 i\frac{\phi'}{\phi}d\tau\,, $$
where $\phi'$ and $\phi$ are evaluated at $t=\tau+ih$ for $0\le\tau\le 1$
with $h$ constant. Briefly, we can write
$$ {\rm Lyap}'(h)~=~ \Re\int_{[0,1]\times\{h\}}i\,\frac{ d\phi}{\phi}\,.$$
Given two numbers $h_0<h_1$ such that $\phi$ has no zeros at
height $h_0$ or $h_1$, we can now compute the difference
${\rm Lyap}'(h_1)-{\rm Lyap}'(h_0)$ as follows. Translating the
parameter $t$ horizontally if necessary we may assume also that $\phi$ has
no zeros on the vertical line
$\Re(t)=0$. Let $R$ be the rectangle consisting of all $t\in\C$ with
$0\le\Re(t)\le 1$ and $h_0\le\Im(t)\le h_1$. Integrating in the positive
direction around the boundary of $R$,
since the integrals around the left and right sides cancel out, we obtain
$$ {\rm Lyap}'(h_0)-{\rm Lyap}'(h_1) ~=~\Re\oint_{\partial R}i
\,\frac{d\phi}{\phi}\,.$$
But the integral $\oint_{\partial R} d\phi/\phi$ is equal to $2\pi i\, N(R)$,
where $N(R)$ is the number of zeros of $\phi$ in $R$, so this equation
reduces to $\;{\rm Lyap}'(h_1)-{\rm Lyap}'(h_0)=2\pi\, N(R)$. This
proves Lemma \ref{lem:Jensen}. \qed

\begin{rem}{\bf Siegel Disks}\label{rm:sd}
There is an analogous statement for the transverse exponent of a Siegel disk
in $\bP^2$. Let us parametrize the Siegel disk by the unit disk
$\mathbb D$ consisting of all $z\in\C$ with $|z|=r<1$, taking $z=0$ for
the center point. Then the transverse Lyapunov exponent can be expressed as a
a convex piecewise linear function of $\log(r)$. In fact, it is determined
up to an additive constant by the identity
$$  \frac{d\;{\rm Lyap}}{d\log(r)}~=~N({\mathbb D}_r)~\ge~0\,,$$
where $N({\mathbb D}_r)$ is the number of zeros of the transverse derivative,
counted with multiplicity, in the disk of radius $r$. The proof, similar
to that above, is essentially Jensen's original computation.
\end{rem}
\smallskip

{\bf Herman Rings for Maps with Real Coefficients.}
Most known examples of Herman rings have been specially constructed.
The surprise in \ref{e4} was to find an apparent example
which appeared out of the blue, with no obvious reason to expect it.
The set of complex rational maps of specified degree with a Herman ring
presumably has measure zero, so that a randomly chosen example will never
have a Herman ring. However, if we consider rational maps with real
coefficients then the situation is different, and the discussion in
Example \ref{e4}
suggests that the set of real parameters which give rise to a complex Herman
ring should have positive Lebesgue measure.

Let $f$ be a rational map of $\bP^2$ with real coefficients, and
suppose that there exists an embedded $f$-invariant circle $\Gamma
\subset\bP^2(\R)$ with irrational
rotation number. If $\Gamma$ is smooth of class $C^2$, then according to
Denjoy's Theorem the restriction $f|_\Gamma$ is topologically conjugate to
a circle rotation. In particular, there is a canonical $f$-invariant
probability measure $d\mu$ with support equal to the entire circle. The
transverse Lyapunov exponent of $\Gamma$ in $\bP^2(\R)$ is then well defined.
Just as in the proof of Lemma \ref{lm:Weyl},  a positive
Lyapunov exponent implies that $\Gamma$ is uniformly repelling, and similarly,
a negative exponent implies that $\Gamma$ is uniformly attracting and hence
is a trapped attractor.
\smallskip

In the real analytic case we can complexify the circle $\Gamma$.

\begin{Quote}\begin{lem}{\bf From Circle to Ring.}\label{lem:circ2hr}
Let $f$ be a rational map of $\;{\bP}^2(\R)$. If $\;\Gamma$ is a real analytic
$f$-invariant circle with Diophantine rotation number, then
the associated map from $\;{\bP}^2(\C)$ to itself possesses a Herman ring
$H\supset \Gamma$. Furthermore, the transverse Lyapunov exponent of $\;\Gamma$
in $\;\bP^2(\R)$ is identical with the transverse Lyapunov exponent of $\;H$
along $\Gamma$.
\end{lem}\end{Quote}

If we exclude the special case where the transverse exponent is exactly zero,
then it follows that $\Gamma$ is repelling (or attracting) in the real
projective plane if and only if a neighborhood of $\Gamma$ in $H$ is repelling
(or attracting) in the complex projective plane.
\smallskip

\begin{rem}{\bf The One-Dimensional Case.}\label{rem:1var}
It is interesting to compare the situation in one variable.
For any odd number $d\ge 3$, the set of degree $d$ rational maps which carry
the real projective line $\bP^1(\R)$ diffeomorphically onto itself is open,
and all possible rotation numbers are realized. If such a map has
Diophantine rotation number, then a similar argument shows that
the corresponding rational map of $\bP^1(\C)$
contains a Herman ring. 
\end{rem}
\smallskip

{\bf Proof of Lemma \ref{lem:circ2hr}.} By a theorem of Herman,
as sharpened by Yoccoz, (see \cite{Y}), any orientation preserving real
analytic diffeomorphism of a circle with Diophantine rotation number $\alpha$
is real analytically conjugate to the rigid rotation
$t\mapsto t+\alpha~({\rm mod}\;\Z)$
of the standard circle $\R/\Z$. That is, there is a real analytic
diffeomorphism $h:\R/\Z\to \Gamma\subset\bP^2(\R)$ so that
$f(h(t))=h(t+\alpha)$. Since $h$ is real analytic,
it extends to a complex analytic
diffeomorphism from a neighborhood of $\R/\Z$ in the cylinder $\C/\Z$ into the
complex projective plane. The image of this extended map on some
neighborhood $\{t\;{\rm mod}\,\Z\, \,; |\Im(t)|<\epsilon\}$
is the required
Herman ring $H\subset\bP^2(\C)$. Evidently the translation $t\mapsto t+\alpha$
on this neighborhood is conjugate to the rational map $f$ on $H$. Further
details are straightforward, since any norm on the normal bundle of $H$ in
$\bP^2(\C)$ will restrict to a norm on the normal bundle of $\Gamma$ in
$\bP^2(\R)$. \qed

\smallskip

Now consider a $C^\infty$-smoothly embedded circle $\Gamma_0$ in a
real 2-dimensional manifold $M$.
 Let $f_0:M\to M$ be a $C^\infty$-smooth map
which restricts to an irrational rotation on $\Gamma_0$,
and which has negative transverse exponent on $\Gamma_0$.
As noted above, this implies that $\Gamma_0$ is a trapped attractor. Let
$N\subset M$  be a trapping neighborhood,
with $f_0(N)\subset{\rm interior}(N)$.
\smallskip

\begin{Quote}
\begin{theo}{\bf Persistence of Invariant Circles.}\label{th:nearhr}
 In this situation,
for any $C^\infty\!$-map $f_\tau$ which is close enough to $f_0$ in the
$C^1$-topology, the intersection
$$\Gamma_\tau~=~\bigcap_nf_\tau^{\,\circ n}(N)$$
of the iterated forward images of $N$ under $f_\tau$
will be a topological circle; 
and $f_\tau$ will map this circle
homeomorphically onto itself with a rotation
number $\rho_\tau$ which varies continuously with $\tau$.
Furthermore, for any finite $k$, if $f_\tau$ is $C^1\!$-sufficiently
close to $f_0$, then $\Gamma_\tau$ will be $C^k\!$-smooth.
\end{theo}\end{Quote}

Here we can expect the continuous function $\tau\mapsto\rho_\tau$ to
have an interval of constancy whenever $\rho_\tau$ takes a rational value.
In fact, a generic map $f_\tau$ with $\rho_\tau=p/q$ will have
an attracting period $q$ orbit contained in $\Gamma_\tau$. In this case,
$\Gamma_\tau$ cannot contain any dense orbit. Evidently such an attracting
orbit will be stable under perturbation.
\medskip



{\bf Proof of Theorem \ref{th:nearhr}.} If a neighborhood of $\Gamma_0$
is orientable, then we can choose local coordinates $t\in\R/\Z$ and
$-\epsilon< x< \epsilon$ throughout some neighborhood of
$\Gamma_0$ so that $\Gamma_0$ is given by the equation $x=0$, and is
mapped to itself by $t\mapsto t+\alpha$ with $\alpha$ irrational.
In the non-orientable case, we can first pass to the 2-fold orientable covering
of a neighborhood and then choose such coordinates.

We will denote the map $f_0$ in these coordinates by $f_0(t,x)=(T,\,U)$.
It will be convenient to use abbreviated notations such as
$$ \partial_tU~=~\frac{\partial}{\partial t}U(t,x)\,,\qquad
 \partial_xU~=~\frac{\partial}{\partial x}U(t,x)\,.$$
Given $0<c<1$, after replacing $f_0$ by some high iterate, we may assume
that $|\partial_x U|< c$ when $x=0$. Furthermore,
after a carefully chosen change of coordinates, we will show
that $\partial_x T=0$ when $x=0$, so that
\begin{equation}\begin{matrix}
\partial_t T~=~1, & \qquad \partial_x T~=~0,\cr
\partial_t U~=~0, &\qquad |\partial_x U|~<~c
\end{matrix}\end{equation}
along the circle $x=0$. In order to construct this change of coordinates,
we first introduce the smooth functions
$$ r_n~=~r_n(t)~=~\partial_x U(t+n\alpha,\,0)\,,\qquad
  s_n~=~s_n(t)~=~\partial_x T(t+n\alpha,\,0)\,,$$
with $|r_n|< c$. Now introduce the change of coordinates
$$ (t,x)~\leftrightarrow~(\,\widehat t,\,x)\qquad{\rm where}\qquad
\widehat t~=~t\,+\,(s_0+r_0s_1+r_0r_1s_2+\cdots)\,x\,.$$
A brief calculation shows that the induced map $(\,\widehat t,x)
 \to(\widehat T,\,U)$ in these new coordinates satisfies the required condition
$\partial_x\widehat T=0$ when $x=0$. Henceforth, we will omit the hat
and simply assume that $\partial_x T=0$ along the circle $x=0$.

Given any $0<\eta<1$, we can first choose a trapping neighborhood
$\;N=\{(t,x)~;~|x|\le b_0\}\;$
for $\Gamma_0$ which is small enough so that the inequalities
\begin{equation}\begin{matrix}\label{ineq:partials}
~|\partial_t T-1|~<\,\eta, & \qquad |\partial_x T|~<~\eta,\cr
|\partial_t U|~<~\eta,~ &\qquad |\partial_x U|~<~c,
\end{matrix}\end{equation}
are valid throughout this neighborhood. We can then
choose $f_\tau$ close enough to $f_0$, so that
$f_\tau(N)\subset{\rm interior}(N)$, so that these inequalities
(\ref{ineq:partials}) remain true for the map $f_\tau(t,x)=(T,U)$.

Now consider a curve $t\mapsto x(t)$ with slope $v(t)=dx/dt$. Setting
$f_\tau\big(t, x(t)\big)$ $=(T,U)$, we have
$$ \frac{d\,T\big(t, x(t)\big)}{dt}~=~\frac{\partial T}{\partial t}+\frac{dx}{dt}
\,\frac{\partial T}{\partial x}\qquad{\rm or~briefly}\qquad
D_tT=\big(\partial_t+v\,\partial_x)T\,,$$
and similarly $\;D_tX=\big(\partial_t+v\,\partial_x)X$.
Given some upper bound $b_1$ for $|v|=|dx/dt|$, we can estimate that
$$ D_tT~>~ 1-\eta-b_1\eta\qquad{\rm and}\qquad |D_tX|~<~\eta+b_1c\,,$$
and hence
$$\left|\frac{dX}{dT}\right|~=~\frac{|D_tX|}{D_tT}~<
~\frac{\eta+b_1c}{1-\eta-b_1\eta}\,.$$
If $\eta$ is sufficiently small, then this upper bound will be strictly
less than $b_1$. More precisely, if
$$ 0~<~\eta~<~\frac{b_1(1-c)}{1+b_1+b_1^2}, $$
then a brief computation shows that $\,D_tT>0\,$ and that
$\,(\eta+b_1c)/(1-\eta-b_1\eta)<b_1$. It will then follow that
the image curve is the graph of a well defined function $X=X(T)$, and
furthermore that the slope of this image curve is bounded by the same constant,
$$|dX/dT|~<~b_1\,.$$
It follows inductively that
each iterated forward image of the initial curve will again
be a well defined curve with $|{\rm slope}|<b_1$.

As trapping neighborhood $N$ we can choose a union of horizontal circles
 $x=\hat x$ with $-b_0\le \hat x\le b_0$. Then each iterated image
$f_\tau^{\,\circ n}(N)\subset N$ will be a corresponding union of a continuum of
curves of the form $x_n=x_n(t)$ with $|dx_n/dt|\le b_1$. If $x_N^-(t)$ is the
infimum of these curves and $x_n^+(t)$ is the supremum, then it follows
easily that the $n$-th forward image of $N$ is given by
$$f_\tau^{\,\circ n}(N)~=~\{(x,t)~~;~~ x_n^-(t)\le x\le x_n^+(t)\}.$$
Here the upper and lower boundaries both satisfy a Lipschitz condition
$$|x_n^\pm(t_1)-x_n^\pm(t_0)|~\le~ b_1|t_1-t_0|.$$
On the other hand, it is easy to check that the
Jacobian determinant is bounded by
$$  |{\rm Jacobian}|~<~ (1+\eta)c+\eta^2\,,$$
and we may assume that this upper bound is strictly less than
one. Hence the areas of these successive images shrink to zero.
Thus $\;\int_0^1 \big(x_n^+(t)-x_n^-(t)\big)\,dt\;$ tends to zero as
$n\to\infty$.
It follows easily that the upper and lower bounding curves tend to a common
Lipschitz limit. {\it This proves that the attracting set
$$\Gamma_\tau~=~\bigcap_{n\ge 0} f_\tau^{\,\circ n}(N)$$
is itself the graph of a function
$\;t\mapsto\lim_{n\to\infty}x_n^{\pm}(t)\;$
which is continuous $($and in fact Lipschitz with Lipschitz constant $b_1)$.}

In order to prove that this attracting curve is $C^1$ smooth, we must estimate
second derivatives. Computations show that
$$ \frac{d^2X}{dT^2}~=~\frac{(D_t^{\,2}X)(D_tT)-(D_t^{\,2}T)(D_tX)}{(D_tT)^3}
\,,$$
where
$$ D_t^{\,2}~=~\partial_t^{\,2}\,+\,2v\partial_t\partial_x\,
+\,v^2\partial_x^{\,2}\,+\,(dv/dt)\partial_x\,.$$
Separating out the $dv/dt=d^2x/dt^2$ terms, we can write this as
\begin{equation}\label{eq:2est}
 \frac{d^2X}{dT^2}~=~A_2\,+\,B_2\,\frac{d^2x}{dt^2}\,,
\end{equation}
where $A_2$ is uniformly bounded and where
$$B_2~=~\frac{\partial_xX}{(D_tT)^2}\;-\;\frac{(\partial_xT)(D_tX)}{(D_tT)^3}
\,,$$
so that
$$|B_2|~<~\frac{c}{(1-\eta-b_1\eta)^2}\;+\;\frac{\eta(\eta+b_1c)}
 {(1-\eta-b_1\eta)^3}\,.$$
Evidently, if $\eta$ is small enough, then $|B_2|<{\rm constant}<1$. We can then
choose a constant
 $$b_2\,>\,|A_2|/(1-|B_2|)\,,
\qquad{\rm so~ that}\qquad b_2\,>\, A_2+B_2\,b_2\,.$$
Then if $|d^2x/dt^2|<b_2$, it follows that $|d^2X/dT^2|<b_2$.
{\it Thus, with these choices, the successive forward images of $\Gamma_0$
will be curves $x_n=x_n(t)$ which converge uniformly to a limit, with both
$|dx_n/dt|$ and $|d^2x_n/dt^2|$ uniformly bounded.}

Now we can continue inductively. By successively differentiating the formula
(\ref{eq:2est}) we find formulas of the form
$$ \frac{d^{k}X}{dT^k}~=~A_k\,+\,B_k\,\frac{d^kx}{dt^k}\,, $$
where $A_k$ depends not only on the iterated partial derivatives of $T(t,x)$
and $X(t,x)$ but also on the derivatives $d^\ell x/dt^\ell$ with $\ell<k$,
and where
$$ B_{k+1}~=~B_k/D_tT~<~B_k/(1-\eta-b_1\eta)\,.$$
Thus, choosing $\eta$ small enough
so that $B_k<{\rm constant}<1$, we can find a suitable upper bound
$b_k$ for $|d^kx/dt^k|$ which is preserved when we replace a curve $\Gamma$ by
$f_\tau(\Gamma)$. Thus, given any finite $k$, we can choose $\eta$ small
enough so that the iterated forward images of a curve $x=x(t)$ which satisfies
the inequalities
$$ |d^\ell x/dt^\ell|~\le~b_\ell\qquad{\rm for ~all}\quad \ell\le k$$
will be a curve which satisfies these same inequalities.
\smallskip

To complete the proof of Theorem \ref{th:nearhr}, we need the following.

\begin{Quote}\begin{lem}{\bf A Derivative Inequality.}\label{lem:squeeze}
 If a $C^2$-smooth function $x=x(t)$ satisfies
uniform inequalities $|x(t)|<\alpha$ and  $|d^2x/dt^2|<\beta$, then it follows
that
$$  |dx/dt|~<~\sqrt{2\alpha\beta} \,.$$
\end{lem}\end{Quote}

{\bf Proof.}
Suppose, for example, that the first derivative $x'=dx/dt$
satisfied $x'(0)\ge\sqrt{2\alpha\beta}$ with $x(0)\ge 0$. Using the lower
bound $x''(t)>-\beta$ and integrating twice, we see that
$$ x'(t)~>~\sqrt{2\alpha\beta}-\beta t\qquad{\rm and}\qquad
 x(t)~>~ \sqrt{2\alpha\beta}\;t\,-\,\beta t^2/2 $$
for $t>0$. In particular, substituting $t_0=\sqrt{2\alpha/\beta}$, it would
follow that $x(t_0)>\alpha$, thus contradicting the hypothesis.
Other cases can be handled similarly.
 \qed
        \smallskip

The proof of Theorem \ref{th:nearhr}
concludes as follows. Again let $x_n=x_n(t)$ be the
$n$-th forward image of $\Gamma_0$ under $f_\tau$. As $m$ and $n$
tend to infinity, the difference $x_n(t)-x_m(t)$ tends to zero, while the
difference $x''_n(t)-x''_m(t)$ remains uniformly bounded. Thus it follows from
Lemma \ref{lem:squeeze} that the difference $x'_n(t)-x'_m(t)$ converges
uniformly to zero. Similarly, it follows inductively that differences of higher
derivatives converge to zero.
This completes the proof that the limit curve $\Gamma_\tau$
is $C^k$-smooth for any specified $k<\infty$. Further details of
the proof are not difficult.
 \qed

\smallskip

Note that this argument does not produce a $C^\infty$ curve,
since we need to impose tighter and tighter restrictions on $f_\tau$ in
order to get successive higher derivatives. The argument certainly does not
produce a real analytic curve, which is what we would need
in order to show that $\Gamma_\tau$
is contained in a Herman ring. We have no idea how to prove real analyticity,
even assuming that the rotation number $\rho_\tau$ is Diophantine.\footnote
{It is interesting to compare the boundaries of Siegel disks for rational
maps of $\bP^1$. These are simple closed curves in all known cases. They can
never be real analytic,
but in the non-Diophantine case they can be $C^\infty$-smooth. (Compare
\cite{ABC}.) In the Diophantine case, such a boundary
necessarily contains a critical point, and hence cannot be smooth. (See \cite
{Gh}.)}

\setcounter{lem}{0}

\section{Open Problems.}\label{s:?}

 The results of this note leave a number of conjectures and open questions.
Here is a brief list.

\begin{conj}\label{cj1} For any $f$-invariant elliptic curve
$\;\cC\subset\bP^2(\C)$, the basin of attraction, consisting of all points
whose orbits converge to $\;\cC$, is contained in the Julia set of $f$.
(Intuitive proof:
Otherwise the basin would have to contain an open set $U$ such that
every sequence of iterates of $f$ on $U$ contains a subsequence converging
to a constant $c\in\cC$. This looks very unlikely, given the fact that $f$
is highly expanding in directions tangent to $\;\cC$.)
\end{conj}

The following two conjectures are closely related. Compare the discussion
in Remark \ref{rem:basin}.

\begin{conj}\label{cj2} Such an attracting basin cannot contain any
nonvacuous open set. In other words, the set of all points {\it not\/}
attracted to $\;\cC$ is always everywhere dense.
\end{conj}

\begin{conj}\label{cj3} Every invariant complex
elliptic curve contains a repelling
periodic point. (In all examples known to us there is a repelling
fixed point.)
\end{conj}

\begin{conj}\label{cj4}
In the space of complex Desboves maps with real coefficients, there is a subset
of positive measure consisting of maps which have a cycle of
attracting Herman rings. (Compare \S\ref{s:hr} and Example \ref{e4}.)
\end{conj}

There are also many questions where we have no idea what to guess.

\begin{itemize}


\item[$\bullet$] To what extent are maps with an attracting periodic orbit common in
the space of all degree $\;d\;$ maps preserving a given elliptic
curve? For example, do they form a dense open set?

\item[$\bullet$] Can an invariant elliptic curve be a global attractor, with an attracting
basin of full measure? We have constructed a number of examples where this
seems to be true empirically; but how can one exclude the possibility of
other attractors with basins of very small measure?

\item[$\bullet$] Can a smooth real elliptic curve be a trapped attractor?



\item[$\bullet$] What can one say about the dynamics when the elliptic curve has positive
transverse Lyapunov exponent? Could such a map have an absolutely continuous
invariant measure?
Is it true that an elliptic curve can never be a measure-theoretic
attractor when its transverse exponent is positive? (Compare
Remarks 1.1 and \ref{rm:Lpos}.)

\item[$\bullet$] What other kinds of attractor can occur for a rational map with invariant
elliptic curve? Can there be fractal attractors? Can there be a set of dense
orbits of positive measure, or even of full measure? Can the Julia set
have positive measure, even when it has no proper sub-attractors?

\end{itemize}

\section*{Acknowledgments}
 We want to thank the National Science Foundation
(DMS 0103646) and the
Clay Mathematics Institute for their support of Dynamical System activities
in Stony Brook.

\end{document}